\documentclass[10pt]
{article}

\usepackage{amsthm}
\usepackage[margin=30mm]{geometry}
\usepackage{amsmath}%
\usepackage{amsfonts}%
\usepackage{amssymb}%
\usepackage{amscd}
\usepackage{graphicx}
\usepackage{color}
\usepackage{xcolor}
\usepackage{bbm}
\usepackage[all,cmtip]{xy}
\usepackage{rotating}
\usepackage{url}

\usepackage{eucal}
\DeclareMathAlphabet{\mathpzc}{OT1}{pzc}{m}{it}
\usepackage{mathrsfs}

\definecolor{grau}{rgb}{0.1,0.1,0.1}
\usepackage[colorlinks=true, urlbordercolor=grau, linkcolor=grau, citecolor=grau, urlcolor=grau]{hyperref}
\usepackage{breakurl}

\usepackage{tikz}
\tikzset{node distance=2cm, auto}

\usepackage{bibspacing}

\newtheorem{theorem}{Theorem}[section]

\newtheorem{lemma}[theorem]{Lemma}
\newtheorem{corollary}[theorem]{Corollary}

\theoremstyle{definition}
\newtheorem{definition}[theorem]{Definition}
\newtheorem{example}[theorem]{Example}
\newtheorem{remark}[theorem]{Remark}
\newtheorem{remarks}[theorem]{Remarks}

\newcommand{\NN}{\mathbbm{N}} 
\newcommand{\ZZ}{\mathbbm{Z}} 
\newcommand{\Zinf}{\overline{\mathbbm{Z}}} 
\newcommand{\wo}{\backslash} 
\newcommand{\To}{\longrightarrow} 
\newcommand{\tof}[1]{\stackrel{#1}{\longrightarrow}} 
\newcommand{\fromf}[1]{\stackrel{#1}{\longleftarrow}} 
\newcommand{\im}{\textnormal{im}} 
\DeclareMathOperator*{\coker}{coker}

\newcommand{\textdef}[1]{\textnormal{\textit{#1}}}
\newcommand{\iso}{\cong} 
\newcommand{\one}{\mathbbm{1}} 

\newcommand{\pt}{\textnormal{pt}} 
\newcommand{\incl}{\hookrightarrow} 
\newcommand{\surj}{\twoheadrightarrow} 
\newcommand{\st}{\ |\ } 
\newcommand{\impl}{\Rightarrow} 
\newcommand{\Hom}{\textnormal{Hom}} 

\newcommand{\mathsmall}[1]{{\begingroup\everymath{\scriptstyle}\small{#1}\endgroup}}
\newcommand{\mathfootnotesize}[1]{{\begingroup\everymath{\scriptstyle}\footnotesize{#1}\endgroup}}

\newcommand{\id}{\textnormal{id}} 

\newcommand{\wt}[1]{\widetilde{#1}}

\newcommand{\ext}{\textnormal{Ext}}
\newcommand{\tor}{\textnormal{Tor}}
\newcommand{\blank}{\underline{\ \ }}
\newcommand{\topsmash}{\wedge}

\newcommand{\homotopic}{\sim}
\newcommand{\GL}{\textnormal{GL}}

\newcommand{\chain}{\textnormal{Ch}}
\newcommand{\hchain}{\textnormal{hCh}}
\newcommand{\complex}[2]{\chain^{#1}(#2)}
\newcommand{\hcomplex}[2]{\hchain^{#1}(#2)}

\newcommand{\tot}{\textnormal{Tot}}
\newcommand{\resolution}{\textnormal{Res}}

\newcommand{\catfontCapital}{\emph}
\newcommand{\catfontSmall}{\emph}

\newcommand{\projClass}{\mathcal}
\newcommand{\injClass}{\projClass}
\newcommand{\cat}{\mathcal}

\newcommand{\sheaves}{\catfontCapital{S}\catfontSmall{h}}
\newcommand{\topcat}{\catfontCapital{T}\catfontSmall{op}}
\newcommand{\spectracat}{\catfontCapital{S}\catfontSmall{p}}
\newcommand{\abcat}{\catfontCapital{A}\catfontSmall{b}}
\newcommand{\grabcat}{\catfontCapital{G}\catfontSmall{r}\catfontCapital{A}\catfontSmall{b}}

\DeclareMathOperator*{\down}{down}

\newcommand{\leqlex}{\leq_{\textnormal{lex}}} 
\newcommand{\lesslex}{<_{\textnormal{lex}}} 
\newcommand{\geqlex}{\geq_{\textnormal{lex}}} 

\begin{document}

\title {Successive Spectral Sequences}

\author{%
Benjamin Matschke%
\setcounter{footnote}{-1}%
}
\date{}

\maketitle

\begin{abstract}
In this paper, we develop a structure theory for generalized spectral sequences, which are derived from chain complexes that are filtered over arbitrary partially ordered sets.
Also, a more general construction method reminiscent of exact couples is studied, together with examples where they arise naturally.
As for ordinary spectral sequences we will see differentials and group extensions, however the real power comes from the appearance of natural isomorphims between pages of differing indices.

The constructions reveal finer invariants than ordinary spectral sequences, and they connect to other fields such as Fary functors and perverse sheaves. They are based on a natural index scheme, which allows us to obtain new results even in the standard case of $\ZZ$-filtered chain complexes, e.g.\ a useful criterion for a product structure for Grothendieck's spectral sequences, and 
news paths to connect the first or second page to the limit.

This turns out to yield the right framework for unifying several spectral sequences that one would usually apply one after another.
Examples that we work out are successive Leray--Serre spectral sequences, the Adams--Novikov spectral sequence following the chromatic spectral sequence, successive Grothen\-dieck spectral sequences, and successive Eilenberg--Moore spectral sequences.

\end{abstract}

%
%
%
%

%
%
%
%


\makeatletter
\renewcommand*\l@section[2]{%
  \ifnum \c@tocdepth >\z@
    \addpenalty\@secpenalty
    \addvspace{0.5em \@plus\p@}%
    \setlength\@tempdima{1.5em}%
    \begingroup
      \parindent \z@ \rightskip \@pnumwidth
      \parfillskip -\@pnumwidth
      \leavevmode 
      \advance\leftskip\@tempdima
      \hskip -\leftskip
      {\normalsize 
      #1}\nobreak\hfil 
      \nobreak\hb@xt@\@pnumwidth{\hss #2}\par\vspace{0.5em}
    \endgroup
  \fi}
\makeatother


\renewcommand\contentsname{Contents}
\setcounter{tocdepth}{2}
\noindent

{\small \tableofcontents}

\section{Introduction}

A spectral sequence is a computational tool in algebraic topology, which relates the homology of a chain complex (or of a space, or a spectrum) to approximations thereof that are constructed from a $\ZZ$-filtration of the chain complex.
In this paper, we more generally consider filtrations over arbitrary posets.
Restricting such a filtration to any chain of the poset yields an associated ordinary spectral sequence.
We study how the pages of these spectral sequences can be related to form one \emph{spectral system} (short for system of spectral sequences).
This includes basic maps induced by inclusion, induced differentials and kernels and cokernels thereof, group extensions, but also natural isomorphisms, splitting principles, and multiplicative structures.
The natural isomorphisms are the first non-trivial property of spectral systems that do not appear in ordinary spectral sequences.
In the frequent case of chain complexes (spaces, spectra) with several $\ZZ$-filtrations, this provided various ways to connect the basic first and second pages of the spectral system to its limit.
\footnote{In a follow-up paper~\cite{Mat14HigherSpectralSequences}, we construct additional differentials for spectral systems that arise from several $\ZZ$-filtrations, which motivates the name ``higher spectral sequence'' for this important special case.}

\paragraph{A motivating example: Successive Leray--Serre spectral sequences.}
When one applies spectral sequences successively 
it fairly often happens that there are several ways to do this.
As a main application of this paper, we will give several general examples of how to unify these ways into one spectral system.

As an example,
consider the following tower of two fibrations,
\begin{equation}
\label{eqVerticalSeqOfFibrations}
\xymatrix{
F \ar[r]^{i_{FE}} & E\ar[d]^{p_{EN}} \\
M \ar[r]^{i_{MN}} & N\ar[d]^{p_{NB}} \\
                  & B.
}
\end{equation}
Our aim is to compute $H_*(E)$ from $H_*(F)$, $H_*(M)$, and $H_*(B)$. We could apply a Leray--Serre spectral sequence to compute $H_*(N)$ from $H_*(M)$ and $H_*(B)$ as an intermediate step, and then apply another Leray--Serre spectral sequence to derive $H_*(E)$ from $H_*(F)$ and $H_*(N)$.
Alternatively we might define $P:=p_{EN}^{-1}(i_{MN}(M))=(p_{NB}\circ p_{EN})^{-1}(\pt)$ and get another sequence of fibrations,
\begin{equation}
\label{eqHorizontalSeqOfFibrations}
\xymatrix{
F \ar[r]^{i_{FP}} & P\ar[r]^{i_{PE}}\ar[d]^{p_{PM}} & E\ar[d]^{p_{EB}} \\
                  & M                               & B. \\
}
\end{equation}
Here we can apply two different but related Leray--Serre spectral sequences with intermediate step $H_*(P)$.
The following diagram illustrates this `associativity law',
\[
\xymatrix{
H_*(F), H_*(M), H_*(B) \ar[rr]^{\ \ \id\times\textnormal{LS}} \ar[d]^{\textnormal{LS}\times\id} & & H_*(F),H_*(N)\ar[d]^{\textnormal{LS}} \\
H_*(P),H_*(B) \ar[rr]^{\ \ \textnormal{LS}}                                                     & & H_*(E). \\
}
\]
Under certain conditions on the fibrations we construct in Section~\ref{secSuccFibrations} a spectral system with ``second page'' $H_*(B,H_*(M,H_*(F)))$ that converges to $H_*(E)$.
It contains considerably finer information than the four separate spectral sequences.


\medskip

\paragraph{Outline.}
In Section~\ref{secGeneralizedFiltration} we construct a spectral system for chain complexes that are filtered over an arbitrary poset. 

In Section~\ref{secSeveralZfiltrations} we treat the special case of chain complexes that are $\ZZ$-filtered in several ways and prove basic properties, some of which have only trivial analogs in usual spectral sequences.
We show that there are several interesting non-trivial connections between the $0$-page and the usual goal of computation $H_*(C)$.

Section~\ref{secExactCoupleAnalogs} is about exact couple systems, which lead to spectral systems that are more general than the ones coming from $I$-filtered chain complexes, in a similar fashion as how exact couples generalize $\ZZ$-filtered chain complexes.
We study basic properties such as differentials, extensions, natural isomorphisms, splitting principles, and multiplicative structures.

We give the following examples of spectral systems:
\begin{enumerate}
\item A basic running example in this paper is the spectral system for the generalized homology theory of a $I$-filtered space $X$.
It has many desirable pro\-per\-ties, and the next two examples are instances of it.
\item As mentioned above there is a spectral system for iterated fibrations, generalizing the Leray--Serre and Atiyah--Hirzebruch spectral sequences, see Section~\ref{secSuccFibrations}.
\item In Section~\ref{secSuccGrothendieckSS} we construct a generalization of Grothendieck's spectral sequence in the setting of a composition of several functors.
We also provide a general and natural condition for the existence of a product structure, which seems new also for the ordinary Grothendieck spectral sequence. 
\item The $E_2$-page of the Adams--Novikov spectral sequence is the limit of the chromatic spectral sequence; we show how to put this into a single spectral system in Section~\ref{secAdamsAndChromaticSS}.
\item In Section~\ref{secSuccEilenbergMoore} we construct a spectral system for a cube of fibrations, where one would usually apply Eilenberg--Moore spectral sequences successively.
\end{enumerate}

Background on spectral sequences can be found in
McCleary~\cite{McC01userGuideToSS},
Spanier \cite{Spanier66algTop},
Cartan--Eilenberg~\cite{CarEil56homologicalAlgebra},
Hatcher~\cite{Hatcher04SS},
Weibel~\cite{Wei94homologicalAlgebra},
Gelfand--Manin \cite{GelfandManin03methHomAlg},
Bott--Tu~\cite{BottTu82diffFormsInAlgTop}, 
Switzer~\cite{Switzer75algTopHomotopyAndHomology}, and many more.
A~detailed account on the early history of spectral sequences is given in McCleary \cite{McC99historyOfSS}.

\paragraph{Previous generalizations of spectral sequences.}

Several very useful extensions of spectral sequences are known: 
Po Hu's~\cite{Hu98thesis,Hu99transfiniteSS} \emph{transfinite spectral sequences} are the closest among them to this paper.
They consists of terms $E^r_{p,q}$ as in the standard setting, except that $p,q,r$ are elements in the Grothendieck group $G(\omega^\alpha)$, which in case $\alpha=n\in\NN$ is isomorphic to $\ZZ^n$ with the lexicographic order.
Thus, indices $P\in G(\omega^\alpha)=\ZZ^n$ correspond to the sets $A(p_1,\ldots,p_n;\id_{\ZZ^n})$ in Section~\ref{secLexicographicConnection}, which makes transfinite spectral sequences basically a substructure of spectral systems, using the lexicographic connection from Lemma~\ref{lemLexConnectionInSmallSteps} (that transfinite spectral sequences need to have $G(\omega^\alpha)$-graded pages as in~\cite{Hu99transfiniteSS} is actually not necessary).
On the other hand, that spectral systems have a richer structure means that not all transfinite spectral sequences will naturally generalize to spectral systems.
Hu gave several examples of transfinite spectral sequences.

Deligne \cite{Del71hodge2} studied in connection with his mixed Hodge structures spectral sequences of chain complexes with respect to several $\ZZ$ -filtrations. 
Compare \cite[\S 1]{Del71hodge2} with Section~\ref{secSeveralZfiltrations} of this paper.

Behrens~\cite{Beh10goodwillieTowerAndEHP} constructed transfinite versions of the Atiyah--Hirzebruch, EHP, Goodwillie, and homotopy spectral sequences; see also Behrens~\cite{Beh06rootInvariantAndAdamsSS}.

The exact couple systems of Section~\ref{secExactCoupleAnalogs} are canonical substructures of perverse sheaves (Be{\u\i}linson--Bernstein--Deligne~\cite{BBD82perverseSheaves}) and more generally of F\'ary functors (Grinberg--MacPherson~\cite{GriMac99EulerCharactAndLagrangIntersect}).
Thus, the latter give also rise to spectral systems over the poset of open sets of the underlying space.

\section{The spectral system of a generalized filtration} \label{secGeneralizedFiltration}

Throughout Sections~\ref{secGeneralizedFiltration}, \ref{secSeveralZfiltrations}, and \ref{secExactCoupleAnalogs}, we will simplify notation by omitting the grading of chain complexes and homology; compare with Remark~\ref{remGrading}.

Let $(C,d)$ be a chain complex, which has several subcomplexes $F_i$, $i\in I$.
Then the family $F:=\{F_i\st i\in I\}$ can be seen as a generalized filtration, since we do not require any inclusion relations in~$F$.
Let us give $I$ the structure of a poset by $i\leq j$ for $i,j\in I$ if and only if $F_i\subseteq F_j$.
We say that $C$ is filtered over the poset $I$, or $I$-filtered for short.

Every chain of countable size in $I$ gives rise to a spectral sequence that converges to $H(C)$ (under some standard assumptions).
The questions we are interested in are: How they are related, and is there a conceptual larger device that contains all those spectral sequences?

Interesting generalized filtrations arise when $C$ is filtered over $\ZZ$ in two or more different ways, see Section~\ref{secSeveralZfiltrations}.

Another basic example to have in mind comes from an $I$-filtered space $X$, that is, from a family of subspaces $(X_i)_{i\in I}$ of some space~$X$, ordered by inclusion.
The singular chain complex $C_*(X)$ of $X$ is then filtered by $F_i:=\im (C_*(X_i)\incl C_*(X))$, $i\in I$.
Similarly, the singular cochain complex $C^*(X)$ of X is filtered by $F^i:=\ker (C^*(X)\to C^*(X_i))$, $i\in I^*$, $I^*$ being the dual poset of $I$ (same elements, reversed relations); compare with Remarks~\ref{remSurjections} and~\ref{remCohomology}.

\subsection{Construction and basic properties} \label{secMainConstruction}

For two subgroups $X$ and $Y$ of a larger abelian group, we write
\[
X/Y := X/(X\cap Y) \iso (X+Y)/Y,
\]
which will simplify notation considerably.

For $p,q,z,b\in I$, we define
\begin{equation}
\label{eqSterms}
S^{pz}_{bq} := \frac{F_p\cap d^{-1}(F_z)}{d(F_b)+F_q}.
\end{equation}
We call $p$ as usual the \textdef{filtration degree}, $z$ the \textdef{cycle degree}, $b$ the \textdef{boundary degree}, and $q$ the \textdef{quotient degree} (note that $q$ has here a different meaning as in the standard notation $E_{p,q}^r$).
In what follows we will only consider the $S$-terms for which $F_q\subseteq F_p$.
If $F$ is closed under taking intersections and sums, then we may assume that
\begin{equation}
\label{eqRelevantSterms}
F_z\subseteq F_p\supseteq F_q\subseteq F_b,
\end{equation}
since we can intersect $F_z$ by $F_p$ and add $F_q$ to $F_b$ without changing the $S$-term.

If the inclusions
\[
F_z\subseteq F_q \textnormal{ and } F_p\subseteq F_b
\]
are fulfilled, then $S^{pz}_{bq}$ can be written as
\begin{align}
\label{eqSintermsofFirstPage}
S^{pz}_{bq} & =\ \im\big(H\left(F_p/F_z\right) \tof{\textnormal{incl}_*} H\left(F_b/F_q\right) \big) \\
& =\ \frac{\ker\big(H(F_p/F_q)\to H(F_q/F_z)\big)}{\im\big(H(F_b/F_p)\to H(F_p/F_q)\big)}.
\end{align}

We may assume that $F$ contains $F_{-\infty}:=0$ (the zero complex) and $F_\infty:=C$.
Then our goal of computation is $H(C)=S^{\infty,-\infty}_{\infty,-\infty}$%
, which we call the \emph{limit} of $S$.
(Note that for us the limit is part of the structure of $S$, it does not imply any kind of convergence in the usual sense.)

As with ordinary spectral sequences of a $\ZZ$-filtration, there are the following ways to relate different $S$-terms.

\subsection{Extension property}
If $F_z\subseteq F_{p_0}\subseteq \ldots\subseteq F_{p_n}\subseteq F_b$, then $S^{p_n,z}_{b,p_0}$ is a successive group extension of the $S^{p_i,z}_{b,p_{i-1}}$, $i=1\ldots n$, since the sequence
\[
0\to S^{p_{i},z}_{b,p_{i-1}}\to S^{p_{i+1},z}_{b,p_{i-1}}\to S^{p_{i+1},z}_{b,p_{i}}\to 0
\]
is exact.

\subsection{Differential}
Suppose that two quadruples $(p_1,z_1,b_1,q_1)$ and $(p_2,z_2,b_2,q_2)$ in $I^4$ satisfy
\begin{equation}
\label{eqDifferentialCondition}
F_{z_2}\subseteq F_{p_1}\textnormal{ and } F_{q_2}\subseteq F_{b_1}.
\end{equation}
Then $d$ induces a well-defined differential
\[
d: S^{p_2z_2}_{b_2q_2} \To S^{p_1z_1}_{b_1q_1},
\]
which we also denote by $d$.%
\footnote{The condition \eqref{eqDifferentialCondition} is not most general.
In fact we only need to assume that $d(F_{p_2})\cap F_{z_2} \subseteq F_{p_1}$ and $d(F_{q_2})\subseteq d(F_{b_1})+F_{q_1}$.}

Now suppose that we have a sequence of such differentials,
\[
S^{p_3z_3}_{b_3q_3} \tof{d_3} S^{p_2z_2}_{b_2q_2} \tof{d_2} S^{p_1z_1}_{b_1q_1},
\]
such that \eqref{eqDifferentialCondition} and the corresponding inclusions for $d_3$ are fulfilled.
It might help to visualize these inclusions as follows:
\newcommand{\subseteqRot}{\stackrel{}{\begin{sideways}$\subseteq$\end{sideways}}}
\begin{equation}
\label{eqInclusionDiagram}
\begin{matrix}
F_{z_1} & \ & F_{q_1} & \subseteq & F_{p_1} & \ & F_{b_1} &&&& &&&&\\
&&&& \subseteqRot & & \subseteqRot &&&& &&&&\\
&&&& F_{z_2} & \ & F_{q_2} & \subseteq & F_{p_2} & \ & F_{b_2} &&&&\\
&&&& &&&& \subseteqRot & & \subseteqRot &&&&\\
&&&& &&&& F_{z_3} & \ & F_{q_3} & \subseteq & F_{p_3} & \ & F_{b_3}\\
\end{matrix}
\end{equation}
If
\begin{equation}
\label{eqKernelCondition}
q_2=b_1 \textnormal{ and } F_{z_2}\cap F_{q_1}=F_{z^*} \textnormal{ for some } z^*\in I,%
\footnote{For example if $z_2=p_1$ we can take $z^*=q_1$.}
\end{equation}
then the kernel of $d_2$ is given by
\[
\ker(d_2) = S^{p_2,z^*}_{b_2,q_2}.
\]
If
\begin{equation}
\label{eqCokernelCondition}
z_3=p_2 \textnormal{ and } F_{b_2}+F_{p_3}=F_{b^*} \textnormal{ for some } b^*\in I,
\footnote{For example if $b_2=q_3$ we can take $b^*=p_3$.}
\end{equation}
then the cokernel of $d_3$ is given by
\[
\coker(d_3) = S^{p_2,z_2}_{b^*,q_2}.
\]
If both conditions \eqref{eqKernelCondition} and \eqref{eqCokernelCondition} hold then we can compute the homology
\begin{equation}
\label{eqHomologyOfEd}
\frac{\ker(d_2)}{\im(d_3)} = S^{p_2,z^*}_{b^*,q_2}.
\end{equation}
Interesting special case: if all four vertical inclusions in Diagram \eqref{eqInclusionDiagram} hold with equality then
\begin{equation}
\label{eqHomologyOfE}
\frac{\ker(d_2)}{\im(d_3)} = S^{p_2,q_1}_{p_3,q_2}.
\end{equation}

\begin{example}[$0$- and $1$-page] 
We call the collection of all $E$-terms of the form
\[
S^{pp}_{qq} = F_p/F_q
\]
the \emph{$0$-page}.
Here, diagram \eqref{eqInclusionDiagram} has the pattern
\[
\begin{matrix}
p & q & p & q &   &   &   &  \\
  &   & p & q & p & q &   &  \\
  &   &   &   & p & q & p & q\\
\end{matrix}
\]
The induced differential $S^{pp}_{qq}\to S^{pp}_{qq}\to S^{pp}_{qq}$ coincides with the boundary map that $d$ induces on $F_p/F_q$.
By \eqref{eqHomologyOfE} the homology of the middle term is
\[
S^{pq}_{pq} = H(F_p/F_q).
\]
We call the collection of all these $S$-terms the \emph{$1$-page}.
By \eqref{eqSintermsofFirstPage}, the $1$-page determines all other relevant $S$-terms as long as we know the maps between them that are induced by inclusion $F_p\incl F_b$.
However, in many applications, only the $S$-terms of the $0$- and $1$-page where $p$ covers $q$ are given. 
\end{example}

\subsection{Remarks}

\begin{remark}[Isomorphic $S$-terms]
For general index posets $I$ it may happen that two $S$-terms with different indices are naturally isomorphic (natural with respect to filtration preserving maps between $I$-filtered chain complexes).
This does occur not in the classic case $I=\ZZ$, but it does when $I=\ZZ^n$, $n\geq 2$, where it turns out the be very useful.
We postpone this connection to Sections~\ref{secReducedSterms} and~\ref{secExactCouplesNatIsos}.
\end{remark}


\begin{remark}[Limits, convergence, comparison]
As in ordinary spectral sequences, depending on the given filtration we might need to take differentials or do extensions an infinite number of times in order to connect the $0$-page to $H(C)$.
The convergence and comparison theorems from Eilenberg--Moore~\cite{EilMoo62limitsAndSS}, Boardman~\cite{Boardman99conditionallyConvergentSS}, McCleary~\cite{McC01userGuideToSS}, and Weibel~\cite{Wei94homologicalAlgebra} go over to the general setting without difficulties.
For simplicity we will ignore this, for example by assuming that the given filtration is of \emph{finite height}, that is, only finitely many different terms $F_p$ appear in any chain of the filtration.
If $C$ is graded then we will require this only degree-wise.
\end{remark}

\begin{remark}[Surjections]
\label{remSurjections}
A generalized filtration of a chain complex is a commutative diagram of chain complexes where all maps are injections.
We could do all constructions of this section dually, starting with a diagram of surjections instead of injections.
This leads essentially to the same spectral system, since one can filter $C$ by the kernels of the given surjections from $C$ to all the chain complexes in the diagram and take the spectral system of this filtration.
\end{remark}

\begin{remark}[Cohomology]
\label{remCohomology}
%
Most of the paper is written in terms of homology.
If the reader prefers cohomology, all one has to do is dualizing the underlying poset $I$, which makes increasing filtrations in this paper decreasing.
For a particular instance, see Example~\ref{exExactCoupleOfFilteredSpace}.
All maps stay the same, except for morphisms between spectral systems, which are now contravariant functors.
\end{remark}

\begin{remark}[Grading]
\label{remGrading}
Instead of omitting the grading of a chain complex in the notation, we can work with a more general ungraded chain complex, that is, an abelian group $C$ (or an object in some other abelian category) together with an endomorphism $d:C\to C$ with $d^2=0$.
All constructions in this and the next two sections work equally in the graded and in the ungraded case.
For example if $C$ has the structure of a graded abelian group (graded over $\ZZ$ or any other abelian group), $d$ is a graded homomorphism of a fixed degree, and the $F_i$ are graded subgroups, then the $S$-terms are graded as well and the differential is graded with the same degree as $d$.
\end{remark}


\subsection{The spectral sequence of a $\ZZ$-filtration}
\label{secSSofZfiltration}

Suppose that the filtration $F$ is indexed over the integers, $F=(F_p)_{p\in\ZZ}$.
Then the well-known \textdef{spectral sequence of the $\ZZ$-filtration} $F$ is defined as
\[
E_{p,q}^r := (S^{p,p-r}_{p+r-1,p-1})_{(p+q)},
\]
together with differentials
\[
d: E_{p,q}^r \to E_{p-r,q+r-1}^r.
\]
In this grading, $d$ is a differential of bidegree $(-r,r-1)$ in $E_{*,*}^r:=\oplus_{p,q} E_{p,q}^r$.
Equation \eqref{eqHomologyOfEd} implies that
\[
H(E_{*,*}^r) = E_{*,*}^{r+1}.
\]

The $E^\infty$-term is
\[
E_{p,q}^\infty = (S^{p,-\infty}_{\infty,p-1})_{(p+q)},
\]
hence the above extension property recovers (or at least relates to) $H_*(C)=S^{\infty,-\infty}_{\infty,-\infty}$ by iterated group extensions.

\begin{remark}
Note that the spectral sequence of a $\ZZ$-filtration contains only a small subset of the $S^{pz}_{bq}$-terms, namely those where $p-q=1$ and $b-q=p-z=r$.
This is definitely a good choice, since it connects the $E^0$-page to $H_*(C)$ by a sequence of differentials and extensions in a direct fashion, which contains often all necessary information.
Often there is also additional structure such as products.
Sometimes it might be reasonable to look at other $S$-terms as well, since for a $\ZZ$-filtration there are many connecting differentials and extensions.%
\footnote{A particular example, where this might be useful, are spectral sequence valued indices of $G$-bundles, see~\cite{Mat12paramGromovWaist}.}

For example it can be reasonable to look at terms
\[
E_{p,q}^{r_b,r_z} := (S^{p,p-r_z}_{p+r_b-1,p-1})_{(p+q)},
\]
that is, the pages are now indexed over all $(r_b,r_z)\in\ZZ_{\geq 0}^2$.
For $r_b=r_z=:r$ this coincides with $E_{p,q}^r$.
There is a connecting differential $d:E_{p_1,q_1}^{r_{b1},r_{z1}}\to E_{p_2,q_2}^{r_{b2},r_{z2}}$ if $p_1+q_1=p_2+q_2+1$, $r_{z1}=r_{b2}$ and $p_1=p_2+r_{z1}$, and we have
\[
\ker d = E_{p_1,q_1}^{r_{b1},r_{z1}+1}
\] and
\[
\coker d = E_{p_2,q_2}^{r_{b2}+1,r_{z2}}.
\]
One can also mix differentials and extensions along the way (but then the $S$-notation seems more convenient).
\end{remark}

%
%
%
%
%
%
%

\section{Several $\ZZ$-filtrations} \label{secSeveralZfiltrations}

Suppose our given chain complex $C$ is $\ZZ$-filtered in $n$ different ways,
\begin{equation}
\label{eqNZfiltrations}
\ldots \subseteq F^{(i)}_{k-1} \subseteq F^{(i)}_k \subseteq F^{(i)}_{k+1} \subseteq \ldots,\ \ \ 1\leq i\leq n.
\end{equation}
In this section we construct a spectral system that incorporates all the spectral sequences of the $n$ filtrations.

We denote $\Zinf:=\ZZ\cup\{-\infty,\infty\}$, which is a lattice ordered by `$\leq$'.
$\Zinf^n$ is the product poset, ordered by $P=(p_1,\ldots,p_n)\geq Q=(q_1,\ldots,q_n)$ if and only $p_i\geq q_i$ for all $i$.
Let $D(\Zinf^n)$ denote the set of downsets $p\subseteq\Zinf^n$, that is, all subsets $p$ that satisfy: if $Q\leq P$ and $P\in p$ then $Q\in p$.
$D(\Zinf^n)$ will be the underlying poset for our spectral system.

We extend the $n$ $\ZZ$-filtrations \eqref{eqNZfiltrations} by
setting $F^{(i)}_{-\infty}:=\bigcap_{k\in\ZZ}F^{(i)}_k$ and $F^{(i)}_{\infty}:=C$.%
\footnote{\label{footDefOfFpmInfty}%
The reader should not wonder why we did not set $F^{(i)}_{-\infty}:=0$ or alternatively $F^{(i)}_{\infty}:=\sum_{k\in\ZZ}F^{(i)}_k$.
The reason for this definition is that it makes \eqref{eqFforDownsetsInZn} very natural. E.g. for $n=1$, we have $F_{\emptyset} = 0$, $F_{\{-\infty\}}=\bigcap_{k\in\ZZ}F^{(1)}_k$, $F_{\ZZ\cup\{-\infty\}}=\sum_{k\in\ZZ}F^{(1)}_k$, and $F_{\Zinf} = C$.
}
For any $P\in \Zinf^n$ we define
\[
F_{P} := F_{p_1}^{(1)}\cap\ldots\cap F_{p_n}^{(n)}.
\]
For any downset $p\in D(\Zinf^n)$ we define
\begin{equation}
\label{eqFforDownsetsInZn}
F_p := \sum_{P\in p} F_{P} \subseteq C,
\end{equation}
where the sum denotes taking interior sums of subcomplexes of~$C$.
Thus, $(F_p)_{p\in D(\Zinf^n)}$ is a $D(\Zinf^n)$-filtration of~$C$.

\begin{definition}
Let $I$ be a distributive lattice with meet denoted `$\cap$' and join denoted~`$\cup$'.
We call an $I$-filtration $(G_p)_{p\in I}$ \emph{distributive} if $G_{a\cap b}=G_a\cap G_b$ and $G_{a\cup b}=G_a + G_b$ for all $a,b\in I$.
\end{definition}

\begin{example}
For distributivity of the above filtration $(F_p)_{p\in D(\Zinf^n)}$ only $F_{a\cap b}=F_a\cap F_b$ needs to be checked.
Dually we could have defined a $D(\Zinf^n)$-filtration $F^p := \bigcap_{P\not\in p-\one} F_{p_1}^{(1)}+\ldots +F_{p_n}^{(n)}$, which in turn is distributive if and only if $F^{a\cup b}=F^a + F^b$ for all $a,b$.
If one of $(F_p)$ and $(F^p)$ is distributive and the underlying $\ZZ$-filtrations~\eqref{eqNZfiltrations} are finite then $(F_p)=(F^p)$.
(This finiteness assumption can be omitted if the filtrations are completely distributive.)
\end{example}

\begin{example}
If $C$ is $\ZZ^n$-graded and \eqref{eqNZfiltrations} are the $n$ canonical filtrations of $C$ then $(F_p)$ is clearly distributive.
\end{example}

\begin{remark}
We could have started this section right away from an arbitrary $D(\Zinf^n)$-filtration $(F_p)$.
All statements below will only require it to be distributive.
(In practice, $\Zinf^n$-filtrations usually arise as the common refinement of $n$ $\ZZ$-filtrations.)
\end{remark}

\subsection{Naturally isomorphic $S$-terms} \label{secReducedSterms}

Let us regard $\Zinf^n$ as the vertices of a graph $G$, where two elements $P,Q\in \Zinf^n$ are connected by an edge if and only if 
$P\leq Q$ or $P\geq Q$.
In this sense, we call a subset $A\subseteq\Zinf^n$ \emph{connected} if and only if the induced subgraph $G[A]$ is connected.

\begin{definition}
To any downsets $z\leq q\leq p\leq b$ in $D(\Zinf^n)$, let
\[
Z(z,q,p,b)\subseteq \Zinf^n
\]
denote the union of all connected components of $p\wo z$ that intersect $p\wo q$, and let
\[
B(z,q,p,b)\subseteq \Zinf^n
\]
denote the union of all connected components of $b\wo q$ that intersect $p\wo q$.

We call $S^{pz}_{bq}$ with $z\leq q\leq p\leq b$ a \emph{reduced $S$-term} if $p\wo z = Z(z,q,p,b)$ and $b\wo q = B(z,q,p,b)$.
\end{definition}


The following lemma shows that any $S^{pz}_{bq}$ with $z\leq q\leq p\leq b$ has an isomorphic reduced $S$-term $S^{p\wt z}_{\wt bq}$ (canonical with respect to maps between $\Zinf^n$-filtered chain complexes).

\begin{lemma}[Reducing $S^{pz}_{bq}$]
\label{lemReducingSterms}
Let $z,q,p,b\in D(\Zinf^n)$ be downsets with $z\leq q\leq p\leq b$.
Let $Z:=Z(z,q,p,b)$ and $B:=B(z,q,p,b)$.
Define $\wt z:=p\wo Z$ and $\wt b:=q\cup B$.
Then the inclusions $z\incl\wt z$ and $\wt b\incl b$ induce a commutative diagram of four natural isomorphisms of $S$-terms,
\begin{equation}
\label{eqReducingSterm}
\xymatrix{
S^{pz}_{\wt bq} \ar[r]^\iso \ar[d]_\iso & S^{p\wt z}_{\wt bq} \ar[d]^\iso\\
S^{pz}_{bq} \ar[r]^\iso                 & S^{p\wt z}_{bq},
}
\end{equation}
and $S^{p\wt z}_{\wt bq}$ is reduced.
\end{lemma}

Later 
we will prove a similar statement for exact couple systems, Lemma~\ref{lemNaturalIsomorphisms}, which has a more conceptual proof, however it does not imply the full Lemma~\ref{lemReducingSterms}.

\begin{proof}
First we prove that the left map in \eqref{eqReducingSterm} is an isomorphism.
Surjectivity is trivial, since on both sides elements are represented by the same group $F_p\cap  d^{-1}(F_z)$.
In order to prove injectivity, let $x\in F_p\cap  d^{-1}(F_z)$ represent an element that goes to zero in $S^{pz}_{bq}$, that is, $x\in F_p\cap  d^{-1}(F_z) \cap (F_q + d(F_b))$.
Note that $F_b=F_{\wt b}+F_C$, where $C$ is the downset of $b\wo \wt b$. 
Thus $x$ is of the form $x_q+d(x_{\wt b})+d(x_C)$ where $x_q\in F_q$, $x_{\wt b}\in F_{\wt b}$, and $x_C\in F_C$.
From $q\subseteq p$ and $x\in F_p$ it follows that $d(x_{\wt b})+d(x_C)\in F_p$.
The key point is that $b\wo \wt b$ and $B$ are disconnected, which implies that $d(x_C)\in F_q$.
Thus $x\in F_q + d(F_{\wt{b}})$.
This means that $x$ represents zero in $S^{pz}_{\wt bq}$.
Analogously, the right map 
is an isomorphism.

For the bottom map in~\eqref{eqReducingSterm}:
Injectivity is immediate, since $z\subseteq\wt z$.
In order to prove surjectivity, let $x\in F_p\cap d^{-1}(F_{\wt z})$ represent an element in $S^{p\wt z}_{bq}$.
Note that $F_p = F_z + F_Y + F_{Z'}$, where $Y$ is the downset of $\wt z\wo z$ and $Z'$ is the downset of $Z$.
Thus $x$ is of the form $x_z + x_Y + x_{Z'}$.
Since $z\leq q$ and $d(x_z)\in F_z$, $x$ and $x_Y+x_{Z'}$ represent the same element in $S^{p\wt z}_{bq}$.
The key is again that $\wt z\wo z$ and $Z$ are disconnected.
Thus, since $d(x_Y)+d(x_{Z'})\in F_{\wt z}$, we have $d(x_{Z'})\in F_z$.
Since also $x_Y\in F_q$, $x$ and $x_{Z'}$ represent the same element in $S^{p\wt z}_{bq}$.
Since $x_{Z'}\in F_p\cap d^{-1}(F_z)$, $x_{Z'}$ already represents an element in $S^{pz}_{bq}$.
This shows surjectivity.
Analogously, the top map 
is an isomorphism.


By definition, $S^{p\wt z}_{\wt bq}$ is reduced.
\end{proof}

\begin{remark}[Other canonical forms for $S$-terms]
Given $z\leq q\leq p\leq b$, the lemma replaces $z$ and $b$ by the largest $\wt z\leq q$ and smallest $\wt b\geq p$ such that $S^{p\wt z}_{\wt bq}$ is isomorphic to $S^{pz}_{bq}$.
Instead we could take also the \emph{smallest} $\wt z\leq z$ and the \emph{largest} $\wt b\geq b$, obtaining a different normal form, which might be useful in some specific cases.
\end{remark}


\begin{lemma}[Natural isomorphisms]
\label{lemNaturallyIsomorphicSterms}
Suppose $(F_p)$ is a distributive $D(\Zinf^n)$-filtration of~$C$.
Let $z,q,p,b\in D(\Zinf^n)$ be downsets with $z\leq q\leq p\leq b$.
Then $Z:=Z(z,q,p,b)$ and $B:=B(z,q,p,b)$ determine $S^{pz}_{bq}$ up to natural isomorphism, that is, all $z\leq q\leq p\leq b$ with the same $Z$ and $B$ give naturally isomorphic $S$-terms.
\end{lemma}

\begin{proof}
This can be proved similarly as the previous lemma.
More elegantly, we can first reduce to the case where $z=p\wo Z$, $z'=p'\wo Z$, $b=q\cup B$, and $b'=q'\cup B$ using the previous lemma.
Then the assertion follows from Lemma~\ref{lemNaturalIsomorphisms} (the associated exact couple system is excisive). 
%
%
\end{proof}

Thus when we are only interested in the isomorphism class we may write $S^{pz}_{bq}$ simply as $S^Z_B$.
When differentials are involved we usually prefer the notation $S^{pz}_{bq}$.

%

\begin{remark}[Splitting $S$-terms]
Suppose $z\leq q\leq p\leq b$, and let $p\wo q = X_1\dot\cup X_2$.
Let $Z_1$ be the union of connected components of $p\wo z$ that intersect $X_1$.
Let $B_1$ be the union of connected components of $b\wo q$ that intersect $X_2$.
Analogously, define $Z_2$ and $B_2$.
Assume further that $Z_1\cap Z_2=\emptyset$ and $B_1\cap B_2=\emptyset$.
Then there is a natural isomorphism
\[
S^{pz}_{bq}\iso S^{Z_1}_{B_1}\oplus S^{Z_2}_{B_2}.
\]
\end{remark}

\begin{remark}
Lemma~\ref{lemNaturallyIsomorphicSterms} also holds more generally under the weaker condition $b\geq q\leq p\geq z$, which includes terms $S^{pp}_{qq}$ of the $0$-page, in the following way: $S^{pz}_{bq}$ is determined up to natural isomorphism by $Z$, $B$ and $p\wo q$.
Note that $p\wo q = Z\cap B$ does not need to hold anymore.
\end{remark}

\subsection{Connections} \label{secConnections}

Consider a distributive $D(\Zinf^n)$-filtration $(F_p)$ of~$C$.
A priori it is not be clear what is the best way to connect the $S$-terms $S^{pp}_{qq}=F_p/F_q$, $p \succ q$, from the $1$-page to the goal of computation $S^{\infty,-\infty}_{\infty,-\infty}=H_*(C)$ by a collection of kernels, cokernels, extensions, limits, and natural isomorphisms.
The answer depends on the particular situation and there are several choices.

\subsubsection{Lexicographic connections} \label{secLexicographicConnection}
Based on the well-known spectral sequence of a $\ZZ$-filtration, the most apparent way to connect the $1$-page to $H(C)$ is the following.

The basic downset we are studying in this section is the following.
For $0 < k \leq n$ and $p_1,\ldots,p_k\in\Zinf$ we define
\[
A(p_1,\ldots,p_k) := \{X\in\Zinf^n \st (x_1,\ldots,x_k) \leqlex (p_1,\ldots,p_k)\},
\]
where $\leqlex$ denotes the lexicographic order.

More generally we will consider linear transforms of $A(p_1,\ldots,p_k)$. Let $\varphi\in\GL(n,\ZZ)$ be a unimodular matrix with only non-negative integer entries.%
\footnote{Actually we only need that in every column of $\varphi$ the top non-zero entry is positive. However this yields no generalization, since $A(P_{1,\ldots,k};\varphi) = A(T_{1,\ldots k}(P);T\varphi)$ for any lower-triangular matrix $T$ with ones on the diagonal.}
Usual choices are $\varphi=\id_{\ZZ^n}$, which means no shearing, and the matrix with all ones on and above the diagonal, which will be used for the definition of the $2$-page in Section~\ref{secSecondaryConnection}.
Let $\varphi_{1\ldots k}(X)$ denote the first $k$ components of $\varphi(X)$.

For $0 < k \leq n$, and $P=(p_1,\ldots,p_k)\in\ZZ^k$ we define
\begin{align*}
A(p_1,\ldots,p_k;\varphi)\ :=\ & \{X\in\Zinf^n \st \varphi_{1\ldots k}(X) \leqlex (p_1,\ldots,p_k)\} \\
=\ & \{X\in\Zinf^n \st \varphi(X) \leqlex (p_1,\ldots,p_k,\infty^{n-k})\},
\end{align*}
which is a downset of $\Zinf^n$ since $\varphi\in(\ZZ_{\geq 0})^{n\times n}$.
Note that for $R\in \ZZ^n$,
\[
A(p_1+r_1,\ldots,p_k+r_k;\varphi) = A(p_1,\ldots,p_k;\varphi) + \varphi^{-1}(R).
\]

Let's fix an offset $Q=(q_1,\ldots,q_n)\in\ZZ^n$ with $Q\geqlex 0$.
Usual choices are $Q=0$, which means no offset, and $Q=\one:=(1,\ldots,1)$, which will be used in Section~\ref{secSecondaryConnection}.
Let $e_1,\ldots,e_n$ denote the standard basis vectors of $\ZZ^n$.
We define
\begin{align*}
p(p_1,\ldots,p_k;r,\varphi,Q) :=\ & A(p_1,\ldots,p_k;\varphi), \\
q(p_1,\ldots,p_k;r,\varphi,Q) :=\ & A(p_1,\ldots,p_k-1;\varphi) \\
                               =\ &  p(p_1,\ldots,p_k;r,\varphi,Q)-\varphi^{-1}(e_k), \\
z(p_1,\ldots,p_k;r,\varphi,Q) :=\ & A(p_1-q_1,\ldots,p_k-q_k-r;\varphi)\\
                               =\ &  p(p_1,\ldots,p_k;r,\varphi,Q)-\varphi^{-1}(Q+re_k),\\
b(p_1,\ldots,p_k;r,\varphi,Q) :=\ & A(p_1+q_1,\ldots,p_k+q_k+r-1;\varphi)\\
                               =\ &  q(p_1,\ldots,p_k;r,\varphi,Q)+\varphi^{-1}(Q+re_k).
\end{align*}
We extend the definitions for $k=0$ and $r=1$ by $p(-;1,\varphi,Q)=b(-;1,\varphi,Q)=C$ and $q(-;1,\varphi,Q)=z(-;1,\varphi,Q)=0$, where ``$-$'' denotes the empty sequence.
We write the associated $S$-terms as $S^{pz}_{bq}(p_1,\ldots,p_k;r,\varphi,Q)$.
Only the first $k$ components of $Q$ matter for this $S$-term.

If one likes to think of usual spectral sequences as a book, then one can think of these $S$-terms as words on page $r$ of chapter $k$.
There are $n$ chapters, which are counted downwards, every chapter starts on page $1$, with a page index shift given by $Q$.

\begin{lemma}[Big steps]
\label{lemLexConnection}
With the above definitions,
\begin{enumerate}
\item $S^{pz}_{bq}(p_1,\ldots,p_n;0,\varphi,0) = F_{A(p_1,\ldots,p_n;\varphi)}/F_{A(p_1,\ldots,p_n-1;\varphi)}$.
\item $S^{pz}_{bq}(p_1,\ldots,p_n;1,\varphi,0) = H_*(F_{A(p_1,\ldots,p_n;\varphi)}/F_{A(p_1,\ldots,p_n-1;\varphi)})$.
\item The differential $d:C\to C$ induces differentials
\[
S^{pz}_{bq}(p_1,\ldots,p_k;r,\varphi,Q)\to S^{pz}_{bq}(p_1-q_1,\ldots,p_k-q_k-r;r,\varphi,Q).
\]
Taking homology at $S^{pz}_{bq}(p_1,\ldots,p_k;r,\varphi,Q)$ with respect to this differential yields $S^{pz}_{bq}(p_1,\ldots,p_k;\allowbreak r+1,\varphi,Q)$.
\item If the filtration is finite, $S^{pz}_{bq}(p_1,\ldots,p_k;r,\varphi,Q)$ equals $S^{pz}_{bq}(p_1,\ldots,p_k;\infty,\varphi,Q)$ for $r$ large enough.
\item The terms $S^{pz}_{bq}(p_1,\ldots,p_k;\infty,\varphi,Q)$ are the filtration quotients of a $\ZZ$-filtration of $S^{pz}_{bq}(p_1,\ldots,p_{k-1};\allowbreak 1,\varphi,Q)$. More precisely, there is a canonical $\ZZ$-filtration $(F_i)_{i\in\ZZ}$
\[
0\subseteq \ldots \subseteq F_i\subseteq F_{i+1}\subseteq \ldots \subseteq S^{pz}_{bq}(p_1,\ldots,p_{k-1};1,\varphi,Q)
\]
such that $S^{pz}_{bq}(p_1,\ldots,p_k;\infty,\varphi,Q)\iso F_{p_k}/F_{p_k-1}$.
\item $S^{pz}_{bq}(-;1,\varphi,Q) = H(C)$.
\end{enumerate}
\end{lemma}

\begin{figure}[htb]
\centerline{
\includegraphics[scale=0.36]{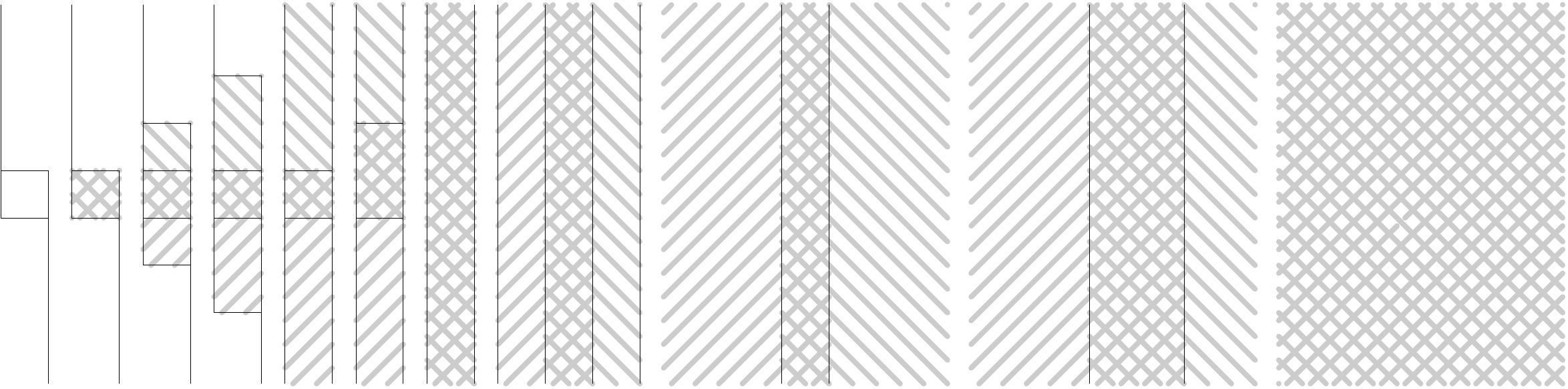}
} 
\caption{Some $S$-terms in the big step connection for $n=2$, $\varphi=\id_{\ZZ^2}$, and $Q=0$.
The first five $S$-terms are for $k=2$ and $r=0, 1, 2, 3$, and $\infty$.
The sixth $S$-term indicates the extension process.
The following three $S$-terms are for $k=1$ and $r=1, 2$, and $\infty$.
The last two $S$-terms indicate the second extension process, leading to $H(C)$.
}
\label{figBigSteps}
\end{figure}

The lemma says that starting with $S^{pz}_{bq}(p_1,\ldots,p_n;0,\varphi,Q)$ or $S^{pz}_{bq}(p_1,\ldots,p_n;1,\varphi,Q)$ (which are part of the $0$-page and $1$-page, respectively, in case the offset $Q$ is zero), we can repeat steps 3.), 4.), and 5.) for $k=n,n-1,\ldots,1$ until we arrive at 6.).
Permuting the coordinates of $\Zinf^n$ gives different connections.
See Figure~\ref{figBigSteps}, where the sets $p\wo z$ and $b\wo q$ are depicted by the ruled areas
\[
\includegraphics[scale=0.40]{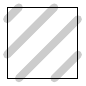}
\qquad \raisebox{0.5em}{\textnormal{and}} \qquad
\includegraphics[scale=0.40]{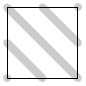},
\]
respectively.
They overlap in $p\wo q$, which is thus ruled in both ways.

Step 3.) in Lemma~\ref{lemLexConnection} increases $B$ by the set of boxes
\[
\begin{split}
b(p_1,\ldots,p_k;r+1,\varphi,Q)\wo b(p_1,\ldots,p_k;r,\varphi,Q) = \\
\{X\in \Zinf^n\st \varphi_{1\ldots k}(X) = (p_1+q_1,\ldots,p_k+q_k+r)\},
\end{split}
\]
which gets bigger and bigger as $k$ decreases from $n$ to $1$, and similarly with $Z$.
Alternatively, we could do small steps only, adding only a box at a time to $B$ and $Z$.

In order to recycle the above notation, in the next lemma we think of $Q$ being the page index, and we always let $r=0$.

\begin{lemma}[Small steps]
\label{lemLexConnectionInSmallSteps}
Suppose that the filtration is finite.
Let $P=(p_1,\ldots,p_n)\in\ZZ^n$ and $Q=(q_1,\ldots,q_n)\in\ZZ^n$.
\begin{enumerate}
\item $S^{pz}_{bq}(P;0,\varphi,0) = F_{A(P;\varphi)}/F_{A(P-e_n;\varphi)}$.
\item $S^{pz}_{bq}(P;0,\varphi,e_n) = H_*(F_{A(P;\varphi)}/F_{A(P-e_n;\varphi)})$.
\item The differential $d:C\to C$ induces differentials
\[
S^{pz}_{bq}(P;0,\varphi,Q)\to S^{pz}_{bq}(P-Q;0,\varphi,Q).
\]
Taking homology at $S^{pz}_{bq}(P;0,\varphi,Q)$ with respect to this differential yields the page $S^{pz}_{bq}(P;0,\varphi,Q+e_n)$.
\item If the filtration is finite and $q_k,\ldots,q_n$, $k\geq 2$, are large enough then $S^{pz}_{bq}(P;0,\varphi,Q)$ is isomorphic to $S^{pz}_{bq}(P;0,\varphi,Q^*)$ with $Q^*=(q_1,\ldots,q_{k-1},q_k+1,-q_{k+1},\ldots,-q_n)$.
\item If the filtration is finite and $q_1,\ldots,q_n$ are large enough then $S^{pz}_{bq}(P;0,\varphi,Q)$ is isomorphic to $S^{pz}_{bq}(P;0,\varphi,\infty^n)$.
\item The terms $S^{pz}_{bq}(P;0,\varphi,\infty^n)$ are filtration quotients of $H(C)$.
More precisely, there is a canonical filtration $(G_P)_{P\in\ZZ^n}$, $G_P\subseteq H(C)$, with $G_P\subseteq G_{P'}$ if $P\leqlex P'$, such that $S^{pz}_{bq}(P;0,\varphi,\infty^n)\iso G_P/G_{P-e_n}$.
\end{enumerate}
\end{lemma}

\begin{figure}[htb]
\centerline{
\includegraphics[scale=0.304]{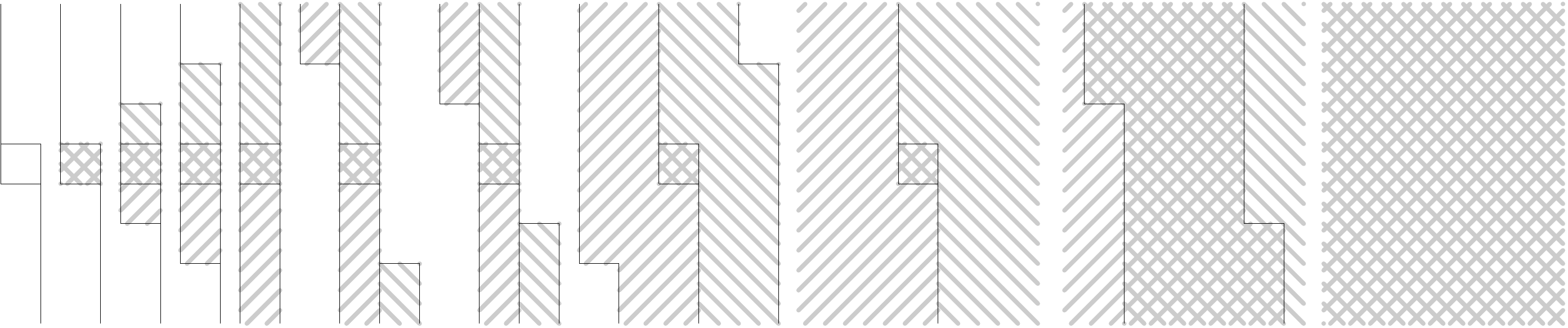}
} 
\caption{Some $S$-terms in the small step connection for $n=2$, $\varphi=\id_{\ZZ^2}$.
The first nine $S$-terms are for $Q=(0,0)$, $(0,1)$, $(0,2)$, $(0,3)$, $(0,\infty)$, $(1,-3)$, $(1,-2)$, $(2,2)$, and $(\infty,\infty)$.
}
\label{figSmallSteps}
\end{figure}

The lemma says that starting with 1.) or 2.) we can repeat steps 3.) and 4.) until we reach 5.), which can be used to compute 6.).
Note that in steps 3.) and 4.), $Q$ increases in the lexicographic order.
Hence for finite filtrations this is indeed a finite connection from the $1$-page to $H_*(C)$.
See Figure~\ref{figSmallSteps}.

Permuting the coordinates of $\Zinf^n$ gives different connections. We can even use different permutations for different $P$, until we reached $Q=\infty^n$. At that point 4.) we reached naturally isomorphic $S$-terms by Lemma~\ref{lemNaturallyIsomorphicSterms}, hence we can proceed with 5.)

Moreover we can merge the ideas of Lemmas~\ref{lemLexConnection} and~\ref{lemLexConnectionInSmallSteps} and get many more connections based on the lexicographic ordering.

\subsubsection{Secondary connection} \label{secSecondaryConnection}

The following connection between the $1$-page and $H(C)$ will be in particular useful for successive spectral sequences for which the limit of one spectral sequence is the \emph{second} page of the next one.

We start with connecting the $1$-page to the `$2$-page', for which we actually need to take $n$ times homology.



Let $\varphi_k:\ZZ^n\to\ZZ^n$, $0\leq k\leq n$, be the automorphism given by
\begin{equation}
\label{eqAutoPhi}
\varphi_k(X) = (x_{k+1},x_{k+2},\ldots,x_{n},\sum_{i=1}^kx_i,\sum_{i=2}^kx_i,\ldots,x_k).
\end{equation}
For $P\in\Zinf^n$ and $1\leq k\leq n$, we define 
\[
T^k_P := A(\varphi_k(P);\varphi_k) = \{X\in\Zinf^n\st\varphi_k(X)\leqlex\varphi_k(P)\}.
\]
In other words, $A(0^n;\id) = \varphi_k(T^k_P-P)$.
Explicitly,
\[
\begin{split}
T^k_P =
\big\{X\in\Zinf^n \st \big((x_{k+1},\ldots,x_n)=(p_{k+1},\ldots,p_n) \textnormal{ and } 
\sum_{i=1}^k x_i + \delta_{X_{1\ldots k}\lesslex P_{1\ldots k}} \leq \sum_{i=1}^k p_i
\big)
\\
\textnormal{ or } (x_{k+1},\ldots,x_n)\lesslex (p_{k+1},\ldots,p_n)\big\}.
\end{split}
\]


Note that
\[
T^k_P \wo \{P\} = T^k_{P+e_{k-1}-e_k}.
\]
Here we set $e_0:=0$, $e_{-1}:=-e_n$, such that we don't need to treat the cases $k=0$, $k=1$, and $k\geq 2$ separately.

Given $P=(p_1,\ldots,p_n)\in\Zinf^n$ and $0\leq k\leq n$, we define
\begin{align*}
p(P;k) :=\ & T^k_P \\
q(P;k) :=\ & p(P;k)\wo \{P\} = T^k_{P+e_{k-1}-e_k} \\
z(P;k) :=\ & p(P;k)-e_k = T^k_{P-e_k} \\
b(P;k) :=\ & q(P;k)+e_k = T^k_{P+e_{k-1}} \\
z^*(P;k) :=\ & z(P;k)\wo \{P-e_k\} = T^k_{P+e_{k-1}-2e_k} \\
b^*(P;k) :=\ & b(P;k)\cup \{P+e_k\} = T^k_{P+e_k}
\end{align*}
We denote the associated $S$-terms as $S^{pz}_{bq}(P;k)$ and $S^{pz^*}_{b^*q}(P;k)$.

\begin{example}
\label{exSecondaryConnectionPages012}
For $k=0$ we have $\varphi_0=\id_{\ZZ^n}$ and the above $S$-terms are
\begin{equation}
\label{eqSecConnS0}
S^{pz}_{bq}(P;0) = F_{A(P)}/F_{A(P)\wo\{P\}},
\end{equation}
and
\begin{equation}
\label{eqSecConnS*0}
S^{pz^*}_{b^*q}(P;0) = H(F_{A(P)}/F_{A(P)\wo\{P\}}).
\end{equation}
For $k=n$, $\varphi_n$ is the matrix with ones on and above the diagonal and zeros otherwise, and we have
\begin{equation}
\label{eqSecConnSn}
S^{pz}_{bq}(P;n) = S^{pz}_{bq}(\varphi_n(P);0,\varphi_n,\one),
\end{equation}
and
\begin{equation}
\label{eqSecConnS*n}
S^{pz^*}_{b^*q}(P;n) = S^{pz}_{bq}(\varphi_n(P);1,\varphi_n,\one).
\end{equation}
We call \eqref{eqSecConnS*n} the \emph{$2$-page}.
\end{example}

\begin{lemma}
\label{lemSecConnDifferentials}
There are differentials ``in direction $-e_k$'':
\[
d:S^{pz}_{bq}(P;k) \To S^{pz}_{bq}(P-e_k;k).
\]
Taking homology at $S^{pz}_{bq}(P;k)$ yields $S^{pz^*}_{b^*q}(P;k)$.
\end{lemma}

\begin{proof}
We only need to check that $q(P)=b(P-e_k)$, $z(P)=p(P-e_k)$, $z^*(P)=q(P-e_k)$, and $b^*(P)=p(P+e_k)$.
\end{proof}

The inverse of $\varphi_k$ is
\[
\varphi_k^{-1}(y) = (y_{k+1}-y_{k+2},y_{k+2}-y_{k+3},\ldots,y_{n-1}-y_{n},y_n, y_1, y_2,\ldots, y_{n-k}).
\]
Since all entries in $\varphi_k$ are non-negative, $\varphi_k$ preserves the order of $\Zinf^n$, however its inverse doesn't.
This means that we can take kernels and cokernels as we did with $A(p_1,\ldots,p_k;\id)$ in Lemma~\ref{lemLexConnection}, but we might expect more natural isomorphisms.
And indeed that's what we will exploit:

Let's denote
\begin{align*}
Z(P;k)   :=\ & Z(z(P;k),q(P;k),p(P;k),b(P;k)), \\
Z^*(P;k) :=\ & Z(z^*(P;k),q(P;k),p(P;k),b^*(P;k)),
\end{align*}
and similarly $B(P;k)$ and $B^*(P;k)$.

Let $V_k\subseteq\{-1,0,1\}^k\times \{0\}^{n-k}$ be the set of all $\{-1,0,1\}$-vectors of length $n$ with the following properties: The last $n-k$ entries are zero,  the last non-zero entry is $+1$, and $+1$ and $-1$ appear alternating.
In other words, $V_0:=\{0\}\subset \ZZ^n$ and
\[
V_{k+1} := V_k~\dot\cup~(e_{k+1}-V_k) = \{0\}~\dot\cup~(e_1-V_0)~\dot\cup~\ldots~\dot\cup~(e_{k+1}-V_k).
\]
Put $V_{-1}:=\emptyset$.

\begin{lemma}
\label{lemSecConnNatIsos}
For $0\leq k\leq n$,
\[
\begin{array}{lll}
Z(P;k) = P-V_{k-1}, & \ \ \ \ \ & B(P;k) = P+V_{k-1}, \\
Z^*(P;k) = P-V_k,   & \ \ \ \ \ & B^*(P;k) = P+V_k.
\end{array}
\]
Thus for $0\leq k\leq n-1$ we have a natural isomorphism
\[
S^{pz^*}_{b^*q}(P;k) \iso S^{pz}_{bq}(P;k+1).
\]
\end{lemma}

\begin{proof}
It suffices to prove the equalities for $B(0;k)$ and~$B^*(0;k)$.
The special cases $B(0;0)$, $B(0;1)$, and $B^*(0;0)$ should be checked separately, since they involve particularly defined terms $e_{-1}$, $e_{0}$, and~$V_{-1}$.

Let's write $b:=b(0;k)$, $b^*:=b^*(0;k)$, and $q:=q(0;k)$.
$B(0;k)$ is the component of $b\wo q$ containing~$0$, where
\begin{align*}
b\wo q
=\ & \{X\st \varphi_k(0)\leqlex \varphi_k(X)\leqlex\varphi_k(e_{k-1})\} \\
=\ & \{X\st 0 \leqlex \varphi_k(X)\leqlex e_{n-k+1}+\ldots+e_{n-1}\}.
\end{align*}
All such $X$ satisfy $x_{k+1}=\ldots=x_n=0$.
Clearly, $V_{k-1}$ is connected and $0\in V_{k-1}\subseteq b\wo q$.
If $V_{k-1}\neq b\wo q$, then there exists an $X\in (b\wo q)\wo V_{k-1}$ that is adjacent to $V_{k-1}$.
Then for exactly one $1\leq j\leq k$, either $x_j=\pm 2$ or ($x_j=\pm 1$ and the next non-zero entry of $X$ equals $x_j$).
For this $j$, $\phi_k(X)_{n-k+j}=\sum_{i=j}^kx_i$ is less than $0$ or larger than $1$, thus $X\not\in b\wo q$, which is the desired contradiction.

A very similar argument works for $B^*(0;k)$ using
\begin{align*}
b^*\wo q
=\ & \{X\st \varphi_k(0)\leqlex \varphi_k(X)\leqlex\varphi_k(e_{k})\} \\
=\ & \{X\st 0 \leqlex \varphi_k(X)\leqlex e_{n-k+1}+\ldots+e_{n}\}.
\end{align*}
\end{proof}

\begin{figure}[htb]
\centerline{
\includegraphics[scale=0.25]{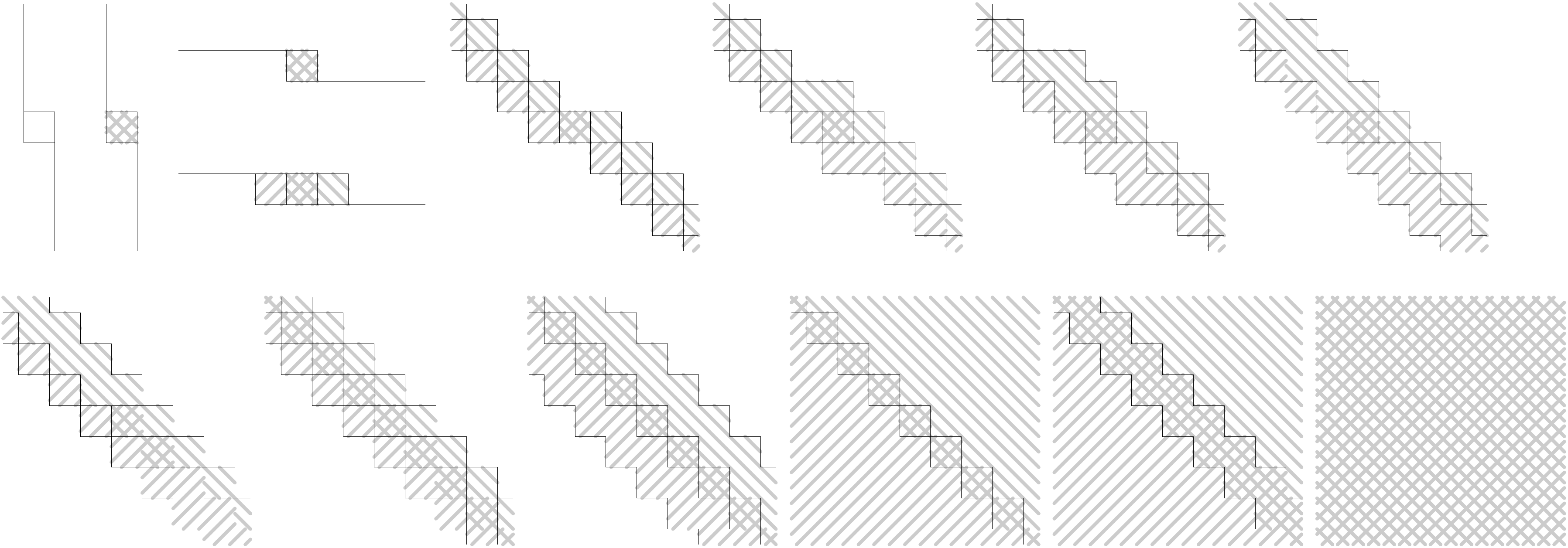}
} 
\caption{The six $S$-terms in the secondary connection for $n=2$, followed by 
some $S$-terms in the big step lexicographic connection with offset $Q=(1,1)$ and shearing matrix $\varphi_2$. The sixth $S$-term on the top is the $2$-page.}
\label{figSecondaryConnection}
\end{figure}

Using the lemmas alternately, we can connect the $0$-page \eqref{eqSecConnS0} and the $1$-page \eqref{eqSecConnS*0} to the $2$-page \eqref{eqSecConnS*n}.
From there we can proceed with Lemma~\ref{lemLexConnection} (or Lemma~\ref{lemLexConnectionInSmallSteps} if small steps are preferred) to connect $S^{pz^*}_{b^*q}(P;n)$ to $H_*(C)$.
See Figure~\ref{figSecondaryConnection}.
In this big step lexicographic connection, on page $r=1,\ldots$ in chapter $k=n,\ldots,1$, the differential has direction $-\varphi_n^{-1}(re_k+Q+\langle e_{k+1},\ldots,e_n\rangle)$.
One particular vector in this $(n-k)$-dimensional affine set is $re_{k-1}-(r+1)e_k = -\varphi_n^{-1}(re_k+(1^k,0^{n-k}))$.

\subsubsection{Generalized secondary connection}

In the secondary connection from the previous section, one starts from $S^{pz^*}_{b^*q}(P;0)\iso S^{pz}_{bq}(P;1)$, takes homology in direction $-e_1$, which yields $S^{pz^*}_{b^*q}(P;1)\iso S^{pz}_{bq}(P;2)$, then one takes homology in direction $-e_2$, and so on, until one arrives at the second page $S^{pz^*}_{b^*q}(P;n)$. In other words, we apply homology and natural isomorphisms alternatingly.

More generally, for arbitrary $Q=(q_1,\ldots,q_n)\geq 0$ we can start with $S^{pz}_{bq}(P;1)$, take $q_1$ times homology in direction $-e_1$, find a natural isomorphic page for which we can apply $q_2$ times homology in a direction close to $-e_2$, and so on.
The secondary connection is the special case $Q=\one$, which might be more useful than all other cases together.
Here is how it works.

Let $\varphi_k^Q\in\GL(n,\ZZ)$, $0\leq k\leq n$, $Q\in\ZZ^n_{\geq 0}$, be given by
\[
\varphi_k^Q(X):=(x_{k+1},\ldots,x_n,x_1+q_1\sum_{i=2}^k x_i,x_2+q_2\sum_{i=3}^k x_i,\ldots,x_k).
\]

We will recycle the notation from Section~\ref{secLexicographicConnection} for the $S$-terms.
Let $Q[n-k]:=(0^{n-k},q_1,\ldots,q_k)$.

\begin{lemma}
\label{lemGeneralizedSecConnNatIsos}
For $1\leq k\leq n$, $P\in\ZZ^n$, there is a natural isomorphism of $S$-terms
\[
S^{pz}_{bq}(\varphi_{k-1}^Q(P);1      ,\varphi_{k-1}^Q,Q[n-k+1])\iso 
S^{pz}_{bq}(\varphi_{k  }^Q(P);1-q_{k},\varphi_{k  }^Q,Q[n-k]).
\]
\end{lemma}

\begin{proof}[Proof (sketch)]
Define $a_0^Q:=0$ and
$
a_k^Q:=e_k-\sum_{i=1}^{k-1}q_ia_i^Q \in\ZZ^n.
$ 
Explicitly, 
\[
a_k^Q=\sum_{i=1}^k(-1)^{k-i}e_i\prod_{j=i}^{k-1}(q_j-\delta_{j>i}).
\]
One can show that
\[
(\varphi_k^Q)^{-1} =
\left(
\begin{array}{c|c}
\begin{matrix}
0^{k\times(n-k)} \\
\hline
\id_{\ZZ^{n-k}}
\end{matrix}
& 
a_1^Q  \cdots  a_k^Q
\end{array}
\right).
\]

Further, let $V_0^Q:=\{0\}$ and
$
V_{k}^Q:=\{0,\ldots,q_{k}\}\cdot a_{k}^Q + V_{k-1}^Q \subseteq\ZZ^n.
$ 
Thus $V_k^Q$ is a `discrete affine cube', 
\begin{equation}
\label{eqGenSecConnExplicitEqForVkQ}
V_k^Q = (\varphi_k^Q)^{-1}(\{0\}^{n-k}\times\{0,\ldots,q_1\}\times\ldots\times \{0,\ldots,q_k\}).
\end{equation}
For example for $Q=\one$ we obtain $\varphi_k^Q=\varphi_k$, $a_k^Q=e_k-e_{k-1}$, and $V_k^Q=V_k$.
By induction one proves that $V_{k-1}^Q = (e_{k}-a_{k}^Q) - V_{k-1}^Q$.
In particular this implies that $V_k^Q$ is connected.

For $S^{pz}_{bq}(0;1,\varphi_{k  }^Q,Q[n-k])$, we have
\begin{align*}
b\wo q
=\ & \{X\st 0\leqlex \varphi_k^Q(X)\leqlex Q[n-k]\}.
\end{align*}
From \eqref{eqGenSecConnExplicitEqForVkQ} it follows $0\in V_k^Q\subseteq b\wo q$.
Let $X\in (b\wo q)$ be adjacent to some $Y=X\pm e_j\in V_k^Q$.
Then $j\leq k$ since otherwise $\varphi_k^Q(e_j)=e_{j-k}$.
Suppose $Y=X+e_j$.
(The case $Y=X-e_j$ is analogous.)
Then $\varphi_k^Q(Y)=\varphi_k^Q(X)+(0^{n-k},q_1,\ldots,q_{j-1},1,0^k-j)$.
Since $X,Y\in b\wo q$, we get $\varphi_k^Q(X)_{n-k+i}=q_i$ for $1\leq i <j$ and $\leq q_i$ for $j\leq i\leq k$.
Hence, $X$ is already in $V_k^Q$.
Therefore, $V_k^Q$ is the connected component of $b\wo q$ that contains $0$.

We get the same statement for $S^{pz}_{bq}(0;1-q_{k+1},\varphi_{k+1}^Q,Q[n-k-1])$ using a similar argument.
Hence both $S$-terms are naturally isomorphic.
\end{proof}

Applying $q_{k}$ times homology to
\[
S^{pz}_{bq}(\varphi_{k}^Q(P);1-q_{k},\varphi_{k}^Q,Q[n-k])
\] 
as in Lemma~\ref{lemLexConnection}.3.) yields
\[
S^{pz}_{bq}(\varphi_{k}^Q(P);1,\varphi_{k}^Q,Q[n-k]).
\]
Therefore we can apply Lemma~\ref{lemLexConnection}.3.) and Lemma~\ref{lemGeneralizedSecConnNatIsos} alternatingly in order to connect the $0$-page \eqref{eqSecConnS0} and the $1$-page \eqref{eqSecConnS*0} to
\begin{equation}
\label{eqGeneralizedSecondPage}
S^{pz}_{bq}(\varphi_{n}^Q(P);0,\varphi_{n}^Q,Q),
\end{equation}
which we call the $2_Q$-page.
Again, using the lexicographic connections we can connect \eqref{eqGeneralizedSecondPage} to $H(C)$.

\subsubsection{Further connections}

There are many more ways to connect the $1$-page to $H_*(C)$, which might be useful sometimes.
Moreover if we know from the given $C$ that certain terms on the $0$-page vanish, more $S$-terms can be identified, which might give even more useful connections.

\begin{lemma}
Let $z\leq q\leq p\leq b$ be downsets in $D(\Zinf^n)$ such that $b\wo z$ is finite.
Then there is a connection from the $1$-page to $S^{pz}_{bq}$ using kernels and cokernels only.
\end{lemma}

\begin{proof}
The proof is by induction on $|b\wo z|$.
If $z=q$ and $p=b$ then $S^{pz}_{bq}=H_*(F_p/F_q)$ is by definition a term on the $1$-page.
If $p\neq b$, then for any coordinate-wise maximal element $X$ in $b\wo p$ we have
\[
S^{pz}_{bq} = \coker\left(S^{b,p}_{b^*,b\wo X}\To S^{p,z}_{b\wo X,q} \right),
\]
$b^*\geq b$ being arbitrary, for example $b^*=b$.
If $z\neq q$, then for any coordinate-wise minimal element $X$ in $q\wo z$ we have
\[
S^{pz}_{bq} = \ker\left(S^{p,z\cup X}_{b,q}\To S^{z\cup X,z^*}_{q,z} \right),
\]
$z^*\leq z$ being arbitrary, for example $z^*=z$.
\end{proof}

This lemma uses kernels and cokernels.
In case one only wants to allow natural isomorphisms and taking homology as in the standard spectral sequence, the following happens:
We start with $S^{pz}_{bq}$-terms from the $1$-page with $|p\wo q|=1$.
Then there are $n$ possible differentials that enlarge $Z$ and $B$ by one box.
We conjecture that this can be iterated, and at every step we have exactly $n$ possible differentials that enlarge $Z$ and $B$ by at least one box when taking homology.
However these operations seem to be not commutative at all, which makes them probably only useful in the above discussed special cases of secondary connections.

Note further that any $S$-term from the $1$-page can be obtained by extensions (and limits in the case of non-finite filtrations) from $1$-page $S$-terms $S^{pp}_{qq}$ with $|p\wo q|=1$.
The latter terms are naturally isomorphic to terms $S^{pp}_{qq}$ with $p=A(p_1,\ldots,p_n)$ and $q=p\wo\{(p_1,\ldots,p_n)\}$.

\begin{lemma}
Let $z\leq q< p\leq b$ be downsets in $D(\Zinf^n)$ such that $p\wo q$ is finite.
Then $S^{pz}_{bq}$ is an iterated extension by terms $S^{p^*z}_{bq^*}$ of the form $q\leq q^*<p^*\leq p$ with $|p^*\wo q^*|=1$.
\end{lemma}

\begin{proof}
This is immediate from the extension property.
\end{proof}

\subsection{Example: $n$ independent subcomplexes}
\label{secNindependentSubcomplexes}

Consider $n$ subcomplexes $F_1,\ldots,F_n$ of a given chain complex $C$.
The goal is to compute $H(C)$.
For any $1\leq i\leq n$ we get a $\ZZ$-filtration $(F_k^{(i)})_{k\in\ZZ}$ of $C$ with $F_0^{(i)}:= F_i$ and $F_{-k}^{(i)}:=0$ and $F_{k}^{(i)}:=C$ for $k\geq 1$.
We further assume that the induced $D(\Zinf^n)$-filtration is distributive.
It induces a spectral system over $D(\Zinf^n)$.

\begin{example}
$C$ might arise as the singular chain complex of a space $X$, and the $F_i$ as the chain complexes of $n$ open subspaces $X_i$.
The induced filtration is in general not distributive, but this problem can be ignored, because the associated exact couple system is excisive (compare with Example~\ref{exNindependentSubspacesAndGeneralizedHomology}, which also deals with the case where the goal of computation is $h_*(X)$ for some generalized homology theory $h_*$).
\end{example}

Of course it makes sense to look only at the subposet $I$ that is given by all downsets of elements in $\{0,1\}^n$ and their unions.
$I$ is isomorphic to the poset of abstract simplicial complexes on $\{1,\ldots,n\}$, ordered by inclusion.
Here we need to distinguish between the empty complex $\{\emptyset\}$, which corresponds to the downset of $0^n$ in $D(\Zinf^n)$, and the void complex $\emptyset$, which corresponds to $-\infty$. 

Let us consider only the secondary connection (see Section~\ref{secSecondaryConnection}).
The only potentially non-zero $S$-terms $S^{pz}_{bq}(P;0)$ are those where $P\in\{0,1\}^n$, and they are given by
\[
S^{pz}_{bq}(P;0)
= F_P/\sum_{Q<P} F_Q 
\iso (\bigcap_{Q>P} F^Q) / F^P,
\]
where $F_P=\bigcap_{i:p_i=0} F_i$, $F^P=\sum_{i:p_i=1} F_i$, $Q$ runs through $\{0,1\}^n$, and empty intersections mean $C$.

Taking homology of each of these $2^n$ chain complexes yields $S^{pz^*}_{b^*q}(P;0)\iso S^{pz}_{bq}(P;1)$, $P\in\{0,1\}^n$.
Then the secondary connection proceeds with an iteration over $k$ from $k=1,\ldots,n$:
At step $k$, we have $2^{n-1}$ chain complexes
\[
0\to S^{pz}_{bq}(P+e_k;k) \tof{d} S^{pz}_{bq}(P;k) \to 0, \ \ \ 
(P\in\{0,1\}^n, p_k=0).
\]
Taking homology yields $S^{pz^*}_{b^*q}(P+e_k;k)\iso S^{pz}_{bq}(P+e_k;k+1)$ and $S^{pz^*}_{b^*q}(P;k)\iso S^{pz}_{bq}(P;k+1)$, respectively.
After step $k=n$, we arrived at the second page $S^{pz^*}_{b^*q}(P;n)=S^{pz}_{bq}(\varphi_n(P);1,\varphi_n,\one)$.

We can proceed with the lexicographic connection in big steps (\eqref{eqSecConnSn} and Lemma~\ref{lemLexConnection}): Iterating over $k$ from $k=n,\ldots,1$, at step $k$ we have to take homology with respect to the differentials in Lemma~\ref{lemLexConnection}.3) only for the pages $r=1,\ldots,n-k$, and then we group the pages and unify these groups using extensions (Lemma~\ref{lemLexConnection}.5).

After step $k=1$, we arrive at one remaining $S$-term
, which is isomorphic to $H(C)$.

\medskip

This spectral system should be compared to the generalized Mayer--Vietoris sequence, which gives rise to the Leray--Mayer--Vietoris spectral sequence.
This spectral sequence is obtained by restricting the above filtration poset $I$ to the skeletons of the simplex on vertex set $\{1,\ldots,n\}$.
It converges to $H(\sum F_i)$ or to $H(C)$, depending on the convention whether or not one includes the  full simplex in the filtration.


\section{Exact couple systems} \label{secExactCoupleAnalogs}

Massey \cite{Mas52exactCouples,Mas52exactCouples2} constructed the framework of exact couples, which give more generally rise to spectral sequences than $\ZZ$-filtrated chain complexes.
So we may ask whether our generalized spectral sequences can also start directly from the $1$-page, without being constructed from the $0$-page.

In this section we construct basic data and axioms similar 
to exact couples that induce a spectral system as above that starts from the $1$-page, that is, it contains only terms $S^{pz}_{bq}$ with $z\leq q\leq p\leq b$.

We call a poset $I$ bounded if it has a minimum and a maximum, which we denote by $-\infty$ and $\infty$.
We regard $I$ as a category whose objects are the elements and whose morphisms $q\to p$ are the relations $q\leq p$.
For $n\geq 2$, let $I_n$ denote the poset of $n$-tuples $(p_1,\ldots,p_n)\in I^n$ with $p_1\geq p_2\geq\ldots\geq p_n$, ordered by $(p_1,\ldots,p_n)\geq (p'_1,\ldots,p'_n)$ if and only if $p_i\geq p'_i$ for all $i$.
As with $I$, $I_n$ is also a category.

\begin{definition}
An \textdef{exact couple system} over a bounded poset $I$ consists of the following data:
\begin{enumerate}
\item 
\label{itECSdataE}
A functor $E:I_2\to \abcat$ from $I_2$ to the category of abelian groups.
We write $E^p_q:=E(p,q)$, $D_p:= E(p,-\infty)$; and $i_{qp}:D_q\to D_p$, $j_{pq}:D_p\to E^p_q$, and $\ell^{pq}_{p'q'}:E^p_q\to E^{p'}_{q'}$ for the maps $E((q,-\infty)\leq (p,-\infty))$, $E((p,-\infty)\leq (p,q))$, and $E((p,q)\leq (p',q'))$, respectively.
\item 
\label{itECSdataK}
Maps $k_{pq}:E^p_q\to D_q$ for all $(p,q)\in I_2$.
\end{enumerate}
We require that it satisfies the following axioms:
\begin{enumerate}
\item 
\label{itECSaxiomTriangle}
The triangles
\[
\xymatrix{
D_q \ar[rr]^{i_{qp}} &                        & D_p \ar[dl]^{j_{pq}} \\
                     & E^p_q \ar[ul]^{k_{pq}} &
}
\]
are exact.
\item 
\label{itECSaxiomSquare}
The diagrams
\[
\xymatrix{
E^p_q \ar[d]_{\ell^{pq}_{p'q'}} \ar[r]^{k_{pq}} & D_q \ar[d]^{i_{qq'}} \\
E^{p'}_{q'} \ar[r]_{k_{p'q'}}                   & D_{q'}
}
\]
commute.
\end{enumerate}
\end{definition}

\begin{remark}
Equivalently, we could define an exact couple system over $I$ as a functor $\wt E: \wt I_2\to \abcat$, where $\wt I_2$ is the same category as $I_2$ except that we add morphisms $\wt k_{pq}:(p,q)\to (q,-\infty)$ for all $(p,q)\in I_2$ (and all resulting compositions) subject to the commutation relation 
$((p,-\infty)\leq (p',-\infty))\circ \wt k_{pq} = \wt k_{p'q'}\circ ((p,q)\leq (p',q'))$.
Datum~\ref{itECSdataK} and Axiom~\ref{itECSaxiomSquare} are then automatically given and we only need to require Axiom~\ref{itECSaxiomTriangle}.

In this wording it is also clear what maps between two exact couple systems $E:\wt I_2\to \abcat$ and $E':\wt I'_2\to \abcat$ should be, namely an order preserving map $f:I\to I'$ together with a natural transformation from $E$ to $E'\circ (f\times f)|_{\wt I_2}$.
In what follows, we call constructions natural if they commute with maps between exact couple systems.

More generally we can take any abelian category in place of $\abcat$.
\end{remark}

\begin{remark}[`Cohomological' definition]
\label{remExactCoupleSystemCohomologicalDefinition}
It will follow from Lemma~\ref{lemExactTriangles} below that 
we could equivalently define an exact couple system as a functor $E:I_2\to\abcat$ with $E^p_q$ and $\ell^{pq}_{p'q'}$ as before, $D^p:=E(\infty,p)$, together with boundary maps $k^{pq}: D^p\to E^p_q$, such that the triangles
\[
\xymatrix{
D^q \ar[rr]^{\ell^{\infty q}_{\infty p}} &   & D^p \ar[dl]^{k^{pq}} \\
 & E^p_q \ar[ul]^{\ell^{pq}_{\infty q}} &
}
\]
are exact and $\ell^{pq}_{p'q'}\circ k^{pq} = k^{p'q'}\circ\ell^{\infty p}_{\infty p'}$ for all $(p,q)\leq (p',q')$.
\end{remark}

\begin{example}[Exact couple system of an $I$-filtered chain complex]
\label{exExactCoupleOfFilteredChainComplex}
The spectral system of an $I$-filtered chain complex $C$ contains an exact couple system over $I$ with $E^p_q := H(F_p/F_q)$, $D_p := H(F_p)$, $i$ and $j$ are induced by inclusion, and $k$ is the connecting homomorphism.
\end{example}

\begin{example}[Exact couple system of an $I$-filtered space]
\label{exExactCoupleOfFilteredSpace}
Let $X$ be a topological space that is filtered by a bounded poset~$I$.
That is, we are given a family of closed subspaces $(X_i)_{i\in I}$ of $X$ with inclusions $X_q\incl X_p$ whenever $q\leq p$.
We may assume that $X_{-\infty}=\emptyset$ and $X_\infty=X$.

Then for any generalized homology theory $h_*$, $E^p_q:=h_*(X_p,X_q)$ is an exact couple system over~$I$.
Here, $D_p=h_*(X_p)$, $i$ and $j$ are induced by inclusions, and $k_{pq}$ is the connecting homomorphism in the long exact sequence of the pair $(X_p,X_q)$.

Analogously, for any generalized cohomology theory $h^*$, $E^p_q:=h^*(X_q,X_p)$ is an exact couple system over~$I^*$, where $I^*$ denotes the dual poset of $I^*$, that is, $p\leq_I q$ if and only if $q\leq_{I^*}p$.
Since $\pm\infty$ denotes minimum and maximum, $\infty_I=-\infty_{I^*}$ and $-\infty_I=\infty_{I^*}$.
Thus, $X_{-\infty}=X$, $X_\infty=\emptyset$, $D_p=h^*(X,X_p)$, $i$ and $j$ are induced by inclusions, and $k_{pq}$ is the connecting homomorphism in the long exact sequence of the triple $(X,X_q,X_p)$.
If we use the exact couple system definition from Remark~\ref{remExactCoupleSystemCohomologicalDefinition}, we have $D^p = h^*(X_p)$.

\end{example}

\begin{example}[Perverse sheaves and F\'ary functors]
Perverse sheaves and more generally F\'ary functors have naturally the structure of an exact couple systems (they have more properties, and $I$ is a set of open sets of a space $X$), see Be{\u\i}linson--Bernstein--Deligne~\cite{BBD82perverseSheaves} and Grinberg--MacPherson~\cite{GriMac99EulerCharactAndLagrangIntersect}.
\end{example}

\begin{example}[Exact couples]
An exact couple system $E$ over $\ZZ$ naturally contains an exact couple, namely the collection of all exact triangles from Axiom~\ref{itECSaxiomTriangle} with $p=q+1$.
All (columns of the) pages of the induced spectral sequence are particular $S$-terms of~$E$, as in Example~\ref{secSSofZfiltration}.
\end{example}

\subsection{Basic properties of exact couple systems}

Now consider an exact couple system $E$ over $I$.
We define differentials
\[
d_{pqz}: E^p_q\to E^q_z.
\]
for all $(p,q,z)\in I_3$ by setting $d_{pqz}:=j_{qz}\circ k_{pq}$.
Composing two such differentials yields the zero map since $k\circ j=0$ by Axiom~\ref{itECSaxiomTriangle}.
Moreover we have $d\circ \ell = \ell\circ d$.

\begin{lemma}[Exact triangles]
\label{lemExactTriangles}
For any $p_1\leq p_2\leq p_3$ in $I$ there is an exact triangle
\[
\xymatrix{
E^{p_2}_{p_1}\ar[r]^\ell & E^{p_3}_{p_1}\ar[r]^\ell & E^{p_3}_{p_2} \ar@/^2pc/[ll]^d.
}
\]
\end{lemma}

The lemma is analogous to the octahedral axiom in triangulated categories, since it says that the center triangle in the following diagram is exact if all three outer triangles are.
\[
\xymatrix{
& & D_{p_1} \ar@/_2pc/[ddll]_i \ar@/^2pc/[ddrr]^i & & \\
& E^{p_2}_{p_1} \ar[ur]^k \ar[rr]_\ell & & E^{p_3}_{p_1} \ar[ul]_k \ar[dl]_\ell & \\
D_{p_2} \ar@/_2pc/[rrrr]_i \ar[ur]^j & & E^{p_3}_{p_2} \ar[ll]^k \ar[ul]_d & & D_{p_3} \ar[ll]^j \ar[ul]_j
}
\]

\begin{proof}
For the exactness at $E^{p_2}_{p_1}$, we chase the diagram
\[
\xymatrix{
                       &                            & D_{p_3} \ar[rd]^j                  & \\
E^{p_3}_{p_2} \ar[r]^k & D_{p_2} \ar[r]^j \ar[ur]^i & E^{p_2}_{p_1} \ar[r]^\ell \ar[d]^k & E^{p_3}_{p_1} \ar[dl]^k \\
                       &                            & D_{p_1} \ar[ul]^i.
}
\]
Let $x\in E^{p_2}_{p_1}$ with $\ell(x)=0\in E^{p_3}_{p_1}$.
Then $k(x)=0\in D_{p_1}$, so $x=j(b)$ for some $b\in D_{p_2}$.
Let $c:=i(b)\in D_{p_3}$.
Now, $j(c)=j(i(b))=\ell(j(b))=\ell(x)=0\in E^{p_3}_{p_1}$.
Thus there is an $a\in D_{p_1}$ with $i(i(a))=c$.
Hence $i(b-i(a))=0\in D_{p_3}$.
So there exists $z\in E^{p_3}_{p_2}$ with $k(z)=b-i(a)$.
Therefore $d(z)=j(k(z))=j(b)-j(i(a))=j(b)=x$.

On the other hand, if $x=d(z)\in E^{p_2}_{p_1}$ with $z\in E^{p_3}_{p_1}$, then $\ell(x)=\ell(d(z))=j(i(k(z))=0\in E^{p_3}_{p_1}$ since $i\circ k=0$.

The proof for the exactness at $E^{p_3}_{p_2}$ works similarly in the diagram
\[
\xymatrix{
                                    & D_{p_1} \ar[rd]^i          & \\
E^{p_3}_{p_1} \ar[r]^\ell \ar[ur]^k & E^{p_3}_{p_2} \ar[r]^k     & D_{p_2} \ar[r]^j \ar[dl]^i & E^{p_2}_{p_1} \\
                                    & D_{p_3} \ar[ul]^j \ar[u]^j.
}
\]

For the exactness at $E^{p_3}_{p_1}$, we chase the commutative diagram
\[
\xymatrix{
D_{p_2} \ar[r]^i \ar[d]_j           & D_{p_3} \ar[d]^j \ar[dr]^j         & \\
E^{p_2}_{p_1} \ar[r]^\ell \ar[dr]_k & E^{p_3}_{p_1} \ar[r]^\ell \ar[d]^k & E^{p_3}_{p_2} \ar[d]^k \\
                                    & D_{p_1} \ar[r]^i                   & D_{p_2}.
}
\]
Let $y\in E^{p_3}_{p_1}$ with $\ell(y)=0\in E^{p_3}_{p_2}$.
Let $a=k(y)\in D_{p_1}$.
Then $i(a)=k(\ell(y))=0\in D_{p_2}$.
Thus $a=k(x)$ for some $x\in E^{p_2}_{p_1}$, and $k(\ell(x)-y)=0\in D_{p_1}$.
Hence there is a $c\in D_{p_3}$ with $j(c)=\ell(x)-y$, which satisfies $\ell(j(c))=0\in E^{p_3}_{p_2}$.
Thus there is a $b\in D_{p_2}$ with $i(b)=c$, and we have $j(i(b))=j(c)=\ell(x)-y$.
Let $x':=j(b)\in E^{p_2}_{p_1}$.
Then $\ell(x')=j(i(b))=\ell(x)-y$.
Therefore $y=\ell(x-x')$.

On the other hand, let $y=\ell(x)\in E^{p_3}_{p_1}$ with $x\in E^{p_2}_{p_1}$.
Let $z:=\ell(y)\in E^{p_3}_{p_2}$.
Then $k(z)=i(k(x))=0\in D_{p_2}$.
Thus there exists $c\in D_{p_3}$ with $j(c)=z\in E^{p_3}_{p_2}$.
Let $y':=j(c)\in E^{p_3}_{p_1}$.
Then $\ell(y'-y)=0\in E^{p_3}_{p_2}$.
From above it follows that $y'-y$ and hence also $y'$ is the image of some element in $E^{p_2}_{p_1}$, say $y'=\ell(x')$.
Then $k(x')=k(y')=k(j(c))=0\in D_{p_1}$, and $x'=j(b)$ for some $b\in D_{p_2}$.
Let $c':=c-i(b)\in D_{p_3}$.
Then $j(c')=0\in E^{p_3}_{p_1}$, hence $j(c')=0\in E^{p_3}_{p_2}$.
Therefore $z=j(c)=j(c')-j(i(b))=j(c')=0$.
\end{proof}

We define naturally associated $S$-terms for all $(b,p,q,z)\in I_4$ by
\begin{equation}
\label{eqStermsOfExactCoupleAnalog}
S^{pz}_{bq}:=\frac{\ker(d_{pqz}:E^p_q\to E^q_z)}{\im(d_{bpq}:E^b_p\to E^p_q)}.
\end{equation}
If $E$ is the exact couple system of an $I$-filtered chain complex $C$, then both definitions for $S$-terms, \eqref{eqSterms}
and \eqref{eqStermsOfExactCoupleAnalog}, coincide for all $(b,p,q,z)\in I_4$.

By Axiom~\ref{itECSaxiomTriangle}, $E^p_p=0$ for all $p\in I$, and hence
\[
S^{pq}_{pq} = E^p_q.
\]
We call the collection of these $S$-terms the \emph{$1$-page}.

For all $(b,p,q,z)\leq (b',p',q',z')$ in $I_4$, $\ell^{pq}_{p'q'}$ induces maps
\[
S^{pz}_{bq}\to S^{p'z'}_{b'q'}.
\]
For a proof, simply chase the diagram
\[
\xymatrix{
E^b_p \ar[r]^d \ar[d]^\ell & E^p_q \ar[r]^d \ar[d]^\ell & E^q_z \ar[d]^\ell \\
E^{b'}_{p'} \ar[r]^d       & E^{p'}_{q'} \ar[r]^d       & E^{q'}_{z'}.
}
\]
which is commutative by Axiom~\ref{itECSaxiomSquare} and functoriality of $E$.
We say that these maps between $S$-terms are \emph{induced by inclusions}.

\begin{lemma}[Extensions]
\label{lemExtensionProperty}
For any $z\leq p_1\leq p_2\leq p_3\leq b$ in $I$, we have a short exact sequence of maps induced by inclusion,
\begin{equation}
\label{eqExtensionProperty}
0\to S^{p_2,z}_{b,p_1}\to S^{p_3,z}_{b,p_1}\to S^{p_3,z}_{b,p_2}\to 0.
\end{equation}
\end{lemma}

\begin{proof}
The exactness can be proved using the following diagram.
\[
\xymatrix{&\\
E^b_{p_2} \ar[r]^d \ar[dd]_\ell    & E^{p_2}_{p_1} \ar[r]^d \ar[d]_\ell   & E^{p_1}_z \ar[dd]^\ell & \ar@/_2pc/[ll]\\
                                   & E^{p_3}_{p_1} \ar[ru]_d \ar[d]^\ell  & \\
E^b_{p_3} \ar[ru]^d \ar[r]_d       & E^{p_3}_{p_2} \ar[r]_d \ar[uul]_\ell|(.333)\hole \ar@/_2pc/@{-}[rr]+0 & E^{p_2}_z \ar[uul]^\ell|(.666)\hole. & \ar@/_1.3pc/@{-}[uu]+0_d \\
&
}
\]
The diagonal and the two horizontal compositions are the defining maps for the $S$-terms in \eqref{eqExtensionProperty}.
The three directed $3$-cycles are exact triangles.
First we show injectivity of the map $S^{p_2,z}_{b,p_1}\to S^{p_3,z}_{b,p_1}$:

Let $x\in E^{p_2}_{p_1}$ with $d(x)=0$ and $\ell(x)=d(y)\in E^{p_3}_{p_1}$ for some $y\in E^b_{p3}$.
Then $d(y)$ is zero in $E^{p_3}_{p_2}$.
Thus there exists $z\in E^b_{p_2}$ with $\ell(z)=y$.
Hence $\ell(d(z)-x)=0\in E^{p_3}_{p1}$.
Thus there exists $a\in E^{p_3}_{p_2}$ with $d(a)=d(z)-x$.
Therefore $z':=z-\ell(a)\in E^b_{p_2}$ has the property that $d(z')=d(z)-\ell(d(a)) = d(z)-(d(z)-x) = x$.
This means that $x$ represents zero in $S^{p_2,z}_{b,p_1}$.

Surjectivity of $S^{p_3,z}_{b,p_1}\to S^{p_3,z}_{b,p_2}$ can be proved similarly.

It remains to prove exactness at $S^{p_3,z}_{b,p_1}$:
Lemma~\ref{lemExactTriangles} shows that the composition of the two maps \eqref{eqExtensionProperty} is zero.
On the other hand, let $x\in E^{p_3}_{p_1}$ with $d(x)=0\in E^{p_1}_z$ and $\ell(x)=d(y)\in E^{p_3}_{p_2}$ for some $y\in E^b_{p_3}$.
Then $x':=x-d(y)\in E^{p_3}_{p_1}$ represents the same element as $x$ in $S^{p_3,z}_{b,p_1}$, and $\ell(x')=0\in E^{p_3}_{p_2}$.
Thus by Lemma~\ref{lemExactTriangles}, there exists $z\in E^{p_2}_{p_1}$ with $\ell(z)=x'$.
We have $d(z)=d(x')=d(d(y))=0\in E^{p_1}_z$.
Therefore $z$ represents an element in $S^{p_2,z}_{b,p_1}$ that maps to the element in $S^{p_3,z}_{b,p_1}$ that $x$ represents.
\end{proof}

As for the spectral system of an $I$-filtration we have differentials between $S$-terms of an exact couple system with the same properties:

\begin{lemma}[Differentials]
For any $(b,p,q,z),(b',p',q',z')\in I_4$ with $z\leq p'$ and $q\leq b'$ there are natural differentials
\begin{equation}
\label{eqDifferential}
d: S^{pz}_{bq}\to S^{p'z'}_{b'q'}
\end{equation}
that commute with $\ell$, that is, $\ell\circ d = d\circ \ell$.
\end{lemma}

\begin{proof}
Chasing the commutative diagram
\begin{equation}
\label{eqConstructionOfDifferential}
\xymatrix{
& E^q_z \\
E^p_q \ar[r]_k \ar[ur]^d                  & D_q \ar[u]_j               &                                & D_{q'} \ar[r]^j         & E^{q'}_{z'}\\
\ker(d_{pqz}) \ar@{^{(}->}[u] \ar@{.>}[r] & D_z \ar[u]^i \ar[r]^i      & D_{p'} \ar[r]^j \ar[ur]^0      & E^{p'}_{q'} \ar[u]^k \ar[ur]_d\\
E^b_p \ar[u]^d \ar@{.>}[r]                & E^q_z \ar[u]^k \ar[r]^\ell & E^{b'}_{p'} \ar[u]^k \ar[ru]_d &
}
\end{equation}
shows that there is a natural and well-defined $d$.
Dotted arrows means that we can choose these maps element-wise such that the diagram commutes.
\end{proof}

\begin{lemma}[Kernels and cokernels]
For any $(b,p,q,z),(b',p',q',z')\in I_4$ with $z=p'$ and $q=b'$ we have
\[
\ker\left(d: S^{pz}_{bq}\to S^{p'z'}_{b'q'}\right) = S^{pq'}_{bq}
\]
and
\[
\coker\left(d: S^{pz}_{bq}\to S^{p'z'}_{b'q'}\right) = S^{p'z'}_{pq'}
\]
\end{lemma}

\begin{proof}
We give only the proof of the first statement, the second one is symmetric.
Let $x\in E^p_q$ represent an element $[x]\in S^{pz}_{bq}$.
$[x]$ lies in $S^{pq'}_{bq}$ if and only if $x\in\ker(d_{pqq'}:E^p_q\to E^q_{q'})$.
By Axiom~\ref{itECSaxiomTriangle}, this is if and only if $k(x)=i(y)\in D_q$ for some $y\in D_{q'}$.

\[
\label{eqConstructionOfTightDifferential}
\xymatrix{
                               & E^q_{q'} \ar[r]^\ell               & E^q_{q'} \\
E^p_q \ar[r]^k \ar[ur]^d       & D_q \ar[u]^j \ar[ur]^j \ar[r]^i    & D_p \\
D_{q'} \ar[r]^i \ar[ur]^i      & D_{p'} \ar[u]^i \ar[r]^j \ar[ur]^i & E^{p'}_{q'} \\
E^p_{p'} \ar[r]^\ell \ar[ur]^k & E^q_{p'} \ar[u]^k \ar[ru]_d        &
}
\]

In this case, the construction \eqref{eqConstructionOfDifferential} of the differential \eqref{eqDifferential} and the triviality of the composition $D_{q'}\to D_{p'}\to E^{p'}_{q'}$ show that \eqref{eqDifferential} sends $[x]$ to zero.

Conversely, assume that $d$ from~\eqref{eqDifferential} sends $[x]$ to zero.
Let $z$ be the choice of the element in $D_{p'}=D_z$ that we made in diagram \eqref{eqConstructionOfDifferential} in order to construct $d([x])$.
Since $d([x])=0$, there is an $r\in E^{b'}_{p'}=E^q_{p'}$ such that $j(k(r))=j(z)\in E^{p'}_{q'}$.
Since $i(k(r))=0\in D_q$, we could have equally well chosen $z':=z-k(r)$ instead of $z$.
Note that $j(z')=0\in E^{p'}_{q'}$.
Thus $z'$ has a preimage $y\in D_{q'}$.
Therefore $k(x)=i(z)=i(z')=i(i(y))\in D_q$.
\end{proof}

\begin{lemma}[$\infty$-page as filtration quotients]
\label{lemExactCouplesInfinityPageFiltrationQuotients}
$D_\infty$ can be $I$-filtered by
\[
G_p:=\im (i_{p,\infty}:D_p\to D_\infty) \iso S^{p,-\infty}_{\infty,-\infty},\ \ \ p\in I.
\]
Furthermore the $S$-terms on the $\infty$-page are filtration quotients
\[
S^{p,-\infty}_{\infty,q}\iso G_p/G_q.
\]
\end{lemma}

\begin{proof}
By definition and Axiom~\ref{itECSaxiomTriangle},
\[
S^{p,-\infty}_{\infty,-\infty}=\frac{\ker(d:D_p\to D_{-\infty})}{\im(d:E^\infty_p\to D_p)}=\frac{D_p}{\ker(i:D_p\to D_\infty)}\iso G_p.
\]
Furthermore, Lemma~\ref{lemExtensionProperty} shows that
$
0\to G_q\to G_p\to S^{p,-\infty}_{\infty,q}\to 0
$
is exact.
\end{proof}

\begin{lemma}[$\infty$-page as quotient kernels]
\label{lemExactCouplesInfinityPageQuotientKernels}
$D_\infty$ has quotients
\[
Q_p:=\frac{D_\infty}{\ker(j:D_\infty\to E^\infty_p)}\iso S^{\infty,-\infty}_{\infty,p},\ \ \ p\in I.
\]
Furthermore the $S$-terms on the $\infty$-page are quotient kernels
\[
S^{p,-\infty}_{\infty,q}\iso \ker(Q_q\to Q_p).
\]
\end{lemma}

\begin{proof}
By definition and Axiom~\ref{itECSaxiomTriangle},
\[
S^{\infty,-\infty}_{\infty,p} = \frac{\ker(d:E^\infty_p\to D_p)}{\im(d:E^\infty_\infty\to E^\infty_p)}
=\im(j:D_\infty\to E^\infty_p) \iso Q_p.
\]
Furthermore, Lemma~\ref{lemExtensionProperty} shows that
$
0\to S^{p,-\infty}_{\infty,q}\to Q_q\to Q_p\to 0
$
is exact.
\end{proof}

\subsection{Natural isomorphisms} \label{secExactCouplesNatIsos}


A lattice is complete if arbitrary meets (`$\cap$') and joins (`$\cup$') exist.
A \emph{closed set system} is a family of sets that is closed under taking arbitrary unions and intersections.
In particular, a closed set system is a complete distributive lattice.

\begin{definition}[Excision]

An exact couple system $E$ over a complete distributive lattice $I$ is called \textdef{excisive} if for all $a,b\in I$,
\[
E^a_{a\cap b} \tof{\ell} E^{a\cup b}_b
\]
is an isomorphism.
\end{definition}

\begin{example}
An exact couple system $E$ over a closed set system $I$ is excisive if and only if for all $(p,q)\leq (p',q')$ in $I_2$ with $p\wo q = p'\wo q'$ the map $E^p_q\tof{\ell} E^{p'}_{q'}$ is an isomorphism.

The exact couple system of a filtered space~\ref{exExactCoupleOfFilteredSpace} is excisive by excision of $h_*$ if $(X_i)_{i\in I}$ is a family of open subsets that is closed under taking arbitrary unions and intersections.

Also, the exact couple system of an $I$-filtered chain complex~\ref{exExactCoupleOfFilteredChainComplex} is excisive if the filtration is distributive and $I$ is complete.
\end{example}

\begin{lemma}[Natural isomorphisms 1]
\label{lemNaturalIsomorphisms}
In an excisive exact couple system $E$ over a closed set system~$I$, $S^{pz}_{bq}$ is uniquely determined up to natural isomorphism by $b\wo p$, $p\wo q$, and $q\wo z$.
\end{lemma}

\begin{proof}
Let $(b,p,q,z), (b',p',q',z')\in I_4$ with $b\wo p$=$b'\wo p'$, and so on.
Let $b'':=b\cup b'$, and so on.
Then the vertical maps in
\[
\xymatrix{
E^b_p \ar[r]^d \ar[d]^\ell & E^p_q \ar[r]^d \ar[d]^\ell & E^q_z \ar[d]^\ell \\
E^{b''}_{p''} \ar[r]^d     & E^{p''}_{q''} \ar[r]^d     & E^{q''}_{z''}
}
\]
are isomorphisms since $E$ is excisive.
Thus $S^{pz}_{bq}\tof{\iso} S^{p''z''}_{b''q''}$ and analogously $S^{p''z''}_{b''q''}\fromf{\iso} S^{p'z'}_{b'q'}$, both maps being induced by inclusion.
\end{proof}

\begin{lemma}[Splitting principle for 1-page]
\label{lemSplitting1page}
In an excisive exact couple system $E$, for any $a,b\in I$ we have a commutative triangle of natural isomorphisms
\begin{equation}
\label{eqSplitting1page}
\xymatrix{
                                                        & E^{a\cup b}_{a\cap b} \ar[dr]^{(\ell,\ell)} & \\
E^a_{a\cap b} \oplus E^b_{a\cap b} \ar[rr]^{\ell\oplus \ell} \ar[ur]^{\ell+\ell} & & E^{a\cup b}_b \oplus E^{a\cup b}_a.
}
\end{equation}
\end{lemma}

\begin{proof}
By Axiom~\ref{itECSaxiomTriangle}, both compositions $E^a_{a\cap b}\to E^{a\cup b}_{a\cap b} \to E^{a\cup b}_a$ and $E^b_{a\cap b}\to E^{a\cup b}_{a\cap b} \to E^{a\cup b}_b$ are zero, which implies that diagram \eqref{eqSplitting1page} commutes.
Since $E$ is excisive the bottom map is an isomorphism. 
We find a section $E^{a\cup b}_a\iso E^b_{a\cap b}\to E^{a\cup b}_{a\cap b}$, which implies that the exact triangle
\[
\xymatrix{E^a_{a\cap b} \ar[r]^\ell & E^{a\cup b}_{a\cap b} \ar[r]^\ell & E^{a\cup b}_a \ar@/^2pc/[ll]^d
}
\]
splits, and this splitting coincides with the top left map $\ell+\ell$ in~\eqref{eqSplitting1page}.
\end{proof}

\begin{lemma}[General splitting principle]
\label{lemSplittingSterms}
In an excisive exact couple system $E$, for any indices $z,q,p,b,x,y\in I$ with $x\cap y\subseteq z\subseteq q\subseteq p\subseteq b\subseteq x\cup y$ we have a commutative triangle of natural isomorphisms
\begin{equation}
\label{eqSplittingSterms}
\xymatrix{
                        & S^{pz}_{bq} \ar[dr]^{(\ell,\ell)} & \\
S^{p\cap x,z\cap x}_{b\cap x,q\cap x} \oplus S^{p\cap y,z\cap y}_{b\cap y,q\cap y} \ar[rr]^{\ell\oplus \ell} \ar[ur]^{\ell+\ell} & & S^{p\cup x,z\cup x}_{b\cup x,q\cup x} \oplus S^{p\cup y,z\cup y}_{b\cup y,q\cup y}.
}
\end{equation}
\end{lemma}

\begin{proof}
Note that $(p\cap x)\cup q=p\cap (x\cup q)$.
By excision, Lemma~\ref{lemSplitting1page}, and then again excision, we have
\[
E^{p\cap x}_{q\cap x} \oplus E^{p\cap y}_{q\cap y}
\tof{\iso} E^{p\cap x\cup q}_q \oplus E^{p\cap y\cup q}_q
\tof{\iso} E^p_q
\tof{\iso} E^p_{p\cap x\cup q} \oplus E^p_{p\cap x\cup q}
\tof{\iso} E^{p\cup x}_{q\cup x} \oplus E^{p\cup y}_{q\cup y}.
\]
We do the same with $(b,p)$ and $(q,z)$ in place of $(p,q)$ and put that together with maps induced by inclusion into a $3\times 5$-diagram, which proves that $\ell+\ell$ and $(\ell,\ell)$ in~\eqref{eqSplittingSterms} are isomorphisms by definition of the $S$-terms.
Commutativity of the diagram follows from commutativity of~\eqref{eqSplitting1page}.
\end{proof}

Let us now suppose that $I=D(J)$ is the complete distributive lattice of downsets of some arbitrary poset $J$. 
As in Section~\ref{secReducedSterms}, we think of $J$ as an undirected graph, whose vertices are the elements of $J$, and $x,y\in J$ are adjacent if they are related, i.e. $x>y$ or $x<y$.
For $(b,p,q,z)\in I_4$,
let $Z(z,q,p,b)\subseteq J$ denote the union of all connected components of $p\wo z$ that intersect $p\wo q$, and let $B(z,q,p,b)\subseteq I$ denote the union of all connected components of $b\wo q$ that intersect $p\wo q$.

\begin{lemma}[Natural isomorphisms 2]
\label{lemNaturalIsomorphisms2}
In an excisive exact couple system $E$ over $I=D(J)$, $S^{pz}_{bq}$ is uniquely determined up to natural isomorphism by $Z:=Z(z,q,p,b)$ and $B:=B(z,q,p,b)$.
\end{lemma}

\begin{proof}
We have $p\wo q = Z\cap B$.
First we prove that we can change $z$ arbitrarily without changing $S^{pz}_{bq}$ as long as $Z$ stays invariant:
Let $z^*:=p\wo Z$ and $z_*:=\down(Z)\wo Z$ be the maximal and minimal elements less or equal to $q$ such that $Z(z_*,q,p,b)=Z=Z(z^*,q,p,b)$ (here we need the completeness of $I$).
Then $z_*\leq z\leq z^*$ implies that the two maps induced by inclusion $S^{pz_*}_{bq}\incl S^{pz}_{bq}\incl S^{pz^*}_{bq}$ are injections.
We claim that their composition is surjective.
\[
\xymatrix{
E^p_{z^*} \ar[r]^\ell \ar[d]_d & E^p_q \ar[r]^d \ar[d]^d & E^q_{z^*} \\
E^{z^*}_{z_*} \ar[r]^\ell & E^q_{z_*}.
}
\]
Let $x\in E^p_q$ with $d(x)=0\in E^q_{z^*}$.
Then there exists $x'\in E^p_{z^*}$ with $\ell(x')=x$.
In order to show that $d(x)=0\in E^q_{z_*}$, it suffices to show that the map $d:E^p_{z^*}\to E^{z^*}_{z_*}$ is zero.
Let $y:=\down(Z)$.
Then $p=z^*\cup y$ and $z_*=z^*\cap y$.
Hence Lemma~\ref{lemSplitting1page} implies that the two maps induced by inclusion $E^{z^*}_{z_*}\incl E^p_{z_*} \surj E^p_{z^*}$ are an injection followed by a surjection.
From the exact triangle for $(p,z^*,z_*)$ we deduce that $d:E^p_{z^*}\to E^{z^*}_{z_*}$ is indeed zero.

Similarly one shows that changing $b$ does not change $S^{pz}_{bq}$ as long as $B$ stays invariant.
Thus we may assume $z=p\wo Z$ and $b=q\cup B$.
The rest follows from Lemma~\ref{lemNaturalIsomorphisms}.
\end{proof}

\subsection{Connections} \label{secExactCoupleAnalogsConnections}

From what we proved so far about kernels, cokernels, and natural isomorphisms, it follows that both, lexicographic connections and the secondary connection from Section~\ref{secConnections}, apply to the spectral systems of excisive exact couple systems over $D(\Zinf^n)$ as well.

The only difference is that such exact couple systems do not give rise to $0$-pages as in Lemmas~\ref{lemLexConnection}(1), \ref{lemLexConnectionInSmallSteps}(1), and 
\eqref{eqSecConnS0}.
They start from the $1$-page
\[
S^{pz}_{bq}(P;1,\varphi,0) = 
S^{pz}_{bq}(P;0,\varphi,e_n) =
E^{A(P;\varphi)}_{A(P-e_n;\varphi)} \iso
E^{A(\varphi^{-1}(P))}_{A(\varphi^{-1}(P)-e_n)},
\]
and respectively
\[
S^{pz^*}_{b^*q}(P;0) = E^{A(P)}_{A(P-e_n)}.
\]

\begin{example}[Spectral system for $n$ open subsets]
\label{exNindependentSubspacesAndGeneralizedHomology}
Suppose $X$ is a topological space with $n$ open subsets $X_i$, and $h$ is a generalized homology theory. 
Then spectral system of $n$ independent subcomplexes from Section~\ref{secNindependentSubcomplexes} generalizes to this setting: It converges to $h(X)$, and its $S^{pz}_{bq}(P,1)$ terms for $P\in \{0,1\}^n$ are given by
\[
S^{pz}_{bq}(P,1) = h\big(X_P,\bigcup_{Q<P} X_Q\big) \iso h\big(\bigcap_{Q>P} X^Q, X^P\big),
\]
where $X_P:=\bigcap_{i:p_i=0} X_i$, $X^P:=\bigcup_{i:p_i=1} X_i$, $Q$ runs through $\{0,1\}^n$, and empty intersections mean~$X$.
If we denote the collection of all $S^{pz}_{bq}(P,k)$, $P\in\{0,1\}^n$, by $S^{pz}_{bq}(\{0,1\}^n,k)$, and similarly with $S^{pz^*}_{b^*q}(P,k)$, then the $2$-page is given by
\[
S^{pz^*}_{b^*q}(\{0,1\}^n,n) = H(\ldots H(S^{pz}_{bq}(\{0,1\}^n,1),d_1),\ldots,d_n),
\]
where $d_i$ denotes the natural differential in direction $-e_i$.
\end{example}

\begin{remark}[Spectral systems over {$D(\Zinf^n)^*$}]
\label{remSSoverDZnDual}
Let $E$ be an exact couple system over $D(\Zinf^n)^*$ (e.g. from Example~\ref{exExactCoupleOfFilteredSpace}).
The most natural analog of the lexicographic and secondary connection in the associated spectral system $S$ comes from regarding $E$ as an exact couple system over $D(\Zinf^n)$ using the identification of $D(\Zinf^n)^*$ with $D(\Zinf^n)$ via $p\mapsto -(\Zinf^n\wo p)$.
\end{remark}

\subsection{Multiplicative structure}

This section about products is not most general, but it is hoped to be sufficient for many situations where spectral systems appear.
Depending on the exact couple system, it may be useful to restrict $I$ to a subposet in order to define a reasonable multiplicative structure.

\begin{definition}
\label{defMultStructure}
A \textdef{multiplicative structure} on an exact couple system $E$ over a distributive lattice $I$ consists of a binary operation
\begin{equation}
\label{eqMultStructSumOnI}
+:I^2\to I
\end{equation}
and homomorphisms
\begin{equation}
\label{eqMultStructProdOn1Page}
\cup: E^{p_1}_{q_1} \otimes E^{p_2}_{q_2} \to E^{p_1+p_2}_{(p_1+q_2)\cup (q_1+p_2)}
\end{equation}
for any $(p_1,q_1),(p_2,q_2)\in I_2$ such that the following properties hold:
\begin{enumerate}
\item (Monotonicity of $+$) If $p\subseteq p'$ and $q\subseteq q'$ then $p+q\subseteq p'+q'$.
\item (Functoriality of $\cup$) If $(p_1,q_1)\leq (p_1',q_1')$ and $(p_2,q_2)\leq (p_2',q_2')$ in $I_2$ then the following diagram commutes:
\begin{equation}
\label{eqMultStructFunctoriality}
\xymatrix{
E^{p_1}_{q_1} \otimes E^{p_2}_{q_2} \ar[r]^{\cup\hspace{1pc}} \ar[d]_\ell & E^{p_1+p_2}_{(p_1+q_2)\cup (q_1+p_2)} \ar[d]^\ell \\
E^{p_1'}_{q_1'} \otimes E^{p_2'}_{q_2'} \ar[r]^{\cup\hspace{1pc}}         & E^{p_1'+p_2'}_{(p_1'+q_2')\cup (q_1'+p_2')}.
}
\end{equation}
\item (Leibniz rule) \label{itLeibnizRule}
For all $(p_1,q_1),(p_2,q_2)\in I_2$ the following diagram commutes:
\[
\xymatrix{
E^{p_1}_{q_1} \otimes E^{p_2}_{q_2} \ar[rr]^{\cup} \ar[d]_{d\otimes\ell \pm \ell\otimes d} & & E^{p_1+p_2}_{(p_1+q_2)\cup(q_1+p_2)} \ar[d]^d\\
D_{q_1} \otimes E^{p_2}_{q_2} \oplus E^{p_1}_{q_1} \otimes D_{q_2} \ar[rr]^{\hspace{0pc}\ell\circ\cup+\ell\circ\cup} & & E^{(p_1+q_2)\cup(q_1+p_2)}_{(q_1+q_2)\cup(p_1+\emptyset)\cup(\emptyset+p_2)}.
}
\]
\end{enumerate}
\end{definition}

By functoriality of `$\cup$' there is no danger in abbreviating any $\ell\circ\cup\circ\ell$ by $\cup$.
Let us omit the sign in the Leibniz rule, since this is not where the meat is in this paper; readers are invited to add and use their favorite sign convention.
Often in applications, $I$ is a closed system of subsets of a semigroup and the operation `$+$' is Minkowski sum, which then implies $p+\emptyset=\emptyset+p=\emptyset$.

\begin{example}[Filtered differential algebras]
Let $C$ be a differential algebra, distributively filtered by $(F_i)_{i\in I}$ over a complete distributive lattice $I$, which has a monotone binary operation `$+$' on it.
Suppose the product on $C$ sends $F_i\otimes F_j$ into $F_{i+j}$ for all $i,j\in I$.
Then this induces naturally a multiplicative structure on the associated exact couple system~\ref{exExactCoupleOfFilteredChainComplex}.
\end{example}

\begin{example}[Filtered spaces]
Let $X$ be filtered by a family of subsets $(X_i)_{i\in I}$ that is closed under taking unions and intersections.
If $h_*$ is a generalized homology theory, then for a multiplicative structure on the associated exact couple system of $(X_i)_{i\in I}$ we need products
\[
\cup: h_*(X_{p_1},X_{q_1}) \otimes h_*(X_{p_2},X_{q_2}) \to h_*(X_{p_1+p_2},X_{(p_1+q_2)\cup (q_1+p_2)})
\]
for some $+:I^2\to I$.
If $h^*$ is a generalized cohomology theory, then for a multiplicative structure on the associated exact couple system of $(X_i)_{i\in I}$ we need products
\[
\cup: h^*(X_{q_1},X_{p_1}) \otimes h^*(X_{q_2},X_{p_2}) \to h^*(X_{(p_1+q_2)\cup (q_1+p_2)},X_{p_1+p_2})
\]
for some $+:I^2\to I$.
These multiplicative structures only exist in suitable cases and then possibly only on (families of) sublattices $J\subset I$.
An example are the Leray--Serre spectral systems of Section~\ref{secSuccFibrations}, when $h^*$ is a generalized cohomology theory with products.

\end{example}

\begin{lemma}[Multiplication]
\label{lemMultiplication}
If $E$ is an exact couple with multiplicative structure, then there is a natural multiplication
\begin{equation}
\label{eqMultiplication}
S^{p_1,z_1}_{b_1,q_1} \otimes S^{p_2,z_2}_{b_2,q_2} \tof{\cup} S^{p',z'}_{b',q'}
\end{equation}
for all $(b_1,p_1,q_1,z_1),(b_2,p_2,q_2,z_2),(b',p',q',z')\in I_4$ that satisfy
\begin{align*}
b'  \supseteq\ & (p_1+b_2)\cup(b_1+p_2), \\
p'  \supseteq\ & p_1+p_2,\\
q'  \supseteq\ & (p_1+q_2)\cup(q_1+p_2)\cup(b_1+z_2)\cup(z_1+b_2), \textnormal{ and} \\
z'  \supseteq\ & (p_1+z_2)\cup(z_1+p_2).
\end{align*}
Moreover the product commutes with maps between $S$-terms induced by inclusion.
It also satisfies a Leibniz rule:
For all indices as above and analogs with a bar, such that $(q_j,z_j)\leq (\bar b_j,\bar p_j)$ for all $j\in\{1,2,\,'\}$, the following diagram commutes
\begin{equation}
\label{eqLeibnizRuleForS}
\xymatrix{
S^{p_1,z_1}_{b_1,q_1} \otimes S^{p_2,z_2}_{b_2,q_2} \ar[rrr]^{\hspace{2pc}\cup} \ar[d]_{d\otimes\ell \pm \ell\otimes d} &&& S^{p',z'}_{b',q'} \ar[d]^d\\
S^{\bar p_1,\bar z_1}_{\bar b_1,\bar q_1} \otimes S^{p_2,z_2}_{b_2,q_2} \oplus S^{p_1,z_1}_{b_1,q_1} \otimes S^{\bar p_2,\bar z_2}_{\bar b_2,\bar q_2} \ar[rrr]^{\hspace{5pc}\cup+\cup} &&& S^{\bar p',\bar z'}_{\bar b',\bar q'}
}
\end{equation}
whenever
\begin{align*}
\bar b'  \supseteq\ & (\bar p_1+b_2)\cup(\bar b_1+p_2)\cup(p_1+\bar b_2)\cup(b_1+\bar p_2), \\
\bar p'  \supseteq\ & (\bar p_1+p_2)\cup(p_1+\bar p_2),\\
\bar q'  \supseteq\ & (\bar p_1+q_2)\cup(\bar q_1+p_2)\cup(\bar b_1+z_2)\cup(\bar z_1+b_2)\, \cup \\
                    & (p_1+\bar q_2)\cup(q_1+\bar p_2)\cup(b_1+\bar z_2)\cup(z_1+\bar b_2), \textnormal{ and} \\
\bar z'  \supseteq\ & (\bar p_1+z_2)\cup(\bar z_1+p_2)\cup(p_1+\bar z_2)\cup(z_1+\bar p_2).
\end{align*}
\end{lemma}

\begin{proof}
The map $\cup:E^{p_1}_{q_1} \otimes E^{p_2}_{q_2} \to E^{p'}_{q'}$ induces \eqref{eqMultiplication}:
It sends pairs of cycles to cycles as we deduce from the commutative diagram
\[
\xymatrix{
E^{p_1}_{z_1} \otimes E^{p_2}_{z_2} \ar[r]^{\hspace{1.5pc}\cup} \ar[d]_{\ell\otimes\ell} & E^{p'}_{z'} \ar[d]^\ell \ar@/^2pc/[dd]^0\\
E^{p_1}_{q_1} \otimes E^{p_2}_{q_2} \ar[r]^{\hspace{1.5pc}\cup}      & E^{p'}_{q'} \ar[d]^d \\
& E^{q'}_{z'}
}
\]
and Lemma~\ref{lemExactTriangles}.
It is well-defined because of the commutative diagram
\[
\xymatrix{
E^{b_1}_{p_1} \otimes E^{p_2}_{z_2} \ar[rr]^{\hspace{1pc}\cup} \ar[d]_{d\otimes\ell \pm \ell\otimes d} && E^{b'}_{p'} \ar[d]^d\\
E^{p_1}_{q_1} \otimes E^{p_2}_{q_2} \oplus E^{b_1}_{p_1} \otimes D_{z_2} \ar[rr]^{\hspace{4pc}\cup+\cup} && E^{p'}_{q'}
}
\]
and its symmetric analog with indices $1$ and $2$ and the tensor order exchanged.
We need the assumptions on the indices $(b',p',q',z')$ in order to assure that the second cup on the bottom is the zero map and that the top map is well-defined.
Naturality and the Leibniz rule \eqref{eqLeibnizRuleForS} follow directly from the same properties of $E$.
\end{proof}

\begin{remark}[Compatibility with natural isomorphisms]
\label{remMultStrucCompatibilityWithNatIsos}
Suppose $E$ is an excisive exact couple system over a closed set system~$I$.

Further suppose that $E$ has a two \emph{compatible} multiplicative structures over sublattices $I'$ and $I''$ in the following sense:
If $(p_1',q_1'), (p_2',q_2')\in I'_2$ and $(p_1'',q_1''), (p_2'',q_2'')\in I''_2$ satisfy $p_1'\wo q_1'=p_1''\wo q_1''$ and $p_2'\wo q_2'=p_2''\wo q_2''$, then also $(p_1'+p_2')\wo((p_1'+q_2')\cup (p_2'+q_1')) = (p_1''+p_2'')\wo((p_1''+q_2'')\cup (p_2''+q_1''))$ and the following diagram with vertical maps being excision commutes,
\[
\xymatrix{
E^{p_1'}_{q_1'} \otimes E^{p_2'}_{q_2'} \ar[r]^{\cup'\hspace{1pc}} \ar[d]_\iso & E^{p_1'+p_2'}_{(p_1'+q_2')\cup (q_1'+p_2')} \ar[d]^\iso \\
E^{p_1''}_{q_1''} \otimes E^{p_2''}_{q_2''} \ar[r]^{\cup''\hspace{1pc}}         & E^{p_1''+p_2''}_{(p_1''+q_2'')\cup (q_1''+p_2'')}.
}
\]


Then also the corresponding products on $S$-terms~\eqref{eqMultiplication} for $I'$ and $I''$ commute with respect to the natural isomorphisms given by Lemma~\ref{lemNaturalIsomorphisms} (and Lemma~\ref{lemNaturalIsomorphisms2} in case $I=D(J)$ for some poset $J$), since~\eqref{eqMultiplication} is induced by~\eqref{eqMultStructProdOn1Page}.
%
%
\end{remark}

\begin{remark}[Compatibility with extensions]
\label{remMultStrucCompatibilityWithExtensions}
Consider any three extension sequences~\eqref{eqExtensionProperty} from Lemma~\ref{lemExtensionProperty}, which we write for short as $0\to A_i\to B_i\to C_i\to 0$ with $i\in\{1,2,'\}$, such that products $X_1\otimes X_2\to X'$ with $X\in\{A,B,C\}$ are defined via Lemma~\eqref{lemMultiplication}.
Then the diagram
\[
\xymatrix{
A_1\otimes A_2 \ar[r]^{\ell\otimes\ell}\ar[d]^\cup & B_1\otimes B_2 \ar[r]^{\ell\otimes\ell}\ar[d]^\cup & C_1\otimes C_2 \ar[d]^\cup \\
A' \ar[r]^\ell & B' \ar[r]^\ell & C'
}
\]
commutes, as the horizontal maps are induced by inclusion.
\end{remark}

\begin{example}[Several $\ZZ$-filtrations]
\label{exMultFornZfiltrations}
Suppose we are given $n$ different $\ZZ$-filtrations of a differential algebra $C$ as in \eqref{eqNZfiltrations}, such that $F_k^{(i)}\cdot F_\ell^{(i)}\subseteq F_{k+\ell}^{(i)}$ for all $i,k,\ell$, and such that the induced $D(\Zinf^n)$-filtration is distributive.
Fix a shearing matrix $\varphi\in \GL(n,\ZZ)$ with non-negative entries.
Let $I$ be the subposet of $D(\Zinf^n)$ consisting of all $A(P;\varphi)$, $P\in\Zinf^n$.
We define \eqref{eqMultStructSumOnI} as the Minkowski sum $A(P;\varphi)+A(P';\varphi):=A(P+P';\varphi)$, using the convention $\infty-\infty:=\infty$ in $\Zinf$ (the reason for this convention is similar to footnote \ref{footDefOfFpmInfty} on page \pageref{footDefOfFpmInfty}).
The product of $C$ induces a multiplicative structure \eqref{eqMultStructProdOn1Page} on the $1$-page.
Notice that in Lemma~\ref{lemMultiplication}, the unions of elements in $I$ are all trivial (that is, of the form $A\cup\ldots\cup A=A$) for all lexicographic connections and the secondary connection.
We obtain cup products
\[
S^{pz}_{bq}(p_1,\ldots,p_k;r,\varphi,Q)\otimes
S^{pz}_{bq}(p'_1,\ldots,p'_k;r,\varphi,Q) \to
S^{pz}_{bq}(p_1+p'_1,\ldots,p_k+p'_k;r,\varphi,Q)
\]
and
\[
S^{pz^*}_{b^*q}(P;k)\otimes 
S^{pz^*}_{b^*q}(P';k) \to
S^{pz^*}_{b^*q}(P+P';k),
\]
which satisfy the obvious Leibniz rules.
Also, the multiplicative structure is compatible with itself and is hence compatible with natural isomorphisms of $S$-terms (Remark~\ref{remMultStrucCompatibilityWithNatIsos}).
Thus these products induce the products on the subsequent pages of the secondary and lexicographic connections.
\end{example}

\begin{remark}[Weak multiplicative structure for $I=D(\Zinf^n)$]
\label{remMultStructureWeakForm}
For general exact couple systems over $D(\Zinf^n)$, a multiplicative structure as in Definition~\ref{defMultStructure} might be too much to ask for.
In order to have a cup product along the lexicographic connections for some fixed shearing matrix $\varphi$, it suffices to have a product~\eqref{eqMultStructProdOn1Page} only for pairs $(p_1,q_1)=(A(P_1;\varphi),A(Q_1;\varphi))$ and $(p_2,q_2)=(A(P_2;\varphi),A(Q_2;\varphi))$ with $P_1+Q_2=Q_1+P_2$.
Here the Leibniz rule needs to be slightly weakened accordingly by replacing $D_{q_1}$ with $E^{q_1}_{A(Q_1+Q_2-P_2;\varphi)}$ and $D_{q_2}$ with $E^{q_2}_{A(Q_1+Q_2-P_1;\varphi)}$.

In order to have a cup product along the secondary connection, we further require the products~\eqref{eqMultStructProdOn1Page} for successive $\varphi_k$ to be compatible; see Remark~\ref{remMultStrucCompatibilityWithNatIsos}.

We call these products along the secondary and lexicographic connections a weak multiplicative structure. 
\end{remark}

\begin{remark}[Cross product]
\label{exMultStructureCrossProduct}
Sometimes one has a more general cross product for three exact couple systems $E$, $E'$, and $E''$ over $I$,
\[
\cup: E^{p_1}_{q_1} \otimes E'{}^{p_2}_{q_2} \to E''{}^{p_1+p_2}_{(p_1+q_2)\cup (q_1+p_2)}.
\]
This section generalizes to this setting without difficulty.
\end{remark}

\section{Successive Leray--Serre spectral sequences} \label{secSuccFibrations}

Leray~\cite{Leray46LanneauDhomologie, Leray46StructureDeLanneauDhomologie, Leray46ProprietesDeLanneauDhomologie} and Serre~\cite{Serre51SS} constructed a spectral sequence that relates the homologies of the base, the fiber, and the total space of a fibration.
Here we study the situation of towers of fibrations.

\subsection{The spectral system}

Suppose we are given a tower of fibrations (always in the sense of Serre)
\begin{equation}
\label{eqSuccLSSSfibrationTower}
\xymatrix{
&& F_{i-1} \ar@{^{(}->}[r]^{i_{i-1}} & E_{i-1} \ar[d]^{f_i} &&&& \\
&&                                   & E_{i} 
\ar@{}[urrrr]|{\textnormal{
\begin{normalsize}
$(1\leq i\leq n)$,
\end{normalsize}}}
}  
\end{equation}
such that $E_1,\ldots,E_{n}$ have the homotopy type of a CW-complex.
We denote the tower with $E_*$.
Let's write $F_{n}:=E_{n}$, $E_{n+k}:=\pt$, $E_{-k}:=E_0$, and $F_{-k}:=\pt$ for $k\geq 1$.

\begin{theorem}
Let $h$ be a generalized homology theory.
Associated to the fibration tower \eqref{eqSuccLSSSfibrationTower} there is a spectral system over $I=D(\ZZ^{n})$ with $2$-page
\begin{equation}
\label{eqLSSSsecondPage}
S^{pz^*}_{b^*q}(P;n) = H_{p_n}(F_n;H_{p_{n-1}}(F_{n-1};\ldots H_{p_1}(F_{1};h(F_0)))))
\end{equation}
and limit $S^{\infty,-\infty}_{\infty,-\infty}=h(E_0)$.
\end{theorem}

Later we will use the shorter notation $H_\bullet(F_n;F_{n-1};\ldots;F_{1};h(F_0))$ for the direct sum of~\eqref{eqLSSSsecondPage} over all~$P$.

\begin{remark}[Local coefficients]
\label{remLSSSlocalCoefficients}
Of course in \eqref{eqLSSSsecondPage} we have local coefficients everywhere:
For $i\leq k$, we write $f_{i,k}:=f_{k}\circ\ldots\circ f_{i+1}\circ\id_{E_i}:E_i\to E_k$, and $F_{i,k}:=f_{i,k}^{-1}(\pt)$.
They also form fibrations $F_{i,k}\incl E_i\to E_k$, and more generally, $F_{i,k}\incl F_{i,\ell}\to F_{k,\ell}$ for $i\leq k\leq \ell$.

An element $\gamma\in\pi_1(F_k)$ induces a map $m^\gamma_{k-1,k}:F_{k-1,k}\to F_{k-1,k}$, and over it a map $m^\gamma_{k-2,k}:F_{k-2,k}\to F_{k-2,k}$, and over it a map $m^\gamma_{k-3,k}:F_{k-3,k}\to F_{k-3,k}$, and so on.
We regard them as a self-map $\wt m^\gamma_k$ on the fibration tower $F_i\incl F_{i,k}\to F_{i+1,k}$, $i<k$, and $\wt m^\gamma_k$ is uniquely given up to fiber homotopy.

For $k=1$, this makes $h(F_0)$ into a local coefficient system over $F_{1}$, as usual.
For $k=2$, $\gamma\in\pi_1(F_{2})$ induces a map $F_{1}\to F_{1}$ that respects the local coefficient system $h(F_0)$ over it; thus $H_*(F_{1};h(F_0))$ becomes a local coefficient system over $F_{2}$.
And so on.


\end{remark}

\begin{remark}[Naturality]
Any map $m:E_0\to E_0'$ between two such fibration towers (that is, it induces well-defined quotient maps $E_i\to E_i'$) naturally induces a
morphisms between the associated spectral systems.

\end{remark}

\begin{remark}[Cross product]
Suppose $h$ is a multiplicative generalized homology theory, that is, it comes from a ring spectrum.
The Cartesian product of two fibration towers $E_*$ and $E'_*$ is again a fibration tower $E''_*$.
As usual (compare with Switzer~\cite[p. 352--353]{Switzer75algTopHomotopyAndHomology}) the composition $h(X_p,X_q)\otimes h(Y_{p'},Y_{q'}) \tof{\times} h((X_p,X_q)\times (Y_{p'},Y_{q'})) \tof{i_*} h((X\times X)_{p+p'},(X\times Y)_{p+q'\cup p'+q})$ induces a cross product between the spectral systems of $E_*$ and $E'_*$ to the one of $E''_*$ for all pages in the lexicographic connections for arbitrary shearing matrix (as in Examples~\ref{exMultFornZfiltrations},~\ref{exMultStructureCrossProduct}), and it commutes with the differentials.
\end{remark}

\begin{proof}
Our following construction follows the idea of the Fadell--Hurewicz construction \cite{FadHur58structureOfHigherDiffInSS} of the Leray--Serre spectral sequence using singular prisms (very similar to Dress' construction~\cite{Dress67SSvonFaserungen}; also compare with Brown~\cite[Sect. 7]{Bro59twistedTensorProductsI}) and Dold's construction of the Leray--Serre spectral sequence for a generalized (co-)homology theory \cite{Dold63RelationsBetwOrdAndExtraordHomology}. 

For $P=(p_1,\ldots,p_{n})\in\ZZ^{n}$, let $\Delta_P:=\Delta_{p_1}\times\ldots\times\Delta_{p_{n}}$ be the product of $p_i$-dimensional simplices.
Let $\Delta_{P,*}$ denote the trivial fibration tower $\Delta_P\to \Delta_{P_{2\ldots n}} \to\ldots\to \Delta_{p_n}$.
Similarly, let $E^k_*$ denote the subtower $E_k\to\ldots\to E_n$ of \eqref{eqSuccLSSSfibrationTower}.
We define $K_P(E^1_*)$ as the set of all fibration tower maps $\Delta_{P,*}\to E^1_*$, that is, collections of maps $\Delta_{P_{i\ldots n}}\to E_i$ ($1\leq i\leq n$) that commute with the projections.
$K:=K_\bullet(E^{1}_*)$ forms in a natural way an $n$-fold simplicial set (even coalgebra), that is, a simplicial set in $n$ ways with commuting face and degeneracy maps, and it has a geometric realization $|K|$.
Let $X\to |K|$ denote the pullback of $f_1$ along the natural map $|K|\to E_{1}$,
\begin{equation}
\label{eqLSSSsimplicialModelForTowerE}
\xymatrix{
X \ar[d]\ar[r]      & E_0 \ar[d]^{f_1} \\
|K| \ar[r]          & E_{1}.
}
\end{equation}

We claim that $\tot(K)$ is chain homotopy equivalent to the singular chain complex $C_*(E_{1})$.
This can be proved similarly to the Eilenberg--Zilber theorem~\cite{EilZil53productsOfComplexes} using acyclic models, compare with \cite[Sect. 2.3.1 and 6.1]{FadHur58structureOfHigherDiffInSS} for the case $n=2$.
A chain homotopy equivalence $\alpha:\tot(K)\to C_*(E_{1})$ is obtained from the standard triangulation of $\Delta_P$.%
\footnote{$\Delta_k$ coincides with the order complex of the poset $\{0,\ldots,k\}$ with the usual order.
Thus we obtain a triangulation of $\Delta_P$ by taking the order complex of the product poset $\prod_i \{0,\ldots,p_i\}$, which is called the standard triangulation of $\Delta_P$.}

Alternatively, one can use inductively the argument in Dress~\cite[Sect. 2]{Dress67SSvonFaserungen} in order to show that $H_{p_{1}}(\ldots H_{p_n}(K,d_n)\ldots,d_{1})$ is zero except when $p_1=\ldots=p_{n-1}=0$, in which case it is naturally isomorphic $H_{p_n}(E_{1})$ via $\alpha$.
Thus, a spectral system argument for the $n$-complex $K$ shows that $\alpha$ is a quasi-isomorphism, and hence a homotopy equivalence since both complexes are free and bounded below.

Thus \eqref{eqLSSSsimplicialModelForTowerE} is a fiber homotopy equivalence and $h(X)=h(E_0)$.

Every $p\in I:=D(\Zinf^{n})$ indexes a skeleton of $K$ and thus a subspace $X_p\subseteq X$.
We define $S$ to be the spectral system of the $I$-filtered space $X$ with respect to $h$ as in Example~\ref{exExactCoupleOfFilteredSpace}.

The secondary connection for $S$ starts with abelian groups
\[
S^{pz}_{bq}(P;1) = h(X_{A(P)},X_{A(P)}) \iso K_P(E^1_*)\otimes h(F_0).
\]
Taking $k$ times homology in directions $-e_1,\ldots,-e_k$, we obtain with the usual arguments (compare with Dress~\cite[Sect. 3, 4]{Dress67SSvonFaserungen}, McCleary~\cite[Thm. 6.47]{McC01userGuideToSS})
\[
S^{pz^*}_{b^*q}(P;k) \iso K_{P_{k+1,\ldots,n}}(E^{k+1}_*) \otimes H_{p_k}(F_k;\ldots H_{p_1}(F_1;h(F_0)))))
\]
as abelian groups.
In particular for $k=n$ we get that the second page is naturally isomorphic to \eqref{eqLSSSsecondPage}.
%
%
%
%
\end{proof}

\subsection{Cohomology version}

As usual there is an analogous version for generalized cohomology theories $h$, using the exact couple system $E^p_q=h(X_q,X_p)$ over $I^*$, where $I:=D(\Zinf^n)$ (Example~\ref{exExactCoupleOfFilteredSpace}).
We identify $I^*\tof{\iso} I$ by $p\mapsto -(\Zinf^n\wo p)$ (Remark~\ref{remSSoverDZnDual}) in order to speak about secondary and lexicographic connections over $I^*$.
The second page is given by
\begin{equation}
\label{eqLSSScohomologyVersionSecondPage}
S^{pz^*}_{b^*q}(-P;n) = H^{p_n}(F_n;H^{p_{n-1}}(F_{n-1};\ldots H^{p_1}(F_{1};h(F_0)))))
\end{equation}
If $h$ is multiplicative, then 
the spectral system $S$ has indeed a natural product structure along the secondary and lexicographic connections as in Example~\ref{exMultFornZfiltrations} and Remark~\ref{remMultStructureWeakForm}.
For this we need to show that $E$ has a multiplicative structure with respect to the subposet of $I$ consisting of all $A(P;\varphi)$, $P\in\Zinf^n$, for certain fixed $\varphi$:

Let $(p_i,q_i)=(A(\varphi(-P_i);\varphi),A(\varphi(-Q_i);\varphi)) \in I_2$ for $i=1,2$, 
and define $M$ (possibly not uniquely) by the equation $A(\varphi(-M);\varphi) = A(\varphi(-Q_1-P_2);\varphi)\ \cup\ A(\varphi(-P_1-Q_2);\varphi)$.
Next, define the set $A^\circ_P:=-(\Zinf^n\wo A(\varphi(-P);\varphi) \in I^*$.
For $P\in\ZZ^n$, this equals $A(\varphi(P);\varphi)\wo\{-P\}$.
We define the sum~\eqref{eqMultStructSumOnI} in $I$ as the Minkowski sum in $I$, \emph{not} in $I^*$.
Thus the desired pairing \eqref{eqMultStructProdOn1Page} reads
\begin{equation}
\label{eqLSSScohomologyVersionPairingOnE1}
h(X_{A^\circ(Q_1)},X_{A^\circ(P_1)})\otimes h(X_{A^\circ(Q_2)},X_{A^\circ(P_2)}) \to
h(X_{A^\circ(M)},X_{A^\circ(P_1+P_2)}).
\end{equation}
Let $D:X\to X\times X$ denote a diagonal approximation (as the vertices of each face of $X$ are consistently ordered, a canonical choice is the Alexander--Whitney diagonal approximation).
The key is that $D$ restricts to a map
\begin{equation}
\label{eqLSSScohomologyVersionRestrictedDiagonalApproximation}
(X_{A^\circ(M)},X_{A^\circ(P_1+P_2)}) \to (X_{A^\circ(Q_1)}\times X_{A^\circ(Q_2)}\cup Y, Y),
\end{equation}
where $Y:=X\times X_{A^\circ(P_2)} \cup X_{A^\circ(P_1)}\times X$.
By excision,
\mathfootnotesize{
\begin{equation}
\label{eqLSSScohomologyVersionExcision}
h(X_{A^\circ(Q_1)}\times X_{A^\circ(Q_2)}\cup Y, Y) \tof{\iso}
h(X_{A^\circ(Q_1)}\times X_{A^\circ(Q_2)}, X_{A^\circ(Q_1)}\times X_{A^\circ(P_2)} \cup X_{A^\circ(P_1)}\times X_{A^\circ(Q_2)}).
\end{equation}
}%
We define the pairing~\eqref{eqLSSScohomologyVersionPairingOnE1} via the cross product of $h$ composed with the inverse of~\eqref{eqLSSScohomologyVersionExcision} and the map induced by~\eqref{eqLSSScohomologyVersionRestrictedDiagonalApproximation}.

Clearly the product of $S$ agrees in the limit $S^{\infty,-\infty}_{\infty,-\infty}=h(E_0)$ with the one in $h(E_0)$.
The multiplicative structures for different $\varphi$ are compatible, thus the products along the secondary connection are compatible (Remark~\ref{remMultStrucCompatibilityWithNatIsos}).
Therefore it agrees on the second page with the cup product of the right hand side of~\eqref{eqLSSScohomologyVersionSecondPage} (up to signs, depending on the convention).


\subsection{Edge homomorphisms}

Let $E_*$ be the fibration tower~\eqref{eqSuccLSSSfibrationTower}.
We use the notation from Remark~\ref{remLSSSlocalCoefficients}.
Clearly, for any $(b,a)\leq (d,c)\in\ZZ_2$, we have a natural map \begin{equation}
\label{eqLSSSnatMapBetweenGenFibers}
F_{a,b}\to F_{c,d}.
\end{equation}
If $b=d$ this map is a fibration with fiber $F_{a,c}$.

Let $x,y:\ZZ_{\geq 0}\to \ZZ_{\geq 0}$ be increasing maps with  $x(i)\leq i\leq y(i)$ for all $i\geq 0$, and $x(n+k)=x(n+1)$ for all $k\geq 1$
We define fibration towers $X_*$ and $Y_*$ by $X_i:=F_{x(i),x(n+1)}$ and $Y_i:=F_{y(i),y(n+1)}=E_{y(i)}$, $0\leq i\leq n$, using the maps \eqref{eqLSSSnatMapBetweenGenFibers} as the projections.
By definition, we get a composition of two fibration tower maps
\[
X_*\to E_*\to Y_*.
\]
It induces a map of spectral systems, which on the second page yields the maps
\mathsmall{
\[
H_\bullet(F_{x(n),x(n+1)};\ldots;h(F_{x(0),x(1)})) \to
H_\bullet(F_{n};\ldots;h(F_{0})) \to
H_\bullet(F_{y(n),y(n+1)};\ldots;h(F_{y(0),y(1)})),
\]
}%
and in the limit $h(F_{0,x(n+1)})\to h(E_0)\to h(E_{y(0)})$.

For a useful special case, choose $0\leq k\leq n$, set $x(i)=0$ and $y(i)=k$ for $i\leq k$, and set $x(i)=k+1$ and $y(i)=n+1$ for $i\geq k+1$.
Then $X_*$ and $Y_*$ have only one non-trivial fiber at $*=k$, and on the second pages we get
\mathsmall{
\[
H_\bullet(\pt;\ldots;F_{0,k+1};\ldots;h(\pt))\to
H_\bullet(F_{n};\ldots;F_k;\ldots;h(F_{0})) \to
H_\bullet(\pt;\ldots;E_k;\ldots;h(\pt)).
\]
}%
In case $h$ is ordinary homology with coefficients in some abelian group, or if $k=0$, then the spectral systems for $X_*$ and $Y_*$ clearly collapse at that second page.
For the limit we get $h(F_{0,k+1})\to h(E_0)\to h(E_{k})$.

There are more $x$ and $y$ for which $X_*$ and $Y_*$ have only one non-trivial fiber, however the corresponding maps to and from $E_*$ factor through the above ones.

For a generalized cohomology theory we get the analogous edge homomorphisms but in reversed direction.

\subsection{Rotating fibration towers}

In the introduction we saw that a vertical sequence of fibrations
\eqref{eqVerticalSeqOfFibrations}
induces a sequence of horizontal fibrations
\eqref{eqHorizontalSeqOfFibrations}.
Under certain conditions we can also go the other way around:

\begin{lemma}
If a horizontal sequence of fibrations up to homotopy \eqref{eqHorizontalSeqOfFibrations} can be delooped, then there is an associated vertical sequence of fibrations up to homotopy \eqref{eqVerticalSeqOfFibrations}.
\end{lemma}

\begin{proof}
Suppose $F=\Omega F'$, $P=\Omega P'$, $i_{FP}=\Omega i'_{FP}$, and so on.
Let $N$ be the homotopy fiber of $F'\incl E'$.
Then $F\incl E\to N$ is a fibration up to homotopy.
Let $N\to B=\Omega B'$ be the map that sends $(f,\gamma)$ to $p'_{EB}\circ\gamma$.
The preimage of the constant loop $\textnormal{const}_{b'_0}\in\Omega B'$ under this map is the homotopy fiber of $i'_{FP}$, which is homotopy equivalent to $\Omega(M)=M$.
%
%
%
\end{proof}


More generally we can consider any proper binary tree of fibrations, in which every node is the total space of a fibration whose left child is fiber and whose right child is the base space:
The root $E$ is the space whose homology we are interested in.
The leaves $F_1,\ldots, F_n$ are the fibers and base spaces whose homology we know.
At every non-leaf we can use a Leray--Serre spectral sequence that calculates the homology of that node.
Thus the tree gives us one way to compute $H_*(E)$ by successive spectral sequences from $H_*(F_1),\ldots,H_*(F_n)$.

The above operations for a diagram of two iterated fibrations can be applied also in the binary tree, giving left- and right rotations of the tree.
The number of such trees and hence ways to use successive spectral sequences (assuming we can always deloop) is the Catalan number $C_{n-1}=\binom{2n-2}{n-1}/n$.

\subsection{Successive Leray spectral sequences}
\label{secLSSSsuccLeraySS}

Another version for sheaves arises from the successive Grothendieck spectral sequences in Section~\ref{secSuccGrothendieckSS}.
Here we start with any sequence of maps $E_0\tof{f_0} \ldots\tof{f_{n-1}} E_{n-1}\tof{f_{n}} E_{n}:=\pt$ and a sheaf of abelian groups $G$ on $E_0$, for example the constant $\ZZ$ sheaf.
Let $\sheaves(X)$ denote the category of sheaves of abelian groups on a space $X$, and let $f_*:\sheaves(X)\to\sheaves(Y)$ denote the direct image functor for a map $f:X\to Y$.
Then we have a sequence of left-exact additive functors
\[
\xymatrix@+1pc {
\sheaves(E_0)\ar[r]^{\ \ (f_1)_*} & \ \ldots\ \ar[r]^{(f_{n-1})_*\ \ \ } & \sheaves(E_{n-1})\ar[r]^{\ \ (f_n)_*} & \sheaves(\pt).
}
\]
The associated spectral system over $D(\Zinf^n)$ converges to
\[
H^*(E_0;G) = R(f_n\circ\ldots\circ f_1)_*(G),
\]
and its $2$-page is given by
\[
H^*(E_{n-1}; R(f_{n-1})_*\circ\ldots\circ R(f_{1})_*(G)) = R(f_n)_*\circ\ldots\circ R(f_1)_*(G).
\]
Note that $R^k(f_i)_*(S)$ coincides with the sheaf associated to the presheaf $(U\subseteq E_{i}) \mapsto H^k(f_i^{-1}(U);S)$.



\section{Successive Grothendieck spectral sequences} \label{secSuccGrothendieckSS}

Let $\cat A_0$, $\cat A_1$, and $\cat A_2$ be abelian categories with enough projectives.
Then Grothendieck's spectral sequence~\cite{Gro57tohoku} computes the left derived functors of a composition of two 
right-exact functors $F_1:\cat A_0\to \cat A_1$ and $F_2:\cat A_1\to \cat A_2$ from the left derived functors of $F_1$ and $F_2$, assuming that $F_1$ sends projective objects to $F_2$-acyclic objects.
More precisely, the second page is given by $E^2_{pq}=L_pF_2\circ L_qF_1(A)$ and it converges to $L_{p+q}(F_2\circ F_1)(A)$, for any object $A\in \cat A_0$.

Now suppose we are given a sequence of $n$ 
right-exact functors
\begin{equation}
\label{eqCompositionOfFunctors}
\cat A_0\tof{F_1} \cat A_1\tof{F_2}\ldots \tof{F_n}\cat A_n
\end{equation}
between abelian categories with enough projectives.
We write $F_{ij}:=F_j\circ\ldots \circ F_{i+1}: \cat A_i\to \cat A_j$.
Assuming that for any $0\leq i<j<k\leq n$, $F_{ij}$ sends projective objects to $F_{jk}$-acyclic objects, we can relate $LF_1,\ldots,LF_n$ to $L(F_n\circ\ldots\circ F_1)$ by applying $n-1$ Grothendieck spectral sequences successively.
There are again $C_{n-1}=\binom{2n-2}{n-1}/n$ many ways to do that, since there are $C_{n-1}$ ways to bracket $F_n\circ\ldots\circ F_1$.

\subsection{Cartan--Eilenberg--Moore resolutions of $n$-complexes}

The construction of Grothendieck's spectral sequence is based on the Cartan--Eilen\-berg resolution for chain complexes, see Cartan--Eilen\-berg~\cite{CarEil56homologicalAlgebra}.
Here we will construct a higher-dimensional analog of Cartan--Eilenberg resolutions for $n$-complexes, which we coin \emph{Cartan--Eilenberg--Moore resolutions}, or \emph{CEM-resolutions} for short.

\paragraph{n-Complexes.}
Let $\cat A$ be an abelian category. 
Let $\chain(\cat A)$ denote the category of chain complexes with objects in $\cat A$, graded over $\ZZ$, the differentials being of degree $-1$.
For $(X,d)\in\chain(\cat A)$, let $C_k(X):=X_k$, $Z_k(X):=\ker(d:X_k\to X_{k-1})$, $Z'_k(X):=\coker(d:X_{k+1}\to X_k)$,  $B_k(X):=\im(d:X_{k+1}\to X_k)$, and $H_k(X):=Z_k(X)/B_k(X)$ denote the graded pieces, cycles, dual cycles, boundaries, and homology groups of $X$.



We define the category $\complex{n}{\cat A}$ of \textdef{$n$-complexes} inductively by $\complex{0}{\cat A}:=\cat A$ and $\complex{n}{\cat A}:=\chain(\complex{n-1}{\cat A})$.
Up to sign conventions, $\complex{2}{\cat A}$ is the category of double complexes over $\cat A$.
We denote the $n$ differentials of an $n$-complex by $d_1,\ldots,d_n$.
The order is important.
They satisfy $d_i\circ d_i=0$ and $d_i\circ d_j=d_j\circ d_i$ for all $i\neq j$.
We can regard an $n$-complex $X$ as chain complex over $\complex{n-1}{\cat A}$ in $n$ different ways, which we denote by $(X,d_i)$.
If we write only $X$, we always mean $(X,d_n)$.
The total complex $\tot(X)$ of an $n$-complex $X$ over $\cat A$ is the chain complex over $\cat A$ with $\tot(X)_k=\oplus_{|P|=k} X_P$, $|P|:=\sum p_i$, whose differential at the summand $X_P$ is $\sum (-1)^{p_1+\ldots+p_{i-1}}d_i$.

A homomorphism between two $n$-complexes $f:X\to Y$ is a $\ZZ^n$-graded homomorphism that commutes with the differentials.
A homotopy between two such homomorphisms $f,g:X\to Y$ is an $n$-tuple $(s_1,\ldots,s_n)$ of $\ZZ^n$-graded homomorphisms of degree $e_1,\ldots,e_n$, respectively, such that $f-g=\sum_{i=1}^n s_id_i-d_is_i$ and $s_id_j=d_js_i$ for all $i\neq j$.

\paragraph{Relative homological algebra.}
Let us fix the basic definitions of relative homological algebra.
(There are notational differences in the literature, see Eilenberg--Moore~\cite{EilMoo65relHomAlg} and Hilton--Stammbach~\cite{HiltonStammbach97CourseHomAlg}, but the concept is always the same.)
We will need this generality for the product structure in Section~\ref{secGSSproductStructure}.

Consider a class $\projClass E$ of epimorphisms in $\cat A$ which is closed under compositions and direct sums, and which contains all isomorphisms as well as all morphisms to $0$ (the standard choice is the class of all epimorphisms).
An object $P\in \cat A$ is called projective with respect to $X\to Y$, if any map $f:P\to Y$ factors over $X$.
$P$ is called $\projClass E$-projective, if it is projective with respect to all epimorphisms in $\projClass E$.
$\projClass E$ is called a projective class of epimorphisms, if for any object $K\in \cat A$ there exists an epimorphism $P\to K$ in $\projClass E$ with $P$ being $\projClass E$-projective.
$\projClass E$ is called closed, if any epimorphism in $\cat A$, such that all $\projClass E$-projective objects are projective with respect to it, lies already in $\projClass E$.

A morphism $f:X\to Y$ in $\cat A$ is called $\projClass E$-admissible if for the canonical factorization of $f$ into an epimorphism and a monomorphism, the epimorphism is in $\projClass E$.
An exact sequence $\ldots\to X_{n+1}\tof{f_n} X_{n}\to \ldots$ is called $\projClass E$-exact if all $f_n$ are $\projClass E$-admissible.
If furthermore $X_n=0$ for $n<0$ then it is called an $\projClass E$-acyclic resolution of $X_0$.
If furthermore $X_n$ are $\projClass E$-projective for all $n>0$ it is called an $\projClass E$-projective resolution of $X_0$.
$\projClass E$-projective resolutions exist for any projective class $\projClass E$ and have the usual properties.
We denote the derived functors of an additive functor $F:\cat A\to \cat B$ with respect to $\projClass E$-projective resolutions by $R_{\projClass E}F$.

\paragraph{CEM-resolutions.}
For all $P\in\ZZ^n$, $J\subseteq[n]$, we have a pair of adjoint functors $Q_{P,J}\dashv Z_{P,J}$:
Here, $Q_{P,J}:\cat A\to\complex{n}{\cat A}$ denotes the ``cube complex functor'' with $Q_{P,J}(K)$ having the object $K$ at all positions in $P+\{0\}^J\times\{0,-1\}^{[n]\wo J}$ and zeros otherwise, with isomorphisms between the $K$-entries.
$Z_{P,J}:\complex{n}{\cat A}\to \cat A$ denotes the functor $Z_{P,J}(X):=\{x\in X_P\st d_j(x)=0 \textnormal{ for } j\in J\}$.
Clearly, $\hom_{\complex{n}{\cat A}}(Q_{P,J}(K),X)\iso \hom_{\cat A}(K,Z_{P,J}(X))$.

Now, fix a closed projective class $\projClass E$ of epimorphisms in the abelian category $\cat A$.
We define $$\projClass E\textnormal{-CEM}_n$$ to be the class of all homomorphism of $n$-complexes $f:X\to Y$ such that the induced maps $Z_{P,J}(X)\to Z_{P,J}(Y)$ are $\projClass E$-surjective for all $P\in\ZZ^n$, $J\subseteq [n]$.
Below we usually omit $n$ from the notation because it will be clear from the context; also compare with Remark~\ref{remStabilizationOfCEMresolutions}.
From the Eilenberg--Moore multiple adjoint theorem \cite[Thm. 3.1]{EilMoo65relHomAlg} we deduce:

\begin{lemma}[$\projClass E$-CEM-projective objects]
\label{lemCEMprojObjects}
$\projClass E$-CEM is a closed projective class.
Furthermore, an $n$-complex $X\in\complex{n}{\cat A}$ is $\projClass E$-CEM-projective, if it is of the form $X=\bigoplus_{P\in\ZZ^n}\bigoplus_{J\subseteq[n]} Q_{P,J}(X_{P,J})$ with $X_{P,J}$ being $\projClass E$-projective objects of $\cat A$.
\end{lemma}

%

\begin{remark}[Canonical $\projClass E$-CEM-resolutions]
In case any $K\in \cat A$ admits a canonical $P\to K$ in $\projClass E$ with $P$ being $\projClass E$-projective, the same holds for $\projClass E$-CEM and there is a canonical $\projClass E$-CEM-resolution for any $X\in \complex{n}{\cat A}$.
\end{remark}


\begin{remark}[Alternative definition]
Equivalently we could inductively define projective classes of exact sequences $\projClass E_n:=\projClass E_n(\cat A)\subseteq\chain(\complex{n}{\cat A})$ as follows.
Let $\projClass E_0\subseteq \chain(\cat A)$ be the class of $\projClass E$-exact sequences in $\cat A$.
For $n\geq 1$, let $\projClass E_n\subseteq \chain(\complex{n}{\cat A})$ be the class of chain complexes $X\in\chain(\complex{n}{\cat A})$ such that $C_k(X,d_1)$ and $Z_k(X,d_1)$ lie in $\projClass E_{n-1}$.
(If $\projClass E$ is the class of all epimorphisms, then $B_k(X,d_1)$ and $H_k(X,d_1)$ will lie in $\projClass E_{n-1}$ as well; but not in general, as differentials might not be $\projClass E$-admissible.)
It follows from \cite[Thm. IV.2.1]{EilMoo65relHomAlg} that $\projClass E_n$ are indeed projective classes, and that the $\projClass E_n$-projective objects are exactly the $\projClass E$-CEM-projective $n$-complexes.
From their symmetry it follows that $\projClass E$-CEM-resolutions are symmetric: We can reorder the differentials of an $\projClass E$-CEM-resolution $Y\to X$ arbitrarily, and $Y\to X$ stays an $\projClass E$-CEM-resolution.
\end{remark}

\begin{lemma}[Maps and homotopies between $\projClass E$-CEM-resolutions]
\label{lemCEMresolutionMapsAndHomotopies}
Let $X$ and $Y$ be $n$-complexes, and let $P\to X$ be an $\projClass E$-CEM-projective resolution and $Q\to Y$ be a $\projClass E$-CEM-acyclic resolution.
Then any homomorphism $f:X\to Y$ admits an extension $P\to Q$ (as a homomorphism of $(n+1)$-complexes).
Any two such extensions are homotopic. 
More generally, if $F,G:P\to Q$ are extensions of two homotopic maps $f,g:X\to Y$, then $F$ and $G$ are homotopic.
\end{lemma}

Therefore, $\projClass E$-CEM resolutions can be regarded as a `resolution functor' $\hcomplex{n}{\cat A}\to\hcomplex{n+1}{\cat A}$ between the homotopy categories of $\complex{n}{\cat A}$ and $\complex{n+1}{\cat A}$, which is well-defined only up to natural isomorphisms.


%

\begin{remark}[Stabilization]
\label{remStabilizationOfCEMresolutions}
In Lemma ~\ref{lemHomologyOfSuccCEMResolution} and below we will regard an $n$-complex $X\in \complex{n}{\cat A}$ also as an $(n+k)$-complex $\wt X$ that is concentrated in the first $n$ coordinates.
Thus, $\complex{n}{\cat A}$ and $\hcomplex{n}{\cat A}$ are full subcategories of $\complex{n+k}{\cat A}$ and $\hcomplex{n+k}{\cat A}$, respectively.
In the same way, an $\projClass E$-CEM$_n$ resolution $Y\to X$ can be regarded as an $\projClass E$-CEM$_{n+k}$ resolution $\wt Y\to \wt X$ in $\chain(\complex{n+k}{\cat A})$.
Thus the resolution functor $\hcomplex{n}{\cat A}\to\hcomplex{n+1}{\cat A}\subseteq \hcomplex{n+k+1}{\cat A}$ coincides (again up to natural isomorphisms) with the restriction of the resolution functor $\hcomplex{n+k}{\cat A}\to\hcomplex{n+K+1}{\cat A}$ to $\hcomplex{n}{\cat A}$.

In general however, not every $\projClass E$-CEM resolution of $\wt X$ is concentrated in the first $n$ plus the resolution coordinates.
For example take a canonical resolution $Y\to X$ (if they exist for $\projClass E$), then $\wt Y\to \wt X$ is not the canonical resolution of $\wt X$, except if $X=0$.
In this sense, taking the canonical resolution is not a stable operation.
\end{remark}

\begin{remark}[Injective resolutions]
As usual one can dualize this and the next section and deal with injective classes of monomorphisms $\injClass M$ instead of $\projClass E$.
We only remark that the corresponding class $\injClass M$-CEM is then given by all maps of $n$-complexes $X\to Y$ such that the induced maps $Z'_{P,J}(X)\to Z'_{P,J}(Y)$ are $\injClass M$-injective for all $P\in\ZZ^n$, $J\subseteq [n]$, where $Z'_{P,J}(X):=X_P / \sum_{j\in J} \im(d_j: X_{P+e_j}\to X_P)$.
$\injClass M$-CEM is indeed a closed injective class since $Z'_{P-\one_{[n]\wo J},J}\dashv Q'_{P,J}$: $\Hom_{\cat A}(Z'_{P-\one_{[n]\wo J},J}(X),K)\iso \Hom_{\complex{n}{\cat A}}(X,Q_{P,J}(K))$, where $\one_{[n]\wo J} := \sum_{j\in[n]\wo J} e_j$.
\end{remark}

If $\projClass E$ is the class of all epimorphisms, then we will omit it from the notation.
In this case, if $X\to A$ is a CEM-resolution, then so is $H_k(X,d_i)\to H_k(A,d_i)$ (and similarly with $C_k$, $Z_k$, $B_k$).


\paragraph{Iterated resolutions.}
It will be convenient to have a symbol for iterated resolutions:
Let $D\in\complex{b_0}{\cat B_0}$, and let
\begin{equation}
\label{eqGSSfunctorSequenceForRes}
G_i:\cat B_{i-1}\to\complex{b_i}{\cat B_i}, \ \ \ 1\leq i\leq m,
\end{equation}
be functors.
We write
\[
\resolution(G_m,\ldots,G_1;D)\in\complex{m-1+\sum b_i}{\cat B_m}
\]
to be a complex that results from taking the CEM-resolution of $G_1(D)$, applying $G_2$, resolving it again, applying $G_3$, and so on, until we applied $G_m$.
Since the CEM-resolution depends on choices, so does $\resolution(G_m\ldots,G_1;D)$.
However a map $D\to D'$ of $b_0$-complexes induces a map between the chosen iterated resolutions $\resolution(G_m\ldots,G_1;D)\to \resolution(G_m\ldots,G_1;D')$, and homotopic maps induce homotopic maps.

$\resolution(G_m\ldots,G_1;D)$ has $m-1$ differentials that come from CEM-resolutions, we denote them by $d_{r_2},\ldots,d_{r_m}$.

\begin{lemma}[Homology of {$\resolution$}]
\label{lemHomologyOfSuccCEMResolution}
Taking homology of $\resolution(G_m\ldots,G_1;D)$ in direction $d_{r_i}$, for some $2\leq i\leq m$, yields
\[
H(\resolution(G_m\ldots,G_1;D),d_{r_i}) = \resolution(G_m,\ldots,G_{i+1},(LG_i)\circ G_{i-1},G_{i-2},\ldots,G_1;D).
\]
Here, $LG_i:\complex{b_{i-1}}{\cat B_{i-1}}\to\complex{b_{i-1}+b_i+1}{\cat B_i}$ denotes the left-derived functor of $G_i$ with respect to CEM-resolutions.

If $G_{i-1}$ sends projective objects to ($b_i$-complexes of) $G_i$-acyclic objects, then $(LG_i)\circ G_{i-1}$ simplifies to $G_i\circ G_{i-1}$, seen as a functor $\cat B_{i-2}\to\complex{b_{i-1}+b_i+1}{\cat B_i}$ that is concentrated where the last coordinate is zero.
\end{lemma}

\begin{proof}
We have
\[
\begin{aligned}
& H(\resolution(G_m,\ldots,G_{1};D),d_{r_i})
& = \\ 
& H\big(\resolution\big(G_m,\ldots,G_{i+1},\id;
\resolution(G_i,G_{i-1},\ldots,G_1;D)
\big),d_{r_i}\big)
& = \\
& \resolution\big(G_m,\ldots,G_{i+1},\id;
H(\resolution(G_i,G_{i-1},\ldots,G_1;D),d_{r_i})
\big)
& = \\
& \resolution\big(G_m,\ldots,G_{i+1},\id;
\resolution((LG_i)\circ G_{i-1},G_{i-2},\ldots,G_1;D)
\big)
& = \\
& \resolution(G_m,\ldots,G_{i+1},(LG_i)\circ G_{i-1},G_{i-2},\ldots,G_1;D).
\end{aligned}
\]
The second equality uses that CEM-resolutions commute with taking homology in one of the complex directions 
and that CEM-injectives are sums of cube complexes by Lemma~\ref{lemCEMprojObjects}.
\end{proof}

We can repeat this lemma for different directions $d_{r_i}$.
There are $(m-1)!$ ways to apply the lemma $m-1$ times successively.


\subsection{The spectral system}
Now suppose we are given a sequence of functors~\eqref{eqCompositionOfFunctors} as above, and let $X^0\in \cat A_0$.
We define
\[
X:= 
\resolution(F_n,\ldots,F_0,X^0),
\]
where $F_0:\cat A_0\to \cat A_0$ is the identity functor.

Let $C$ be the total complex of $X$.
$X$ and $C$ can be naturally $\ZZ$-filtered in $n$ different ways:
For $1\leq i\leq n$ and $p\in\ZZ$, let $F^{(i)}_p:=\bigoplus_{k\leq p}(X,d_i)_k$.
Thus from Section~\ref{secSeveralZfiltrations} we obtain a spectral system $S$.

\begin{lemma}[Limit]
If $F_i$ sends projective objects of $\cat A_{i-1}$ to $F_{i,n}$-acyclic objects for all $1\leq i\leq n-1$, then the limit of $S$ is
\begin{equation}
\label{eqGSSlimit}
S^{\infty,-\infty}_{\infty,-\infty} = H_*(C) = L_*F_{0,n}(X^0).
\end{equation}
\end{lemma}

\begin{proof}
We apply the secondary connection from Section~\ref{secSecondaryConnection} to $S$ with reversed coordinates, that is, we will first take differentials in direction $d_n$, then $d_{n-1}$, and so on.
Let us denote the corresponding $S$-terms from Section~\ref{secSecondaryConnection} by $\bar S^{pz}_{bq}(P;k)$.
We have $\bar S^{pz}_{bq}(P;0) = \bar S^{pz}_{bq}(P;1) = X_P = \resolution(F_n,\ldots,F_0,X^0)_P$ for $P\in\ZZ^n$.
Note that the graded pieces of $\projClass E_n(\cat A_i)$-projective objects are projective objects in $\cat A_i$, and $F_{n-1}$ sends projective objects of $A_{n-2}$ to $F_n$-acyclic objects in $\cat A_{n-1}$.
Thus taking differentials in direction $d_n$ and using Lemmas~\ref{lemSecConnNatIsos}, \ref{lemSecConnDifferentials}, and~\ref{lemHomologyOfSuccCEMResolution} yields
\[
\bar S^{pz}_{bq}(P;2) = \resolution(F_n\circ F_{n-1},F_{n-2},\ldots,F_0;X^0)_P.
\]
Iterating this we get
\[
\bar S^{pz}_{bq}(P;k) = \resolution(F_{n-k,n},F_{n-k},\ldots,F_0;X^0)_P
\]
for $k=1,\ldots,n$, using that $F_{n-k}$ sends projective objects of $\cat A_{n-k-1}$ to $F_{n-k,n}$-acyclic objects.
Taking homology with respect to the last differential in direction $d_1$ yields 
\[
\bar S^{pz^*}_{b^*q}(P;n) = \resolution((LF_{0,n})\circ F_0; X^0),
\]
which at positions $P\in\ZZ_{\geq 0}\times \{0\}^{n-1}$ is $L_P F_{0,n}(X^0)$, and zero otherwise.
These non-zero $S$-terms are all homogeneous in different total degree given by $p_1$.
Thus when we proceed with the lexicographic connection, the result~\eqref{eqGSSlimit} follows.
There is no problem with limits, since $X$ has only finitely many non-zero graded pieces in every degree.
\end{proof}

\begin{lemma}[Second page]
\label{lemGSSsecondPage}
The second page \eqref{eqSecConnS*n} of $S$ at position $P$ is
\begin{equation}
\label{eqGSSsecondPage}
S^{pz^*}_{b^*q}(P;n) = 
(L_{p_n}F_n)\circ\ldots\circ (L_{p_1}F_1)(X^0).
\end{equation}
\end{lemma}

\begin{proof}
The secondary connection starts with
\[
S^{pz}_{bq}(P;0) = S^{pz}_{bq}(P;1) = X_P = \resolution(F_n,\ldots,F_0,X^0)_P.
\]
By Lemma~\ref{lemHomologyOfSuccCEMResolution}, taking homology in direction $d_1$ yields
\begin{equation}
\label{eqGSSsecondaryConnectionFirstStep}
S^{pz^*}_{b^*q}(P;1) = \resolution(F_n,\ldots,F_2,LF_1;X^0)_P.
\end{equation}
Iterating this by taking homology in directions $d_2$, $d_3$, and so on, we obtain
\[
S^{pz^*}_{b^*q}(P;k) = \resolution(F_n,\ldots,F_{k+1},(L_*F_{k})\circ\ldots\circ(L_*F_1);X^0)_P,
\] 
for $k=1,\ldots,n$.
For $k=n$ we obtain the second page.
\end{proof}

There are many more second pages given by permuting the coordinates of $X$.

\begin{remark}[Naturality]
Lemma~\ref{lemHomologyOfSuccCEMResolution} yields explicit formulas for these $n!$ second pages in terms of left-derived functors only.
Thus the second pages do not depend anymore on the particular choice of the CEM-resolutions that we made in order to define $X$.


Any morphism $f_0:X^0\to (X^0)'$ in $\cat A_0$ induces a well-defined and natural morphism between second pages~\eqref{eqGSSsecondPage} of the spectral systems for $X^0$ and $(X^0)'$ by Lemma~\ref{lemCEMresolutionMapsAndHomotopies}.
Hence all $S$-terms that follow the respective second page by the lexicographic connections are well-defined, and morphisms $f_0$ induce natural morphisms between these $S$-terms.
\end{remark}


\subsection{Product structure}
\label{secGSSproductStructure}

Constructing a product structure for the Grothendieck spectral sequence is not a completely trivial matter, even if the second page has a natural one.

In this section we build upon the work of Swan~\cite{Swan99cupProductsInSheafCohomology}.
He constructs a product structure on the Leray--spectral sequence $E_2^{pq}=H^p(Y;R^qf_*(S))\impl H^{p+q}(X;S)$ (where $f:X\to Y$ is continuous and $S$ is a sheaf of rings on $X$), by regarding it as a hypercohomology spectral sequence for $\ext(\underline\ZZ,f_* I^*)$, where $I^*$ is an injective resolution of $S$ and $\underline\ZZ$ is the constant $\ZZ$-sheaf.
He obtains a product structure using a final system of flat resolutions of $\underline\ZZ$.
At the end of~\cite[Sect. 3]{Swan99cupProductsInSheafCohomology}, he remarks that a suitable pure injective Cartan--Eilenberg resolution of $f_*(I)$ would suffice as well if one could construct it. That is exactly what we will do. (Note that in an earlier paper \cite{Swan96hochschildCohomology} Swan defined a version of pure injective CE-resolutions, which differs from ours for $n=2$.)

\medskip

For this section, we consider the sequence of additive functors \eqref{eqCompositionOfFunctors}, this time with the following extra structure:
On the abelian categories $\cat A_0,\ldots,\cat A_n$ we have biadditive functors $\otimes:\cat A_i\times \cat A_i\to \cat A_i$ with natural isomorphisms $X\otimes Y\iso Y\otimes X$ and $(X\otimes Y)\otimes Z\iso X\otimes(Y\otimes Z)$ and a right adjoint biadditive functor $\hom:\cat A\times \cat A\to \cat A$ such that $\Hom(X\otimes Y,Z)\iso\Hom(X,\hom(Y,Z))$.
Furthermore we require that $\cat A_0,\ldots,\cat A_{n-1}$ (all but $\cat A_n$) are Baer--Grothendieck categories, that is, they have generators and filtered colimits exist and are exact.
By Grothendieck~\cite[Thm. 1.10.1]{Gro57tohoku}, $\cat A_0,\ldots,\cat A_{n-1}$ have enough injectives.
Further we require that the functors $F_i$ commute with $\otimes$, that is, there are natural maps
\begin{equation}
\label{eqGSSproductFiCommutesWithTensor}
F_i(X)\otimes F_i(Y)\to F_i(X\otimes Y), \ \ \ X,Y\in \cat A_{i-1}, \ 0\leq i\leq n-1.
\end{equation}
Note that $\otimes:\cat A_i\times \cat A_i \to \cat A_i$ extends naturally to a biadditive functor $\complex{k}{\cat A_i}\otimes\complex{k}{\cat A_i}\to \complex{k}{\cat A_i}$, which we also denote by $\otimes$.

\begin{remark}
More generally, we could start with three sequences of functors
$\cat A_0\tof{F_1}\ldots\tof{F_n}\cat A_n$, $\cat A_0\tof{F'_1}\ldots\tof{F'_n}\cat A_n$, and $\cat A_0\tof{F''_1}\ldots\tof{F''_n}\cat A_n$, together with natural maps 
$F_i(X)\otimes F'_i(Y)\to F''_i(X\otimes Y)$, $X,Y\in \cat A_{i-1}$, $0\leq i\leq n-1$.
As above, this yields a pairing between the first two spectral systems to the third one.
\end{remark}

A monomorphism $X\to Y$ in $\cat A_i$ is called pure if $X\otimes Z\to Y\otimes Z$ is a monomorphism for all $Z\in \cat A_i$.
Let $\injClass M_i$ be the class of all pure monomorphisms in $\cat A_i$.
From \cite[Thm. 2.1]{Swan99cupProductsInSheafCohomology} it follows that the $\injClass M_i$ are injective classes for $0\leq i\leq n-1$, that is, $\cat A_i$ has enough pure injectives.

Let $X,Y\in\complex{k}{\cat A_i}$ for some $0\leq i\leq n-1$.
We regard $X\otimes Y$ also as a $k$-complex over $\cat A_i$.
Choose $\injClass M_i$-CEM-injective resolutions $X\to I$, $Y\to J$ and $X\otimes Y\to K$.
We regard $I\otimes J$ as a sequence of $k$-complexes (these $k$-complexes are in general not injective anymore), together with a map $X\otimes Y\to I\otimes J$.
As in \cite[Sect. 2, 3]{Swan99cupProductsInSheafCohomology} one obtains:

\begin{lemma}
$X\otimes Y\to I\otimes J$ is an $\injClass M_i$-CEM-acyclic resolution. 
\end{lemma}

Thus we can choose a map $I\otimes J\to K$ over $\id_{X\otimes Y}$.
Applying $F_{i+1}$ and precomposing with \eqref{eqGSSproductFiCommutesWithTensor} yields $F(I)\otimes F(J)\to F(I\otimes J)\to F(K)$.
Taking homology yields 
\[
(R_{\injClass M_i\textnormal{-CEM}}F_{i+1})(X)\otimes (R_{\injClass M_i\textnormal{-CEM}}F_{i+1})(Y)\to (R_{\injClass M_i\textnormal{-CEM}}F_{i+1})(X\otimes Y).
\]
If furthermore we are given a map of $k$-complexes $X\otimes Y\to Z$, we can take an $\injClass M_i$-CEM-resolution $Z\to L$, choose an extension $K\to L$, and get $F_{i+1}(K)\to F_{i+1}(L)$ over $F_{i+1}(X\otimes Y)\to F_{i+1}(Z)$.
This induces a cup product
\[
(R_{\injClass M_i\textnormal{-CEM}}F_{i+1})(X)\otimes (R_{\injClass M_i\textnormal{-CEM}}F_{i+1})(Y)\to (R_{\injClass M_i\textnormal{-CEM}}F_{i+1})(Z).
\]
By Lemma~\ref{lemCEMresolutionMapsAndHomotopies} it is well-defined, and homotopic maps $X\otimes Y\to Z$ define the same cup product.

More generally, suppose we are given a sequence of functors \eqref{eqGSSfunctorSequenceForRes} with the analogous extra structure: $\cat B_0,\ldots,\cat B_{m-1}$ are Baer--Grothendieck categories, $\cat B_0,\ldots,\cat B_m$ have $\otimes$ and $\hom$, and the functors $G_i$ preserve $\otimes$.

Suppose we are given $X,Y,Z\in\complex{b_0}{\cat B_0}$ and a map $X\otimes Y\to Z$ of $b_0$-complexes (which might be given only up to homotopy).
Iterating the above procedure, we obtain a map of $(m-1+\sum b_i)$-complexes,
\begin{equation}
\label{eqGSScupProductForPureRes}
\resolution_{\injClass M}(G_m,\ldots,G_1;X)\otimes \resolution_{\injClass M}(G_m,\ldots,G_1;Y)\to \resolution_{\injClass M}(G_m,\ldots,G_1;Z),
\end{equation}
which is well-defined up to homotopy, where the subscript $\injClass M$ means that we always take $\injClass M_i$-CEM-injective resolutions.
Therefore we obtain a multiplicative structure on the associated spectral systems over $I=D(\Zinf^n)$ as in Example~\ref{exMultFornZfiltrations} (note that here the $n$-complexes are concentrated in the negative quadrant).
We still need to discuss when this induces a multiplication on the spectral system $S$ from the previous section.


If $\injClass M\subseteq \injClass M'$ are two injective classes of monomorphisms in an abelian category $\cat A$, then $\injClass M$-injective resolutions are $\injClass M'$-acyclic resolutions.
Thus if $X\to I$ and $X\to I'$ are an $\injClass M$- and an $\injClass M'$-injective resolution of $X$, we obtain a map $I\to I'$ over $\id_X$, which is natural up to homotopy.
Iterating this, we obtain a map
\begin{equation}
\label{eqGSSprodMapsBetweenResMandRes}
\resolution_M(G_m,\ldots,G_1;X) \to \resolution (G_m,\ldots,G_1;X),
\end{equation}
and analogous maps for $Y$ and $Z$.
We want to compose them with \eqref{eqGSScupProductForPureRes} in order to obtain a product structure on the spectral system from the previous section, however the maps for $X$ and $Y$ go in the wrong direction.

Let $S$ denote (as in the previous section) the spectral system for $\resolution (G_m,\ldots,G_1;X)$, and $S_{\injClass M}$ the spectral system for $\resolution_{\injClass M}(G_m,\ldots,G_1;X)$.
Then \eqref{eqGSSprodMapsBetweenResMandRes} induces maps $S_{\injClass M}\to S$.

\begin{lemma}
\label{lemGSSisoBetweenSMandS}
Suppose that all pure injective objects in $\cat A_{i-1}$ are $F_i$-acyclic for all $1\leq i\leq n$.
Then \eqref{eqGSSprodMapsBetweenResMandRes} induces an isomorphism between the second pages of $S_{\injClass M}$ and $S$, both of which are thus given by \eqref{eqGSSsecondPage}.
\end{lemma}

\begin{proof}
The proof idea of Lemma~\ref{lemGSSsecondPage} works here as well using the following extra facts:
Any $\injClass M_{i-1}$-CEM-injective resolutions $X\to D$ of any complex $D\in\complex{n}{\cat A_{i-1}}$ is a CEM-acyclic resolution of $D$ with $F_i$-acyclic entries.
Let $1\leq j\leq n$.
Using the exact sequences $0\to Z_k(\blank,d_j)\to C_k(\blank,d_j)\to B_{k-1}(\blank,d_j)\to 0$ and $0\to B_k(\blank,d_j)\to Z_k(\blank,d_j)\to H_k(\blank,d_j)\to 0$, we deduce using induced long exact sequences (as with usual Cartan--Eilenberg resolutions) that $H(X,d_j)\to H(D,d_j)$ is also a CEM-acyclic resolution, whose entries are again $F_i$-acyclic by the special form of $\injClass M_i$-injectives.
And as usual, $LF_i$ can be constructed from $F_i$-acyclic resolutions.
\end{proof}

\begin{corollary}
If all pure injective objects in $\cat A_{i-1}$ are $F_i$-acyclic for all $1\leq i\leq n$, then the product structure on $S_{\injClass M}$ induces a product structure on $S$ on the second page and all pages following it via the lexicographic connections.
\end{corollary}

\begin{example}[Successive Leray spectral sequences]
By Swan~\cite[Lemma 3.6]{Swan99cupProductsInSheafCohomology}, in $\sheaves(X)$ all pure injective sheaves are flasque, and flasque sheaves are $f_*$-acyclic for any $f:X\to Y$.
Thus, Lemma~\ref{lemGSSisoBetweenSMandS} implies that the Leray spectral system from Section~\ref{secLSSSsuccLeraySS} has a product structure from the second page on.
\end{example}


\section{Adams--Novikov and chromatic spectral sequences} \label{secAdamsAndChromaticSS}

The chromatic spectral sequence converges to the second page of the Adams--Novikov spectral sequence, which in turn converges to the $p$-component of the stable homotopy groups of spheres $\pi_*(S^0)$.
For background see Ravenel \cite{Rav86complexCobordismAndStableHomotopy}.
In this section we show how to unify these two spectral sequences into one spectral system over $D(\Zinf^2)$.


\medskip

Let $BP$ be the Brown--Peterson spectrum, $\Gamma:=BP_*(BP)$, and let $\ext(\blank):=\ext_\Gamma(BP_*,\blank)$.
Miller, Ravenel, and Wilson~\cite{MRW77periodicPhenomenaANSS} constructed the \emph{chromatic resolution}, which is a long exact sequence of $BP_*(BP)$-comodules,
\begin{equation}
\label{eqChromaticResolution}
0\to BP_*\to M^0\to M^1 \to\ldots
\end{equation}
Ravenel~\cite{Rav87geometricRealizationOfChromRes} realized the chromatic resolution~\eqref{eqChromaticResolution} geometrically in the following sense (note that our $N_n$ is Ravenel's $\Sigma^{-n}N_n$):
There are spectra $N_n$ and $M_n$, $n\geq 0$, such that
$M^n=BP_*(M_n)$, $N_0:=S^0$ the sphere spectrum, 
and there are fibrations $N_{n+1}\to N_n\to \Sigma^{-n}M_n$ that induce short exact sequences
$0\to BP_*(N_n) \to BP_*(\Sigma^{-n} M_n) \to BP_*(\Sigma N_{n-1}) \to 0$.
Thus this gives a diagram
\begin{equation}
\label{eqChromatricResolutionInFormOfAdamsResolution}
\xymatrix{
S^0=N_0 \ar[d] & N_1 \ar[d] \ar[l] & N_2 \ar[d] \ar[l] & \ldots \ar[l]\\
M_0            & \Sigma^{-1} M_1        & \Sigma^{-2} M_2.
}
\end{equation}
%
Let
\begin{equation}
\label{eqAdamsResolutionOfX}
\xymatrix{
S^0=X_0 \ar[d] & X_1 \ar[d] \ar[l] & X_2 \ar[d] \ar[l] & \ldots \ar[l] \\
I_0            & \Sigma^{-1} I_1        & \Sigma^{-2} I_2,
}
\end{equation}
be the canonical $BP$-Adams resolution of $S^0$.
That is,
\[
\Sigma^{-n}I_n:=BP\wedge X_n
\]
and $X_{n+1}$ is defined as the fiber of $X_n\to \Sigma^{-n}I_n$.
We obtain the canonical $BP$-Adams resolution of any spectrum $X$ by smashing \eqref{eqAdamsResolutionOfX} with $X$.
As usual we obtain short exact sequences
\begin{equation}
\label{eqANSS_SES}
0\to BP_*(X\wedge X_n) \to BP_*(X\wedge \Sigma^{-n}I_n) \to BP_*(X\wedge \Sigma X_{n+1})\to 0
\end{equation}
that splice together to long exact sequences of $\Gamma$-comodules
\begin{equation}
\label{eqANSS_LES}
0\to BP_*(X) \to BP_*(X\wedge I_0) \to BP_*(X\wedge I_1) \to \ldots,
\end{equation}
which is an $\ext$-acyclic $\Gamma$-resolution of $BP_*(X)$, since
\begin{equation}
\label{eqANSS_ExtBPXwedgeIn}
\ext^t(BP_*(X\wedge I_n))=
\begin{cases}
\pi_*(X\wedge I_n),& \textnormal{if }t=0,\\
0,& \textnormal{otherwise.}
\end{cases}
\end{equation}.


The Adams--Novikov spectral sequence (Adams~\cite{Adams64stableHomotopyTheory}, Novikov~\cite{Nov67methOfAlgTopCobordismTh}) for a spectrum $X$ is derived from the exact couple given by
\[
\xymatrix{
\pi_*(X\wedge X_n) \ar[r] & \pi_*(X\wedge \Sigma^{-n} I_n) \ar[r] & \pi_*(X\wedge \Sigma X_{n+1}) \ar@/^2pc/[ll]^{\deg =-1}
}
\]
and thus the $E_1$-page consists of terms
\[
E_1^{n,*}(X,BP_*):=\pi_*(X\wedge \Sigma^{-n} I_n)=\ext^0(BP_*(X\wedge \Sigma^{-n} I_n)).
\]
From the resolution~\eqref{eqANSS_LES} and~\eqref{eqANSS_ExtBPXwedgeIn} it follows that the $E_2$-page is given by
\[
E_2^{n,*}(X,BP_*) := \ext^*(BP_*(X)).
\]

Following Miller~\cite[Section 5]{Mil81adamsSS}, we smash \eqref{eqChromatricResolutionInFormOfAdamsResolution} with \eqref{eqAdamsResolutionOfX} in order to construct a ``double complex of spectra'':
First by a telescope argument we may assume that $X_0\supseteq X_1\supseteq\ldots$ and $N_0\supseteq N_1\supseteq\ldots$ are decreasing filtrations of $CW$ spectra, $X_n/X_{n+1}=\Sigma^{-n} I_n$, and $N_n/N_{n+1}=\Sigma^{-n} M_n$.
We extend these definitions to negative indices simply by $X_n:=X_0$ and $N_n:=N_0$ for $n\leq 0$.
Let $I=U(\ZZ^2)$ be the poset of upsets in $\ZZ^2$, that is, the set of all subset $p\subset\ZZ^2$ with the property: if $(a,b)\leq (c,d)\in\ZZ^2$ and $(a,b)\in p$ then $(c,d)\in p$.
In $I$ we have $-\infty=\emptyset$ and $\infty=\ZZ^2$.
We define a filtration of $Z:=X_0\wedge N_0\iso S^0$ by
\[
Z_p:= \bigcup_{(i,j)\in p} X_i \wedge N_j, \ \ \ p\in I,
\]
with $Z_{-\infty}=\pt$ and $Z_\infty=Z$.
From that define an exact couple system $E$ over $I$ by
\[
E^p_q:=\pi_*(Z_p/Z_q)\ \ \ \textnormal{and}\ \ \ D_p:=\pi_*(Z_p).
\]
\begin{lemma}
The associated spectral system has the following properties.
\begin{enumerate}
\item Let $(p,q)\in I_2$. If $p\wo q=\{(i,j)\}$ with $i,j\geq 0$ then
\begin{equation}
\label{eqASSEpqWithOnlyOneDifference}
E^p_q = \pi_*(\Sigma^{-i}I_i \wedge \Sigma^{-j}M_j) = \ext^0(BP_*(\Sigma^{-i}I_i \wedge \Sigma^{-j}M_j)). 
\end{equation}
\item Taking homology of \eqref{eqASSEpqWithOnlyOneDifference} in direction $d_2$ yields
\[
S^{pz}_{bq}= \ext^i(BP_*(\Sigma^{-j}M_j)) = \ext^i(M^j[-j]).
\]
whenever $b\wo p=\{(i-1,j)\}$, $p\wo q=\{(i,j)\}$, $q\wo z=\{(i+1,j)\}$.
\item Let $f_n:=\{(i,j)\st i+j\geq n\}\in I$.
Then
\begin{equation}
\label{eqChrANSSdiagonalE1}
E^{f_n}_{f_{n+1}} = \bigoplus_{\substack{i,j\geq 0,\\i+j=n}} \pi_*(\Sigma^{-i}I_i \wedge \Sigma^{-j}M_j)
\iso E_1^{n,*}(S^0,BP_*) = \ext^0(BP_*(\Sigma^{-n} I_n)). 
\end{equation}
\item Taking homology of \eqref{eqChrANSSdiagonalE1} in the only possible anti-diagonal direction yields
\[
S^{f_n,f_{n+2}}_{f_{n-1},f_{n+1}} \iso E_2^{n,*}(S^0,BP_*) = \ext^n(BP_*).
\]
\item $S^{\infty,-\infty}_{\infty,-\infty} = \pi_*(S^0)$.
\end{enumerate}
\end{lemma}

More generally we can smash $Z$ and all $Z_p$ with a spectrum $X$ and obtain a spectral system with limit $S^{\infty,-\infty}_{\infty,-\infty} = \pi_*(X)$.
Note that this limit does not imply anything about convergence in the usual sense (e.g.\ the sub spectral system given by the $f_n$'s coincides with the ordinary Adams--Novikov spectral sequence and in general it does not converge to $\pi_*(X)$, see Ravenel~\cite[Thm. 4.4.1.(b)]{Rav86complexCobordismAndStableHomotopy}).

\begin{proof}
3.) follows from the fact that one can replace $\Sigma^{-n}I_n$ by $Z_n/Z_{n+1}$ in~\eqref{eqAdamsResolutionOfX} and $X_n$ by $Z_n$ and obtains another (non-standard) $BP$-Adams resolution of $S^0$.
The natural inclusion from the standard resolution~\eqref{eqAdamsResolutionOfX} into the new one induces the claimed isomorphism.
\end{proof}

Analogously to the above construction, 
we can also unify Miller's generalizations of May's and Mahowald's spectral sequences~\cite{Mil81adamsSS} into a spectral system.
This follows closely his construction in~\cite[Section 5]{Mil81adamsSS}.

\section{Successive Eilenberg--Moore spectral sequences}
\label{secSuccEilenbergMoore}

Eilenberg and Moore~\cite{EilMoo66homologyAndFibrations1} constructed a second quadrant cohomology spectral sequence for the following setting.
Suppose we are given a pullback diagram of fibrations,
\[
\xymatrix{
X \ar[r] \ar[d]_{f_1'} & E_1 \ar[d]^{f_1} \\
E_2 \ar[r]^{f_2}      & B,
}
\]
$f_1'$ being the pullback of the fibration $f_1$ along $f_2$.
Let $H(\blank):=H^*(\blank;k)$ denote singular cohomology with coefficients in some field $k$.
In this section we assume that all spaces have the homotopy type of a $CW$-complex and that $H(B)$ and all $H(E_i)$ are finite dimensional in every dimension.
For any map $f:X\to Y$, $H(Y)$ acts on $H(X)$ 
via $f^*$.
Assume that $\pi_1(B_0)$ acts trivially on $H(E_0)$.
Then there is an Eilenberg--Moore spectral sequence,
\[
E_2^{p,q} = \tor_{H(B)}^{p,q}(H(E_1),H(E_2))\ \impl \ H(X).
\]
The index $p\leq 0$ in $\tor_{H(B)}^{p,q}$ is the resolution index and $q\geq 0$ is the grading index.
Smith~\cite{Smith70LecturesEMSS} and Hodgkin~\cite{Hod75eqKuenneth} constructed a geometric Eilenberg--Moore spectral sequence, see also Hatcher~\cite[Chapter 3]{Hatcher04SS}. (So far it is still open whether the geometric construction yields a spectral sequence that is isomorphic to the original one, which Eilenberg and Moore defined algebraically.)
We will follow and extend their approach to the following setting.

\subsection{A cube of pullbacks}
Suppose we are given $n$ fibrations $f_i:E_i\to B$, $1\leq 1\leq n$.
For any subset $I\subseteq [n]$, let $X_I$ be the pullback of the maps $\{f_i\st i\in I\}$.
Thus, $E_i=X_{\{i\}}$ and $B=X_\emptyset$.
Let $X:=X_{[n]}$.
For $I\subseteq J\subseteq [n]$, $X_J\to X_I$ is a fibration.
Assume that $\pi(B)=0$ (as usual, this assumption can be weakened). 

We will construct a spectral system over $D(\ZZ^n)$ with limit
\[
S^{\infty,-\infty}_{\infty,-\infty}=H(X)
\]
and second page
$S^{pz^*}_{b^*q}((p_1,\ldots,p_n);n) =$
\[
\tor^{p_n}_{HB}\Big(\ldots \tor^{p_2}_{HB}\big(\tor^{p_1}_{HB}(HB,HE_1),HE_2\big) \ldots,HE_n\Big).
\]
Here, $\tor^{p_1}_{HB}(HB,HE_1)$ is equal to $HE_1$ if $p_1=0$ and zero otherwise.
Actually there are $n!$ different second pages, one for every total order of the $n$ coordinates of $\ZZ^n$ and thus the order of the $E_i$. 

\subsection{The spectral system}
For a fixed topological space $B$, let $\topcat/B$ denote the category of spaces over $B$, whose objects are $p_X:X\to B$, and morphisms are the obvious commutative triangles.
Let $\topcat_B$ (Smith writes $(\topcat/B)_*$ instead) denote the category of ex-spaces over $B$, whose objects are spaces $X$ together with a map $p_X: X\to B$ and a section $s_X:B\to X$ such that $p_X\circ s_X=\id_B$, and morphisms are the obvious commutative double triangles.
We have a functor $\topcat/B\to \topcat_B$ that sends $X$ to $X_+:=X\sqcup B$.
Define $S^n_B:=B\times S^n\in \topcat_B$, where the section $B\to S^n_B$ is coming from picking a basepoint of~$S^n$.
Products $X\times_B Y$, smash products $X\topsmash_B Y$, suspensions $\Sigma_B X=S^n_B\topsmash_B X$, homotopies $\alpha\homotopic_B \beta$, mapping cylinders $M_B(\alpha)$, mapping cones $C_B(\alpha)$, and quotients $X/_B Y$ are all defined fiberwise; see Smith~\cite{Smith70LecturesEMSS}.

In order to simplify indices quite a bit, we move into a naive category of spectra over $B$, $\spectracat_B$, whose objects are sequences $(X_n)_{n\in\ZZ}$ of ex-spaces over $B$ together with maps $\Sigma_BX_n\to X_{n+1}$ over $B$.
There is a functor $\Sigma^\infty_B:\topcat_B\to \spectracat_B$ that sends $X$ to $X_n:=\Sigma^n_BX$ for $n\geq 0$, and $X_n:=B$ for $n<0$.
Let $S_B:=\Sigma^\infty_B(S^0_B)$.
This naive definition is all we need; for a thorough treatment of parametrized spectra see May--Sigurdsson~\cite[Section 2]{MaSi06paramHomotopyTheory}.
A reader who prefers spaces can proceed as in Smith~\cite{Smith70LecturesEMSS}.

We construct diagrams in $\spectracat_B$ for all $1\leq k\leq n$, similar to the Adams resolution,
\begin{equation}
\label{eqEMSSresolutionOfEk}
\xymatrix{
\Sigma^\infty_B(E_{k+})=:E_k^0 \ar[d] & E_k^{-1} \ar[d] \ar[l] & E_k^{-2} \ar[d] \ar[l] & \ \ldots \ar[l]    \\
I_k^0 & \Sigma^{-1}I_k^{-1} & \Sigma^{-2}I_k^{-2},
}
\end{equation}
where $\Sigma^{-i}I_k^{-i}$ is defined as 
the spectrum over $B$ with $(\Sigma^{-i}I_k^{-i})_n := ((E_k^{-i})_n/s(B))\times B$ and the natural maps (the projection to $B$ is the projection to the second factor),
and $\Sigma E_k^{-i-1}$ is defined as the cofiber of $E_k^{-i}\to \Sigma^{-i}I_k^{-i}$.

Let $H_B(\blank):\spectracat_B\to \grabcat$ denote cohomology relative to $\Sigma^\infty_B(B)$ with coefficients in $k$.
Since $H_B(\Sigma^{-i}I_k^{-i})=H_B(S_B)\otimes_k H_B(E_k^{-i})$, $H_B(\Sigma^{-i}I_k^{-i})$ is a free $H_B(S_B)$-module and $H_B(\Sigma^{-i}I_k^{-i})\to H_B(E_k^{-i})$ is surjective.
Thus the associated long exact sequences in cohomology become short exact sequences
\[
0\to H_B(\Sigma E_k^{-i-1})\to H_B(\Sigma^{-i}I_k^{-i})\to H_B(E_k^{-i})\to 0
\]
They splice together to a free $H_B(S_B)$-resolution of $H_B(\Sigma^\infty_B(E_{k+}))$,
\[
\ldots\ \to H_B(I_k^2)\to H_B(I_k^1)\to H_B(I_k^0)\to H_B(\Sigma^\infty_B(E_{k+})) \to 0.
\]

For $i>0$ we define $E_k^i:=E_k^0$ and correspondingly $I_k^i:=\Sigma^\infty_B(B)$.
As in Section~\ref{secAdamsAndChromaticSS}, following Miller~\cite[Section 5]{Mil81adamsSS}, we can first assume by a telescope argument that the horizontal maps in~\eqref{eqEMSSresolutionOfEk} are inclusions of spectra over $B$, and $I_k^{-i}$ the corresponding quotients.
Then we smash the $n$ diagrams~\eqref{eqEMSSresolutionOfEk}:
For $p\in I:=D(\ZZ^n)$ define
\[
Z_p:=\bigcup_{P\in p} E_1^{p_1}\topsmash_B\ldots\topsmash_B E_n^{p_n},
\]
which gives an $I$-filtration of $Z:=E_1^0\topsmash_B\ldots\topsmash_B E_n^0$.
Let $S$ be the spectral system associated to this filtered spectrum and the cohomology theory $H_B$, as in Example~\ref{exExactCoupleOfFilteredSpace}.
In the notation of Section~\ref{secSecondaryConnection}, we have
\[
S^{pz}_{bq}(P;1)=H_B(I^{p_1}_1\topsmash_B\ldots\topsmash_B I^{p_n}_n).
\]
Since all $H_B(I^{p_i}_i)$ are of finite type (that is, finite dimensional in every dimension) and projective over $H_B(S_B)$, it follows as in~\cite[Prop. 5.1]{Smith70LecturesEMSS} that
\[
S^{pz}_{bq}(P;0)=H_B(I^{p_1}_1)\otimes_{H_B(S_B)}\ldots\otimes_{H_B(S_B)} H_B(I^{p_n}_n).
\]
Since the functors $\blank\otimes_{H_B(S_B)} H_B(I_k^{p_k})$ are exact, taking homology with respect to the differential in direction $-e_1$ yields
\[
S^{pz^*}_{b^*q}(P;1)=H_B(E_1^0)\otimes_{H_B(S_B)}H_B(I^{p_2}_2)\otimes_{H_B(S_B)}\ldots\otimes_{H_B(S_B)}H_B(I^{p_n}_n)
\]
if $p_1=0$ and zero otherwise.
Continuing taking homology with respect to differentials in direction $-e_2$,\ldots,$-e_n$ yields eventually
\[
S^{pz^*}_{b^*q}(P;n)=\tor^{p_n}_{H_B(S_B)}\Big(\ldots \tor^{p_2}_{H_B(S_B)}\big(H_B(E_1^0),H_B(E_2^0)\big)\ldots,H_B(E_n^0)\Big)
\]
if $p_1=0$ and zero otherwise.

\paragraph{Acknowledgements.}

I want to thank
Michael Andrews,
Tony Bahri,
Mark Behrens,
Ofer Gabber,
Mark Goresky,
Jesper Grodal,
Bernhard Hanke,
Bob Mac\-Pher\-son,
John McCleary,
Haynes Miller,
and in particular Pierre Deligne for very useful discussions, as well as the referee for very useful remarks.

This work was supported by Deutsche Telekom Stiftung at Technische Universit\"at Berlin and Freie Universit\"at Berlin, by NSF Grant DMS-0635607 at Institute for Advanced Study, by an EPDI fellowship at Institut des Hautes \'Etudes Scientifiques,  Forschungsinstitut f\"ur Mathematik (ETH Z\"urich), and the Isaac Newton Institute for Mathematical Sciences, and by Simons Foundation grant {\#}550023 at Boston University (in chronological order).%

%
%
%


%
%

\small

\bibliographystyle{plain}


\bibliography{../mybib07}

\medskip

\small
\noindent
Boston University\\
matschke@bu.edu

\end{document}